\documentclass[10pt,a4paper]{article}
\usepackage{setspace}
\usepackage{titlesec}
\usepackage{url}
\usepackage[T1]{fontenc}
\usepackage[all]{xy}
\usepackage[pdftex]{graphicx}
\usepackage{palatino,amsmath,amsfonts,latexsym,amssymb,amsmath}
\usepackage{graphics}
\usepackage{subfigure}
\usepackage[english]{babel}
\usepackage[usenames,dvipsnames]{color}
\usepackage{accents}
\usepackage{mathtools}
\usepackage{float}
\usepackage{multicol}
\usepackage{comment}
\usepackage{stmaryrd}
\usepackage[font={footnotesize,it}]{caption}
\usepackage{upgreek}
\usepackage[pdftex]{hyperref}
\usepackage[Conny]{fncychap}
\usepackage{imakeidx}
\usepackage{appendix}
\usepackage{pdfpages}
\usepackage{geometry}
\usepackage{tikz}
\usepackage{xcolor}
\usepackage{textcomp}

\usepackage{todonotes}

\makeindex

\newtheorem{teo}{Theorem}
\newtheorem{ex}{Example}

\newtheorem{cor}[teo]{Corollary}
\newtheorem{lem}[teo]{Lemma}
\newtheorem{prop}[teo]{Proposition}
\newtheorem{defn}[teo]{Definition}
\newtheorem{obs}[teo]{Remark}

\newenvironment{dem}{\smallskip \noindent{\bf Proof}: }
{\hfill $\blacksquare$\hspace{0in}\medskip}

\newenvironment{dem*}{\smallskip \noindent{\bf Proof}: }
{\hfill $\square$\hspace{0in}\medskip}

\newenvironment{demtranslem}{\smallskip \noindent{\bf Proof of Lemma \ref{transv}}: }
{\hfill $\blacksquare$\hspace{0in}\medskip}

\newenvironment{demjordan}{\smallskip \noindent{\bf Proof of Proposition \ref{jordan}}: }
{\hfill $\blacksquare$\hspace{0in}\medskip}

\newenvironment{demmain}{\smallskip \noindent{\bf Proof of Theorem \ref{mainformal}}: }
{\hfill $\blacksquare$\hspace{0in}\medskip}

\newenvironment{dempropdegreedis}{\smallskip \noindent{\bf Proof of Proposition \ref{degreedisgivesvers}}: }
{\hfill $\blacksquare$\hspace{0in}\medskip}

\newenvironment{demgap}{\smallskip \noindent{\bf Proof of Proposition \ref{eigenvaluegaptozero}}: }
{\hfill $\blacksquare$\hspace{0in}\medskip}

\DeclareMathOperator{\Span}{span}
\DeclareMathOperator{\spec}{Spec}
\DeclareMathOperator{\id}{Id}
\DeclareMathOperator{\diag}{diag}

\title{\LARGE Chaotic behavior in diffusively coupled systems}

\author
{Eddie Nijholt$^{1\ast}$, Tiago Pereira$^{1,2}$, Fernando C. Queiroz$^{1}$,  and  Dmitry Turaev$^{2}$\\
\\
\normalsize{$^{1}$Instituto de Ci\^encias Matem\'aticas e Computa\c{c}\~ao, Universidade de S\~ao Paulo, S\~ao Carlos, Brazil}\\
\normalsize{$^{2}$Department of Mathematics, Imperial College London, London, UK}\\
\\
\normalsize{$^\ast$To whom correspondence should be addressed; E-mail: eddie.nijholt@gmail.com}
}

\begin{document}
\maketitle	
\begin{abstract}
We study emergent oscillatory behavior in networks of diffusively coupled nonlinear ordinary differential equations. Starting from a situation where each isolated node possesses a globally attracting equilibrium point, we give, for an arbitrary network configuration, general conditions for the existence of the diffusive coupling of a homogeneous strength which makes the network dynamics chaotic. The method is based on the theory of local bifurcations we develop for diffusively coupled networks. We, in particular, introduce the class of the so-called versatile network configurations and prove that the Taylor coefficients of the reduction to the center manifold for any versatile network can take any given value.	
\end{abstract}

\tableofcontents
		
\section{Introduction}
		
Coupled dynamical systems play a prominent role in biology \cite{IzhikevichBook}, chemistry \cite{Kuramoto}, physics and other fields of science \cite{RMP}. Understanding the emergent dynamics of such systems is a challenging problem depending starkly on the underlying interaction structure \cite{Nijholt,TiagoJEMS,Li,Ricard, keller2006uniqueness}.

In the early fifties, Turing thought of the emergent oscillatory behavior due to diffusive interaction as a model for morphogenesis \cite{turing1952}.
We note that weak coupling of globally stable individual systems cannot alter the stability of the homogeneous regime, this is the globally attracting state. At the same time, no matter what the individual dynamics are, the strong diffusive coupling by itself stabilizes the homogeneous regime. Therefore the idea that the intermediate strength diffusive coupling can create a non-trivial collective behavior is quite paradoxical. However, in the mid-seventies, Smale \cite{smale1976} proposed an example of diffusion-driven oscillations. He considered two 4th order diffusively coupled differential equations which by themselves have globally asymptotically stable equilibrium points. Once the diffusive interaction is strong enough, the coupled system exhibits oscillatory behavior. Smale posed a problem to find conditions under which diffusively coupled systems would oscillate. 

Tomberg and Yakubovich \cite{Tomberg} proposed a solution to this problem for diffusive interaction of two systems with scalar nonlinearity. For networks, Pogromsky, Glad and Nijmeijer \cite{Pogromsky} showed that diffusion-driven oscillations can result from an Andronov-Hopf bifurcation. Moreover, they presented conditions to ensure the emergence of oscillations for general graphs. While this provides a good picture of the instability leading to periodic oscillations, there is evidence that the diffusive coupling may also lead to chaotic oscillations.  Indeed, Kocarev and Janic \cite{Kocarev} provided numerical evidence that two isolated Chua circuits having globally stable fixed points may exhibit chaotic behavior when diffusively coupled.  Along the same lines, Perlikowski and co-authors \cite{perlikowski2010routes} investigated numerically the dynamics of rings of unidirectionally coupled Duffing oscillators. Starting from the situation where each oscillator has an exponentially stable equilibrium point, once the oscillators are coupled akin to diffusion the authors found a great variety of phenomena such as rotating waves, the birth of periodic dynamics, as well as chaotic dynamics. 

Drubi, Ibanez and Rodriguez \cite{Ibanez} studied two diffusively coupled Brusselators. Starting from a situation where the isolated systems have a globally stable fixed point, they proved that the unfolding of the diffusively coupled system can display a homoclinic loop with an invariant set of positive entropy.

In this paper, we provide general conditions for diffusively coupled identical systems to exhibit chaotic oscillations. We describe necessary and sufficient conditions (the so-called skewness condition) on the linearization matrix  at an exponentially stable equilibrium point of the isolated system such that for any network of such systems there exists a diffusive coupling matrix such that the network has a nilpotent singularity and thus a nontrivial center manifold. When the network structure satisfies an extra condition, which we call versatility, we show that Taylor coefficients of the vector field on the center manifold are in general position. This allows us to employ the theory of bifurcations of nilpotent singularities due to Arneodo, Coullet, Spiegel and Tresser \cite{arneodo1985asymptotic} and Ibanez and Rodr\'{i}guez \cite{ibanezconfig} and to show that when the isolated system is at least four-dimensional, invariant sets of positive entropy (i.e., chaos) emerge in such networks. 

The paper is organized as follows: In Section 2, we formulate the main theorem. We introduce basic concepts about graph theory, and introduce the notion of versatility. In Section 3, we present examples and constructions of $\rho$-versatile graphs and illustrate some of the other concepts appearing in the main theorem. In Section 4, we show the existence of a positive-definite coupling matrix $D$ that yields a nilpotent singularity in the network system. In Section 5, we discuss the stability of the center manifold. Finally, in Section 6, we prove that chaotic behavior emerges in the coupled system by investigating the dynamics on the center manifold.  
		
\subsection{The model}
	
We consider ordinary differential equations $\dot{x}=f(x)$ with  $f\in\mathcal{C}^{\infty}(U,\mathbb{R}^{n}),\ n\in \mathbb{N}$ for some open set $U\subset{\mathbb{R}^{n}}$. We assume that $f$ has an exponentially stable fixed point in $U$; with no loss of generality we put the origin of coordinates to this point. We study a network of such systems coupled together according to a given graph structure by means of a diffusive interaction. Namely, we consider the following equation:
\begin{equation}\label{eq1}
\dot{x}_{i}=f(x_{i})+\alpha\sum_{j=1}^{N}w_{ij}D(x_{j}-x_{i}),\quad i=1,\dots,N,
\end{equation}
where $\alpha>0$ is the coupling strength, ${W} = (w_{ij})$ is the adjacency matrix of the graph, thus, $w_{ij} = 1$ if nodes $i$ and $j$ are connected and zero otherwise. Moreover, $D$ is a positive-definite matrix (that is $x^TDx > 0$ for all non-zero vectors $x$). 

The homogeneous regime $x=0$ persists for every value of the coupling strength $\alpha$. It keeps its stability for small $\alpha$ and is, typically, stable at sufficiently large $\alpha$. However, at intermediate values of the coupling strength, the stability of the homogeneous regime can be lost. Our goal is to investigate the accompanying bifurcations. The difficulty is that the structure of system (\ref{eq1}) is quite rigid: all network nodes are the same (are described by the same function $f$) and the diffusion coupling $\alpha D$ is the same for any pair of nodes. Therefore, the genericity arguments, standard for the bifurcation theory, cannot be readily applied and must be re-examined.
	
\subsection{Informal statement of main results}
	
Our main goal is to give conditions for the emergence of non-periodic dynamics in system (\ref{eq1}). Denote by $A = Df(0)$ the linearization matrix $n\times n$ of the individual uncoupled system at zero. Recall that matrix $A$ is Hurwitz when all its eigenvalues have strictly negative real parts. Our main result  can be stated as follows.
    
\begin{quotation}
\textit{Suppose that for some orthogonal basis the Hurwitz matrix $A$ has $m$ positive entries on the diagonal. Then, there exists a positive-definite matrix $D$ such that the linearization of system \eqref{eq1} at the homogeneous equilibrium at zero has a zero eigenvalue of multiplicity at least $m$ for a certain value of the coupling parameter  $\alpha>0$. If the network satisfies a condition we call versatility, for an appropriate choice of the nonlinearity of $f$, the corresponding center manifold has dimension precisely $m$ and the Taylor coefficients of the restriction of the system on the center manifold can take on any prescribed value.} 
\end{quotation}

The last statement means that the bifurcations of the homogeneous state of a versatile network follow the same scenarios as general dynamical systems. Applying the results for triple instability 
\cite{Ibanez, ibanezconfig} we obtain the following result.
	
\begin{quotation}
\textit{For $n\geqslant 4$, for any generic 2-parameter family of nonlinearities $f$ and any versatile network graph,
one can find the positive-definite matrix $D$ such that the homogeneous state of the coupled system \eqref{eq1} 
has a triple instability at certain value of the coupling strength $\alpha$, leading to chaotic dynamics for 
a certain region of parameter values.}	
\end{quotation}
	
The condition on the Jacobian of the isolated dynamics can be understood in a geometric sense as follows. We write $\dot{x} = f(x) = Ax + \mathcal{O}(|x|^2)$. We claim that if a nonzero vector $x_0 \in \mathbb{R}^n$ exists for which $\langle x_0, Ax_0\rangle > 0$, then there are points arbitrarily close to the origin, whose forward orbit has its Euclidian distance to the origin increasing for some time, before coming closer to the (stable) origin again.  To see why, consider $\|x(t)\|^2 = \langle x(t), x(t) \rangle$, then it follows that
$
\frac{d}{dt}\|x\|^2 =  2 \langle Ax, x\rangle + \mathcal{O}(|x|^3)
$, so $\langle Ax, x\rangle >0$ implies the growth of this derivative.
		
{The property of versatility holds for graphs with heterogeneous degrees - the simplest example is a star network. In a sense, versatility means that the network is not very symmetric. Given a graph, one verifies whether the versatility property holds by evaluating the eigenvectors of the graph's Laplacian matrix, so it is an effectively verifiable property.}
	
It is possible that a similar theory can be developed for the Andronov-Hopf bifurcation in diffusively coupled networks (an analysis of diffusion-driven Andronov-Hopf bifurcations was undertaken in \cite{Pogromsky} but the question of genericity of the restriction of the network system to the central manifold was not addressed there).

We also point out that symmetry is often instrumental in explaining and predicting anomalous behavior in network dynamical systems \cite{antoneli2008hopf,Nijholt, quiver}.

The network of just two symmetrically coupled systems has the corresponding graph Laplacian that is not versatile, yet the emergence of chaos via the triple instability has been established in \cite{Ibanez} for the system of two diffusively coupled Brusselators. In general, we do not know when the genericity of the Taylor coefficients of the center-manifold reduced vector field would hold if the graph is not versatile or when graph symmetries would impose conditions on the dynamics that forbid the existence of limiting sets of positive entropy.

\section{Main results}
	
We start by introducing the basic concepts involved in the set-up of the problem.
	
\subsection{Graphs}
	
A graph $G$ is an ordered pair $(V,E)$, where $V$ is a non empty set of vertices and $E$ is a set of edges connecting the vertices. We assume both to be finite and the graph to be undirected. The order of the graph $G$ is $|V|=N$, its number of vertices, and the size is $|E|$, its number of edges. We will not consider graphs with self-loops. The degree of a vertex is the number of edges that are connected to it.
	
\begin{equation}\label{degree}
k_{i}=\sum_{j=1}^{N}w_{ij},
\end{equation}
for $i=1,\dots,N.$ We also define ${K}=\diag\{k_{1},\dots,k_{N}\}$ to be the diagonal matrix of vertex degrees.
	
We only consider undirected graphs $G$, meaning that a vertex $i$ is connected with a vertex $j$ if and only if it is  \textit{vice-versa}. Thus, the adjacency matrix ${W}$ is a symmetric matrix.
In this context there is another important matrix related to the graph $G,$ which is the well-known \textit{Laplacian} discrete matrix $L_{G}$. It is defined by:
\begin{equation}
\nonumber L_{G}={K}-{W},
\end{equation}
so that each entry $l_{ij}$ of $L_{G}$ can be written as
\begin{equation}\label{entry}
l_{ij}=\delta_{ij}k_{i}-w_{ij},\quad i,j=1,\dots,N,
\end{equation}
where $\delta_{ij}$ is \textit{Kronecker\textquotesingle s} delta. The matrix $L_{G}$ provides us with important information about connectivity and synchronization of the network. It follows from Gershgorin disk  theorem \cite{gerschgorin31} that $L_{G}$ is positive semi-definite and thus its eigenvalues can be ordered as 
$$
0=\lambda_{1}\leq\lambda_{2}\leq\cdots\leq\lambda_{N-1}\leq\lambda_{N},
$$ 
and let $\{v_{1},\dots,v_{N}\}$ be the corresponding eigenvectors. {We assume the network is connected. This implies that the eigenvalue $\lambda_1 = 0$ is simple}.	
	
We are interested in a well-behaved class of graphs $G$ whose structure induces a special property of the associated Laplacian matrix. This property will be the existence of an eigenvector where the sum of certain coordinate powers is non-vanishing, which corresponds to a simple eigenvalue of $L_{G}$. To this end, we define:
	
\begin{defn}[$\rho$-versatile graphs]\label{versatile}	
Let $G=(V,E)$ be a graph and $\rho \in \mathbb{N}$ a positive integer. We say that $G$ is $\rho$-versatile for the eigenvalue-eigenvector pair $(\lambda,v)$ with $\lambda>0$, if the  Laplacian matrix $L_{G}$ has a simple eigenvalue $\lambda$ with corresponding eigenvector $v=(\nu_{1},\dots,\nu_{N})$, satisfying
\begin{equation}\label{specialproperty}
\sum_{i=1}^{N}\nu^{\ell}_{i}\neq0,\quad\forall\ell=2,\dots,\rho+1.
\end{equation}
\end{defn}
Note that any eigenvector $v=(\nu_{1},\dots,\nu_{N})$ for a non-zero eigenvalue necessarily satisfies $\sum_{i=1}^{N}\nu_{i} = 0$.  This is because $\nu$ is orthogonal to the eigenvector $(1, \dots, 1)$ for the eigenvalue $0.$
	
\subsection{Parametrizations}
	
We show that a system of diffusively coupled stable systems can display a wide variety of dynamical behavior, including the onset of chaos. As the coupling strength $\alpha$ increases, a non-trivial center manifold can emerge with no general restrictions on the Taylor coefficients of the reduced dynamics.
	
Note that we may alternatively write Equation (\ref{eq1}) in terms of the Laplacian:	
\begin{equation}\label{eanalyseqlaplacian}
\dot{x}_{i}=f(x_{i})-\alpha\sum_{j=1}^{N}l_{ij}Dx_{j},\quad i=1,\dots,N.
\end{equation}
Let $X:=col(x_{1},\dots,x_{N})$ denote the vector formed by stacking $x_{i}$'s in a single column vector. In the same way we define $F(X):=col(f(x_{1}),\dots,f(x_{N}))$. We obtain the compact form for equations \eqref{eq1} and  \eqref{eanalyseqlaplacian} given by
\begin{equation}\label{eqcompact}
\dot{X}=F(X)-\alpha(L_G\otimes D)X,
\end{equation}
where $\otimes$ stands for the Kronecker product.   In order to analyze systems of the form \eqref{eq1}, we allow $f$ to depend on a parameter $\varepsilon$ taking values in some open neighborhood of the origin $\Omega \subseteq \mathbb{R}^d$. For simplicity, we assume the fixed point at the origin persists: 
\begin{equation}
f(0;\varepsilon) = 0 \qquad \forall \, \varepsilon \in \Omega.
\end{equation}
We assume the origin to be exponentially stable for $\varepsilon = 0$, from which stability follows for sufficiently small $\varepsilon$ as well. Note that the non-linear diagonal map $F$ now depends on the parameter $\varepsilon$ as well. 

We start with our working definition of center manifold reduction.
\begin{defn}
Let
\begin{align}\label{parametrsyst}
H: \mathbb{R}^n \times \Omega \rightarrow \mathbb{R}^n
\end{align}
be a family of vector fields on $\mathbb{R}^n$, parameterized by a variable $\varepsilon$ in an open neighborhood of the origin $\Omega \subseteq \mathbb{R}^d$.  Assume that $H(0;\varepsilon) = 0$ for all $\varepsilon \in \Omega$, and denote by $\mathcal{E}^c \subseteq \mathbb{R}^n$ the center subspace of the Jacobian $D_xH(0;0)$ in the direction of $\mathbb{R}^n$. A (local) \emph{parameterized center manifold} of the system \eqref{parametrsyst} is a (local) center manifold of the unparameterized system $\tilde{H}$ on $\mathbb{R}^n \times \Omega$, given by
\begin{align}\label{parametrsystaug} 
\tilde{H}(x; \varepsilon) = (H(x; \varepsilon), 0) \in \mathbb{R}^n \times \mathbb{R}^d\, ,
\end{align}
for $x \in \mathbb{R}^n$ and $\varepsilon \in \Omega$. We  say that the parameterized center manifold is of dimension $\dim(\mathcal{E}^c)$, and is parameterized by $d$ variables. Under the assumptions on $H$, the center subspace of $\tilde{H}$ at the origin is equal to $\mathcal{E}^c \times \mathbb{R}^d$. We can show that the dynamics on the center manifold of Equation \eqref{parametrsystaug} is conjugate to that of a locally defined system
\begin{align}\label{parametrsystred}
\tilde{R}(x_c;\varepsilon) = (R(x_c; \varepsilon), 0)\, ,
\end{align}
on $\mathcal{E}^c \times \Omega$, where the conjugation respects the constant-$\varepsilon$ fibers. The map $R$  satisfies $R(0; \varepsilon) = 0$ for all $\varepsilon$ for which this local expression is defined, and we have $D_{x_c}R(0;0) = D_xH(0;0)|_{\mathcal{E}^c}: \mathcal{E}^c \rightarrow \mathcal{E}^c$. We will refer to $R: \mathcal{E}^c \times \Omega \rightarrow \mathcal{E}^c$ as a \emph{parameterized reduced vector field} of $H.$ 
 \end{defn}
 	
In the definition above, the constant and linear terms of the parameterized reduced vector field $R$ are given. Motivated by this, we will write $H^{[2,\rho]}$ for any map $H$ to denote the non-constant, non-linear terms in the Taylor expansion around the origin of $H$, up to terms of order $\rho$. In other words, we have
\begin{equation}
\nonumber H(x)= H(0) + DH(0)x + H^{[2,\rho]}(x) + \mathcal{O}(||x||^{\rho+1}).
\end{equation}
Given vector spaces $W$ and $W'$, we will use $\mathcal{P}_2^l(W;W')$ to denote the linear space of polynomial maps from $W$ to $W'$ with terms of degree $2$ through $l.$ It follows that $H^{[2,l]} \in \mathcal{P}_2^l(W;W')$ for $H: W \rightarrow W'.$
 	 	
We are  interested in the situation where the domain of $H$ involves some parameter space $\Omega$, in which case $H^{[2,\rho]}$ involves all non-constant, non-linear terms up to order $\rho$ in both types of variables (parameter and phase space). For instance, if $H$ is a map from $\mathbb{R} \times \Omega$ to $\mathbb{R}$ with $\Omega \subseteq \mathbb{R}$, then $H^{[2,3]}(x; \varepsilon)$ involves the terms $$a_{1}x^2, a_{2}x\varepsilon, a_{3}\varepsilon^2, a_{4}x^3, a_{5}x^2\varepsilon, a_{6}x\varepsilon^2 \text{ and } a_{7}\varepsilon^3,$$ with some constants $a_{i}$. Note that  a condition on $H$ might put restraints on $H^{[2,\rho]}$ as well. For instance, if $H(0;\varepsilon) = 0$ for all $\varepsilon \in \Omega$, then $H^{[2,3]}(x; \varepsilon)$ does not involve the terms $\varepsilon^2$ and  $\varepsilon^3.$
 
\subsection{Main theorems}
 	
We  now formulate the main theorem, along with an important corollary. 
	
\begin{teo}[Main Theorem]\label{mainformal}
For any $\alpha \geq 0$, consider the $\varepsilon$-family of network dynamical systems given by 
\begin{equation}\label{systeminainthr}
\dot{X}=F(X;\varepsilon)-\alpha(L_{G}\otimes D)X.
\end{equation}
Denote by $A = D_xf(0;0)$ the Jacobian of the isolated dynamics. If there exist $m$ mutually orthogonal vectors $x_1,\dots,x_m$ such that $\langle x_i, Ax_i \rangle > 0$, then there exists a positive-definite matrix $D$ together with a number $\alpha^*>0$ such that the system of Equation \eqref{systeminainthr} has a local parameterized center manifold of dimension at least $m$ for $\alpha = \alpha^*$.
 
Suppose that the graph $G$  is $\rho$-versatile for the pair $(\lambda,v)$. After an arbitrarily small perturbation to $A$ if needed,  there exists a positive-definite matrix $D$ and a number  $\alpha^*>0$ such that the following holds:
\begin{enumerate}
\item The system of Equation \eqref{systeminainthr} has a local parameterized center manifold of dimension exactly $m$ for $\alpha = \alpha^*$.
\item Denote by $R: \mathcal{E}^c \times \Omega \rightarrow \mathcal{E}^c$ the corresponding parameterized reduced vector field, then $R(0; \varepsilon) = 0$ for all $\varepsilon \in \Omega$ and $D_xR(0;0):  \mathcal{E}^c  \rightarrow \mathcal{E}^c$ is nilpotent.
\item The higher order terms $R^{[2,\rho]}$ can take on any value in $\mathcal{P}_2^\rho(\mathcal{E}^c \times \Omega; \mathcal{E}^c)$ (subject to $R^{[2,\rho]}(0; \varepsilon) = 0$) as $f^{[2,\rho]}$ is varied (subject to $f^{[2,\rho]}(0; \varepsilon) = 0$).
\end{enumerate}	
\end{teo}

The above result guarantees the existence of the center manifold and the reduced vector field. When the dimension of the isolated dynamics is at least 4, the reduced vector field can exhibit invariant sets of positive entropy as the following result shows.
	
\begin{cor}[Chaos]\label{chaos}
Assume the conditions of Theorem \ref{mainformal} to hold for $m=3$ and $\rho = 2$. Then, in a generic $3$-parameter system we have the emergence of chaos through the formation of a Shilnikov loop on the center manifold. In particular, chaos can form this way in a system of $4$-dimensional nodes coupled diffusively in a network.
\end{cor}

In the Appendix, we show that the conditions of Theorem \ref{mainformal} are natural, by constructing multiple classes of networks that are $\rho$-versatile for any $\rho \in \mathbb{N}$, as well as by giving examples of matrices $A$ that satisfy the conditions of the theorem. In Subsection \ref{means_of_the_complement_graph}, we present a geometric way of constructing $\rho$-versatile graphs, by means of the so-called complement graph. In Subsection \ref{stargraphs}, we then show using direct estimates that star graphs satisfy $\rho$-versatility. Finally, in Subsection \ref{InternalDynamics}, we present examples of matrices that satisfy the conditions of Theorem \ref{mainformal}. In particular, we will see that being Hurwitz is no obstruction.
        	
\section{Proof of main theorem}
	
In this section we present the proof of Theorem \ref{mainformal}. We start by analyzing the linearized system in Subsection \ref{linearization}, after which we perform center manifold reduction and have a detailed look at the reduced vector field in Subsection \ref{Centermanifoldreduction}.
	
\subsection{Linearization}\label{linearization}
	
In this subsection we investigate the linear part of the system 
\begin{equation}\label{systeminchapter3}
\dot{X}=F(X;\varepsilon)-\alpha(L_G\otimes D)X,
\end{equation}
from Theorem \ref{mainformal}. Writing $A \in \mathbb{R}^{n \times n}$ for the Jacobian of $f$ at the origin, we see that the linearization of Equation \eqref{systeminchapter3} at the origin is given by
\begin{equation}\label{systeminchapter3lin}
\dot{Y}=(\id_{N} \otimes A-\alpha L_G\otimes D)Y.
\end{equation}

An important observation is the following: if $v$ is an eigenvector of $L_{G}$ with eigenvector $\lambda \in \mathbb{R}$, then the linearization above sends a vector $v \otimes x$ with $x \in \mathbb{R}^n$ to
\begin{align}\label{calcinvariance}
(\id_N \otimes A-\alpha L_G\otimes D)(v \otimes x) &= v \otimes Ax  - \alpha L_Gv \otimes Dx \\ \nonumber
&= v \otimes (A-\alpha \lambda D)x \, .
\end{align}

It follows that the space
\begin{equation}
v \otimes \mathbb{R}^n := \{v \otimes x \mid x \in \mathbb{R}^n\} \subseteq \mathbb{R}^N \otimes  \mathbb{R}^n
\end{equation}
is kept invariant by the linear map of Equation \eqref{systeminchapter3lin}. We claim that $v \otimes \mathbb{R}^n$ is in fact a linear subspace. 

Equation \eqref{calcinvariance}  tells us that the linearization \eqref{systeminchapter3lin} restricted to $v \otimes \mathbb{R}^n$ is conjugate to $A-\alpha \lambda D :\mathbb{R}^n \rightarrow \mathbb{R}^n$.

Recall that eigenvectors $v_p$ and $v_q$ of $L_G$ corresponding to distinct eigenvalues can be chosen to be orthonormal. Finally, we write 
\begin{equation}
V_p := v_p \otimes \mathbb{R}^n
\end{equation}
for the corresponding linear subspaces of $\mathbb{R}^N \otimes \mathbb{R}^n$. We thus have a direct sum decomposition
\begin{equation}
\mathbb{R}^N \otimes \mathbb{R}^n = \bigoplus_{p=1}^N V_p ,
\end{equation}
where each component is respected by the linearization \eqref{systeminchapter3lin}, with the restriction to $V_p$ conjugate to $A-\alpha\lambda_p D :\mathbb{R}^n \rightarrow \mathbb{R}^n$. We see that the spectrum of the linearization \eqref{systeminchapter3lin} is given by the union of the spectra of the maps $A-\alpha\lambda_p D$, with a straightforward relation between the respective algebraic and geometric multiplicities.
	
This observation motivates the main result of this subsection,  Proposition \ref{jordan} below. 

In what follows, we denote by $\mathbb{M}_{n}(\mathbb{R})$ the space of $n$ by $n$ matrices over the field $\mathbb{R}$. We write  $\langle x,y\rangle:=x^{T}y$ for the Euclidean inner product between vectors $x,y \in \mathbb{R}^n$.
	
We state a technical lemma. Suppose we are given a block matrix 
\begin{equation}
M=\left(\begin{array}{cc}
M_{11} &  M_{12}\\
M_{21}& M_{22} 
\end{array}\right)
\end{equation} 
with blocks $M_{11}, \dots, M_{22}$. If $M_{22}$ is invertible then we may form the \emph{Schur complement} of $M$, given by
\begin{equation}
M/{M_{22}}:=M_{11}-M_{12}M_{22}^{-1}M_{21}.
\end{equation}
This expression has various useful properties. We are interested in the situation where $M$ is symmetric, so that $M_{12}^T = M_{21}$, $M_{11}^T = M_{11}$ and $M_{22}^T = M_{22}$. In that case the matrix $M$ is positive-definite if and only if both $M_{22}$ and $M/{M_{22}}$ are positive-definite.  Using this result, we may prove:
	
\begin{lem}\label{Schurlemma}
Let
\begin{equation}
D=\left(\begin{array}{cc}
A_{11}&A_{12}\\
A_{21}&A_{22} 
\end{array}\right)
\end{equation}
be a block matrix and assume $A_{22}=c \id$ for some scalar $c \in \mathbb{R}_{>0}$. Suppose furthermore that $A_{11}$ is positive-definite. Then, for $c$ sufficiently large the matrix $D$ is positive-definite as well.
\end{lem}
\begin{dem}
Clearly $D$ is positive-definite if and only if the symmetric matrix 
\begin{equation}
H := D+D^T =\left(\begin{array}{cc}
A_{11}+A_{11}^T &A_{12}+A_{21}^T\\
A_{21}+A_{12}^T&2c \id 
\end{array}\right) = \left(\begin{array}{cc}
H_{11} &H_{12}\\
H_{12}^T&H_{22}
\end{array}\right)
\end{equation}
is positive-definite. We set $H_{22} := 2c \id$, $H_{12} := A_{12}+A_{21}^T$ and $H_{11} := A_{11}+A_{11}^T$, the third of which is positive-definite as $A_{11}$ is. As $H_{22} = 2c \id$ is invertible with inverse $1/(2c)\id$, we may form the Schur complement 
\begin{eqnarray}
H/H_{22}&=&H_{11}-H_{12}H_{22}^{-1}H_{12}^T\\ \nonumber
&=&H_{11}-\frac{H_{12}H_{12}^T}{2c}.
\end{eqnarray}

As $H_{22}$ is positive-definite, it follows from the above discussion that $H$ is positive-definite if and only $H/H_{22}$ is. However, as $c \rightarrow \infty$  we have $H/H_{22} \rightarrow  H_{11}$, so $H/H_{22}$ behaves like a small perturbation of $H_{11},$ then it is positive-definite for $c>0$ large enough. This shows that $D$ is likewise positive-definite for large enough $c$. 
\end{dem}

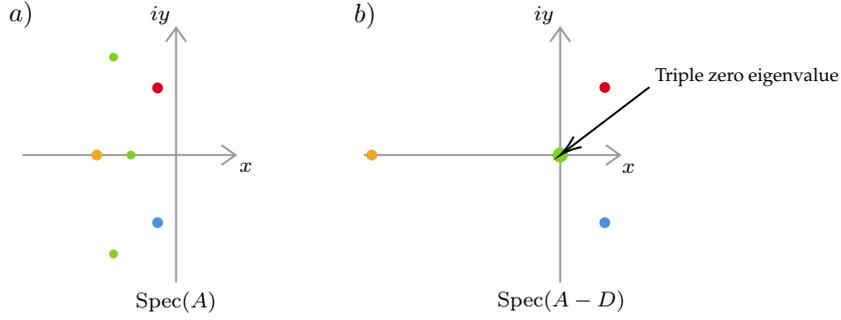
\begin{figure}[H]
\centering	
\tikzset{every picture/.style={line width=0.75pt}} 
\begin{tikzpicture}[x=0.75pt,y=0.75pt,yscale=-1,xscale=1]
\draw [color={rgb, 255:red, 155; green, 155; blue, 155 }  ,draw opacity=1 ] (183.05,97) -- (310.88,97)(281.05,33) -- (281.05,161) (303.88,92) -- (310.88,97) -- (303.88,102) (276.05,40) -- (281.05,33) -- (286.05,40)  ;
\draw  [color={rgb, 255:red, 126; green, 211; blue, 33 }  ,draw opacity=1 ][fill={rgb, 255:red, 126; green, 211; blue, 33 }  ,fill opacity=1 ] (277.8,97) .. controls (277.8,95.21) and (279.26,93.75) .. (281.05,93.75) .. controls (282.84,93.75) and (284.3,95.21) .. (284.3,97) .. controls (284.3,98.79) and (282.84,100.25) .. (281.05,100.25) .. controls (279.26,100.25) and (277.8,98.79) .. (277.8,97) -- cycle ;
\draw    (325.6,62.85) -- (282.64,95.78) ;
\draw [shift={(281.05,97)}, rotate = 322.53] [color={rgb, 255:red, 0; green, 0; blue, 0 }  ][line width=0.75]    (10.93,-3.29) .. controls (6.95,-1.4) and (3.31,-0.3) .. (0,0) .. controls (3.31,0.3) and (6.95,1.4) .. (10.93,3.29)   ;
\draw  [color={rgb, 255:red, 245; green, 166; blue, 35 }  ,draw opacity=1 ][fill={rgb, 255:red, 245; green, 166; blue, 35 }  ,fill opacity=1 ] (184.78,97) .. controls (184.78,95.78) and (185.77,94.79) .. (186.99,94.79) .. controls (188.21,94.79) and (189.2,95.78) .. (189.2,97) .. controls (189.2,98.22) and (188.21,99.21) .. (186.99,99.21) .. controls (185.77,99.21) and (184.78,98.22) .. (184.78,97) -- cycle ;
\draw  [color={rgb, 255:red, 208; green, 2; blue, 27 }  ,draw opacity=1 ][fill={rgb, 255:red, 208; green, 2; blue, 27 }  ,fill opacity=1 ] (301.12,63) .. controls (301.12,61.78) and (302.11,60.79) .. (303.33,60.79) .. controls (304.54,60.79) and (305.53,61.78) .. (305.53,63) .. controls (305.53,64.22) and (304.54,65.21) .. (303.33,65.21) .. controls (302.11,65.21) and (301.12,64.22) .. (301.12,63) -- cycle ;
\draw  [color={rgb, 255:red, 74; green, 144; blue, 226 }  ,draw opacity=1 ][fill={rgb, 255:red, 74; green, 144; blue, 226 }  ,fill opacity=1 ] (301.12,131) .. controls (301.12,129.78) and (302.11,128.79) .. (303.33,128.79) .. controls (304.54,128.79) and (305.53,129.78) .. (305.53,131) .. controls (305.53,132.22) and (304.54,133.21) .. (303.33,133.21) .. controls (302.11,133.21) and (301.12,132.22) .. (301.12,131) -- cycle ;
\draw [color={rgb, 255:red, 155; green, 155; blue, 155 }  ,draw opacity=1 ] (12.6,97) -- (119,97)(89.32,33) -- (89.32,161) (112,92) -- (119,97) -- (112,102) (84.32,40) -- (89.32,33) -- (94.32,40)  ;
\draw  [draw opacity=0][fill={rgb, 255:red, 126; green, 211; blue, 33 }  ,fill opacity=1 ] (55.89,47.75) .. controls (55.89,46.53) and (56.88,45.54) .. (58.1,45.54) .. controls (59.32,45.54) and (60.31,46.53) .. (60.31,47.75) .. controls (60.31,48.97) and (59.32,49.96) .. (58.1,49.96) .. controls (56.88,49.96) and (55.89,48.97) .. (55.89,47.75) -- cycle ;
\draw  [draw opacity=0][fill={rgb, 255:red, 126; green, 211; blue, 33 }  ,fill opacity=1 ] (55.89,146.75) .. controls (55.89,145.53) and (56.88,144.54) .. (58.1,144.54) .. controls (59.32,144.54) and (60.31,145.53) .. (60.31,146.75) .. controls (60.31,147.97) and (59.32,148.96) .. (58.1,148.96) .. controls (56.88,148.96) and (55.89,147.97) .. (55.89,146.75) -- cycle ;
\draw  [color={rgb, 255:red, 208; green, 2; blue, 27 }  ,draw opacity=1 ][fill={rgb, 255:red, 208; green, 2; blue, 27 }  ,fill opacity=1 ] (77.89,63.25) .. controls (77.89,62.03) and (78.88,61.04) .. (80.1,61.04) .. controls (81.32,61.04) and (82.31,62.03) .. (82.31,63.25) .. controls (82.31,64.47) and (81.32,65.46) .. (80.1,65.46) .. controls (78.88,65.46) and (77.89,64.47) .. (77.89,63.25) -- cycle ;
\draw  [color={rgb, 255:red, 74; green, 144; blue, 226 }  ,draw opacity=1 ][fill={rgb, 255:red, 74; green, 144; blue, 226 }  ,fill opacity=1 ] (77.89,131) .. controls (77.89,129.78) and (78.88,128.79) .. (80.1,128.79) .. controls (81.32,128.79) and (82.31,129.78) .. (82.31,131) .. controls (82.31,132.22) and (81.32,133.21) .. (80.1,133.21) .. controls (78.88,133.21) and (77.89,132.22) .. (77.89,131) -- cycle ;
\draw  [draw opacity=0][fill={rgb, 255:red, 126; green, 211; blue, 33 }  ,fill opacity=1 ] (64.56,97) .. controls (64.56,95.78) and (65.55,94.79) .. (66.77,94.79) .. controls (67.99,94.79) and (68.97,95.78) .. (68.97,97) .. controls (68.97,98.22) and (67.99,99.21) .. (66.77,99.21) .. controls (65.55,99.21) and (64.56,98.22) .. (64.56,97) -- cycle ;
\draw  [color={rgb, 255:red, 245; green, 166; blue, 35 }  ,draw opacity=1 ][fill={rgb, 255:red, 245; green, 166; blue, 35 }  ,fill opacity=1 ] (47.56,97) .. controls (47.56,95.78) and (48.55,94.79) .. (49.77,94.79) .. controls (50.99,94.79) and (51.97,95.78) .. (51.97,97) .. controls (51.97,98.22) and (50.99,99.21) .. (49.77,99.21) .. controls (48.55,99.21) and (47.56,98.22) .. (47.56,97) -- cycle ;

\draw (281.6,163.8) node [anchor=north] [inner sep=0.75pt]  [font=\footnotesize]  {$\spec( A-D)$};
\draw (327.1,51.75) node [anchor=north west][inner sep=0.75pt]  [font=\scriptsize] [align=left] {Triple zero eigenvalue};
\draw (310.5,99.65) node [anchor=north west][inner sep=0.75pt]  [font=\footnotesize]  {$x$};
\draw (266.5,20.65) node [anchor=north west][inner sep=0.75pt]  [font=\footnotesize]  {$iy$};
\draw (176.78,18.65) node [anchor=north west][inner sep=0.75pt]    {$b)$};
\draw (89.32,163.95) node [anchor=north] [inner sep=0.75pt]  [font=\footnotesize]  {$\spec( A)$};
\draw (119.5,99.4) node [anchor=north west][inner sep=0.75pt]  [font=\footnotesize]  {$x$};
\draw (75,20.4) node [anchor=north west][inner sep=0.75pt]  [font=\footnotesize]  {$iy$};
\draw (4.5,18.4) node [anchor=north west][inner sep=0.75pt]    {$a)$};
\end{tikzpicture}
\caption{\label{fig8}An illustration of Proposition \ref{jordan}. Figure $a)$ shows the eigenvalues of a particular matrix $A$, which might be the linear part of the isolated dynamics $f$ at the origin. The matrix $A$ has $6$ eigenvalues all strictly on the left half of the complex plane. The existence of $m=3$ mutually orthogonal vectors satisfying $\eqref{posit}$ ensures that a positive-definite matrix $D$ exists such that $A-D$ has a $3$-dimensional generalized kernel. In other words, the subtraction of $D$ has moved $3$ eigenvalues to the origin, see Figure $b).$}
\end{figure}

\begin{prop}\label{jordan}
Let $A\in\mathbb{M}_{n}(\mathbb{R})$ be a matrix.  There exist $m$ mutually orthogonal vectors $x_{1},\dots,x_{m}$ such that
\begin{equation}\label{posit}
\langle x_{i},Ax_{i} \rangle>0 \text{ for all } i=1,\dots,m,
\end{equation} 
{if and only if} there exists a positive-definite matrix $D$ such that $A-D$ has at least $m$ zero eigenvalues, counted with algebraic multiplicity.
\end{prop}
Figure \ref{fig8} shows an illustration of Proposition \ref{jordan}.

\begin{obs}\label{number_m_skewness}
		Note that any Hurwitz matrix $A$ has a negative trace, as this number equals the sum of its eigenvalues. It follows that Equation~\eqref{posit} can then only hold when $m<n$, where $n$ is the size of $A.$ In this case, if our goal is to find a $3$-dimensional center manifold for the network, we need $3$ zero eigenvalues and so we must have at least $n=4$. 
\end{obs}

\begin{obs}
		The number of zero eigenvalues for $A-D$ is directly connected with the number of mutually orthogonal vectors for Equation~\eqref{posit}. Moreover, the positive-definite matrix $D$ that we construct is in general not unique. Hence, if we have $m$ such orthogonal vectors, then for each $k \leq m$ we might construct a different matrix $D$ such that $A-D$ has $k$ zero eigenvalues. 
\end{obs}
	
\begin{demjordan}
Suppose first that we have $m$ mutually orthogonal vectors $x_{i}$ such that for each $i=1,\dots,m$ we have
\begin{equation}\label{posit2}
\langle x_{i},Ax_{i} \rangle=x_{i}^{T}Ax_{i}>0.
\end{equation}
Note that we then also have 
\begin{equation}
\langle x_{i},(A+A^{T})x_{i} \rangle = 2\langle x_{i},Ax_{i} \rangle> 0
\end{equation} 
for all $i$. We may re-scale the $x_i$ by any non-zero factor, so that we will now assume without loss of generality that $\|x_i \| = 1$ for all $i$. We start by constructing an auxiliary upper-diagonal $(m\times m)$-matrix $P$ as follows:
\begin{equation}
\nonumber P=P_{A}=\left(\begin{array}{cccc}
0&p_{1,2}&\cdots&p_{1,m} \\
0&0&\cdots&p_{2,m}\\
\vdots&\vdots&\ddots&\vdots\\
0&0&\cdots&0
\end{array}\right)_{m\times m}\, ,
\end{equation}
where each entry $p_{i,j}$ is defined by the rule:
\begin{eqnarray}
\nonumber p_{i,j}&=&x_{i}^{T}(A+A^{T})x_{j},\quad\text{ for all $i<j$;}\\
\nonumber p_{i,j}&=&0,\quad\text{ for all $i\geq j$.}
\end{eqnarray}
We construct $D$ by first defining it on the mutually orthogonal vectors $x_{1},\dots,x_{m}$ as:
\begin{eqnarray}\label{iterativPDA}
Dx_{1}&=&Ax_{1}\\
\nonumber Dx_{2}&=&Ax_{2}-p_{1,2}x_{1}\\
\nonumber &\vdots&\\
\nonumber Dx_{m}&=&Ax_{m}-p_{1,m}x_{1}-\cdots-p_{m-1,m}x_{m-1}.
\end{eqnarray}
Note that $(A-D)x_1 = 0$, whereas $(A-D)x_2 \in \Span(x_1)$, $(A-D)x_3 \in \Span(x_1, x_2)$ and so forth. This shows that the restriction of $A-D$ to $\Span(x_1, \dots, x_m)$ is nilpotent. Equation \eqref{iterativPDA} can be rewritten as 
\begin{equation}\label{formulaforD}
(A-D)_{n\times n}{X}_{n\times m}={X}_{n\times m}P_{m\times m}\, ,
\end{equation}
with $X=(x_{1} \cdots x_{m})$ the $(n\times m)$-matrix with  columns given by the vectors $x_1, \dots, x_m$.
		
To complete our construction of $D$, we let $y_{m+1},\dots, y_{n} \in \mathbb{R}^n$ be mutually orthogonal vectors of norm $1$ such that $y_{k}\perp x_{i}$ for all $i=1,\dots,m$ and $k=m+1,\dots, n$. We define $D$ on $\Span(y_{m+1}, \dots, y_n )$ by simply setting $Dy_{k}=c y_{k}$ for all $k$ and some constant $c >0$ that will be determined later.
	
To show that $c$ can be chosen such that $D$ is positive-definite, let $z \in \mathbb{R}^n$ be any non-zero vector and write
\begin{equation}
\nonumber z={X}a+{Y}b
\end{equation}
where ${Y}=(y_{m+1}\cdots y_{n})$ is the $(n\times (m-n))$-matrix with columns the vectors $y_{k}$, and where $a \in \mathbb{R}^m , b \in \mathbb{R}^{n-m}$ express the components of $z$ with respect to the basis $\{x_1, \dots, x_m, y_{m+1}, \dots, y_n\}$. Note that we have 
\begin{align}\label{identitiesinproofDA}
&{X}^T{X} = I_{m \times m},  &&{Y}^T{Y} = I_{(n-m) \times (n- m)}, &&  \\ \nonumber
&{X}^T{Y} = 0_{m \times (n-m)},  &&{Y}^T{X} = 0_{(n-m) \times m},  &&D{Y} = c{Y}, 
\end{align}
by construction. We calculate

\begin{eqnarray}
z^{T}Dz&=&({X}a+{Y}b)^{T}D({X}a+{Y}b)\\
\nonumber &=&a^{T}{X}^{T}D{X}a+a^{T}{X}^{T}D{Y}b+b^{T}{Y}^{T}D{X}a+b^{T}{Y}^{T}D{Y}b\\
\nonumber &=&a^{T}{X}^{T}(A{X}-{X}P)a+a^{T}{X}^{T}D{Y}b+b^{T}{Y}^{T}(A{X}-{X}P)a+b^{T}{Y}^{T}D{Y}b\\ \nonumber
&=&a^{T}({X}^{T}A{X}-{X}^{T}{X}P)a+a^{T}{X}^{T}c{Y}b+b^{T}({Y}^{T}A{X}-{Y}^{T}{X}P)a+b^{T}{Y}^{T}c{Y}b\\ \nonumber
&=&a^{T}({X}^{T}A{X}-P)a+b^{T}{Y}^{T}A{X}a+cb^{T}b,\label{eqposit}
\end{eqnarray}

where in the third step we have used Equation \eqref{formulaforD}, and where me make use of the identities in \eqref{identitiesinproofDA}. We see that $D$ is positive-definite if the same holds for the matrix 
\begin{equation}
\tilde{D}=\left(\begin{array}{cc}
{X}^{T}A{X}-P & 0\\
{Y}^{T}A{X}& c\id 
\end{array}\right)\, .
\end{equation}
Next, we claim that the $(m\times m)$-matrix ${X}^{T}A{X}-P$ is positive-definite. Indeed, by definition of $P$ we have
\begin{eqnarray}
\nonumber ({X}^{T}A{X}-P)+({X}^{T}A{X}-P)^{T}&=&{X}^{T}(A+A^{T}){X}-(P+P^{T})\\
\nonumber  &=&\diag(2x_{1}^{T}Ax_{1},\dots,2x_{m}^{T}Ax_{m}),
\end{eqnarray}
which is a diagonal matrix and positive-definite by the hypothesis (\ref{posit}). We may thus apply Lemma \ref{Schurlemma} to $\tilde{D}$, so that for $c>0$ sufficiently large $\tilde{D}$ and $D$ are indeed positive-definite.
		
Conversely, suppose there exists a positive-definite matrix $D$ such that
\begin{equation}
\nonumber A-D\quad\text{has $m$ zero eigenvalues. }
\end{equation}

We will prove that $m$ mutually orthogonal vectors $x_1, \dots, x_m$ exist satisfying
\begin{equation}
\langle x_{i},Ax_{i}\rangle>0\text{  for all $i=1,\dots,m.$ }
\end{equation}
By assumption, we may choose $m$ linearly independent vectors $y_{1},\dots,y_{m}$ such that
\begin{equation}
\nonumber  (A-D)y_{1}=0\quad\text{ and }\quad (A-D)y_{i}=\iota_i y_{i-1}\text{ for $i=2,\dots,m$},
\end{equation}
where $\iota_i \in \{0,1\}$ for all $i>1$. Next, we apply the  Gram-Schmidt orthonormalization process to the vectors $y_i$. That is, we set
\begin{eqnarray}\label{gramschmidt}
\nonumber x_{1}&=&y_{1}\\
\nonumber x_{2}&=&y_{2} - \alpha_{2,1} \cdot x_1\\
\vdots\\
\nonumber x_{i}&=&y_{i}-\alpha_{i,1}\cdot x_{1}-\alpha_{i,2}\cdot x_{2}-\cdots-\alpha_{i,i-1}\cdot x_{i-1}\, ,
\end{eqnarray}
where each coefficient is given by
\begin{equation}
\nonumber \alpha_{i,j}=\frac{\langle y_{i},x_{j} \rangle}{\langle x_{j},x_{j} \rangle} \quad \text{ for } j<i. 
\end{equation}
It follows that $\langle x_{i},x_{j} \rangle = 0$ whenever $i\neq j$. Moreover, we see from Equation \eqref{gramschmidt} that we may write
\begin{equation}\label{betas}
x_{i} = y_{i} + \sum_{j<i} \beta_{i,j} y_j\quad \text{ and thus } \quad y_{i} = x_{i} + \sum_{j<i} \beta'_{i,j} x_j
\end{equation}
for some coefficients $\beta_{i,j}, \beta'_{i,j} \in \mathbb{R}$. We therefore have $(A-D)x_1 = 0$, and for $2 \leq i \leq  m$ we find
\begin{align}
(A-D)x_i &=  (A-D)\left(y_{i} + \sum_{j<i} \beta_{i,j} y_j\right) = (A-D)y_{i} + \sum_{j<i} \beta_{i,j} (A-D)y_j \\ &= \nonumber \iota_i y_{i-1} + \sum_{1<j<i} \beta_{i,j} \iota_j y_{j-1} = \sum_{j<i}\gamma_{i,j}x_j\, ,
\end{align}
for certain $\gamma_{i,j} \in \mathbb{R}$. By orthogonality of the $x_i$ we get
\begin{eqnarray}
\nonumber x^{T}_{1}(A-D)x_{1}&=&0 \quad \text{ and }\\
\nonumber x^{T}_{i}(A-D)x_{i}&=&\sum_{j<i}\gamma_{i,j}x_i^Tx_j = 0
\end{eqnarray}
for all $i=2,\dots,m.$ Finally, it follows that
\begin{equation}
\nonumber x^{T}_{i}Ax_{i}=x^{T}_{i}Dx_{i}>0\quad \text{for all $i=1,\dots,m,$}
\end{equation}
which completes the proof.
\end{demjordan}
	
\begin{cor}\label{obsaboutremainingeig}
	The proof of Proposition \ref{jordan} tells us that the remaining eigenvalues of $A-D$ may be assumed to have (large) negative real parts.
\end{cor} 
\begin{dem}
If $m$ mutually orthogonal vectors $x_1, \dots, x_m$ exist such that
\begin{equation}\label{eqposit2}
\langle x_{i},Ax_{i}\rangle>0\text{  for all $i=1,\dots,m,$ }
\end{equation}
then a positive-definite matrix $D$ is constructed such that the restriction of $A-D$ to $\Span(x_1, \dots, x_m)$ is nilpotent. In particular, $A-D$ maps the space $\Span(x_1, \dots, x_m)$ into itself. It follows that the remaining eigenvalues of $A-D$ are given by those of the `other' diagonal block $P_U(A-D)|_U: U \rightarrow U$, where $U$ is some complement to $\Span(x_1, \dots, x_m)$ and $P_U$ is the projection onto $U$ along $\Span(x_1, \dots, x_m)$. If we choose $U = \Span(y_{m+1}, \dots, y_n)$ as in the proof of Proposition \ref{jordan}, then we see that $P_U(A-D)|_U= P_UA|_U - c\id_Y$. Choosing $c>0$ large enough then ensures that the remaining eigenvalues of $A-D$ are stable.

\end{dem}
	
\begin{cor}\label{obsaboutremainingeig2sdf}
	It follows from the proof of Proposition \ref{jordan} that $A-D$ can generically be assumed to have a one-dimensional kernel. In other words, whereas the algebraic multiplicity of the eigenvalue $0$ is $m$, its geometric multiplicity is generically equal to $1$.
\end{cor}	
\begin{dem}
To see why, assume $c>0$ is large enough so that $A-D$ has a generalized kernel of dimension precisely $m$, see Corollary \ref{obsaboutremainingeig}. From the proof of Proposition \ref{jordan} we see that the restriction of $A-D$ to its generalized kernel is conjugate to $P$. It follows that the dimension of the kernel of $A-D$ is equal to $1$ if 
\begin{equation}
p_{i,i+1}=x_{i}^{T}(A+A^{T})x_{i+1} \not= 0\quad\text{ for all $i \in \{1, \dots, m-1\}$}\, .
\end{equation}
This may be assumed to hold after a perturbation of the form
\begin{equation}
A \mapsto A+\varepsilon_1x_{1}x_{2}^{T} + \dots + \varepsilon_{m-1}x_{m-1}x_{m}^{T}	 \, ,
\end{equation}
for some arbitrarily small $\varepsilon_1, \dots, \varepsilon_{m-1}>0$, if necessary.  As a result, the matrix $A-D$ has a single Jordan block of size $m$ for the eigenvalue $0.$

\end{dem}

Example \ref{ex8} below shows that the condition \eqref{posit} imposed on $A$ does not exclude Hurwitz matrices. This might seem surprising, as for any eigenvector $x$ corresponding to a real eigenvalue $\lambda < 0$ we have $x^TAx = \lambda \|x\|^2 < 0$. Moreover, it holds that any positive-definite matrix $D$ has only eigenvalues with positive real part, see Lemma \ref{posdefpos} below. This result is well-known, but included here for completeness.
	
\begin{lem}\label{posdefpos}
Let $D \in \mathbb{M}_{n}(\mathbb{R})$ be a real positive-definite matrix (though not necessarily symmetric). That is, assume we have $x^TDx > 0$ for all non-zero $x \in \mathbb{R}^n$. Then, any eigenvalue of $D$ has positive real part.
\end{lem}
	
\begin{dem}
Let $\lambda$ be an eigenvalue of $D$ with corresponding eigenvector $x$. We write $\lambda = \xi + i\zeta$ and $x = u+iv$ for their decomposition into real and imaginary parts. On the one hand, we find
\begin{equation}\label{firsteqposdef}
\bar{x}^TDx = \bar{x}^T\lambda x = \|x\|^2\lambda = \|x\|^2(\xi + i\zeta)\, .
\end{equation}
On the other, we have
\begin{equation}\label{secondeqposdef}
\bar{x}^TDx = (u-iv)^TD(u+iv) = u^TDu + v^TDv + i(u^TDv - v^TDu)\, .
\end{equation}
Comparing the real parts of equations \eqref{firsteqposdef} and  \eqref{secondeqposdef}, we conclude that 
\begin{equation}
\|x\|^2\xi = u^TDu + v^TDv > 0\, ,
\end{equation}
where we use that $u$ and $v$ cannot both be zero. Hence, we see that indeed $\xi > 0$.
\end{dem}
	
	\begin{ex}\label{ex8}
		Consider the $(4 \times 4)$ matrix:
		\begin{equation}
			\nonumber A=\left(\begin{array}{rrrr}
				1&1&0&0\\
				-1&1&1&0\\
				0&-1&1&16.94\\
				1&-4.24&-4.24&-17.94
			\end{array}\right).
		\end{equation}
		The matrix $A$ is Hurwitz and Equation \eqref{posit} holds for $m=3$ with $e_{1}=(1,0,0,0)^T,e_{2}=(0,1,0,0)^T$ and $e_{3}=(0,0,1,0)^T$. We may determine the upper-diagonal $(3\times3)$-matrix $P$ from the proof of Proposition \ref{jordan} by calculating
		\begin{eqnarray}
		\nonumber p_{1,2}&=&e^{T}_{1}(A+A^{T})e_{2}\\
		\nonumber p_{1,3}&=&e^{T}_{1}(A+A^{T})e_{3}\\
		\nonumber p_{2,3}&=&e^{T}_{2}(A+A^{T})e_{3}
		\end{eqnarray}
		where
		\begin{equation}
			\nonumber A+A^{T}=\left(\begin{array}{rrrr}
				2&0&0&1\\
				0&2&0&-4.24\\
				0&0&2&12.7\\
				1&-4.24&12.7&-35.88
			\end{array}\right).
		\end{equation}
		It follows that $p_{1,2}=p_{1,3}=p_{2,3}=0$, which implies we have $P=0$. As in the proof of Proposition  \ref{jordan}, we first define $D\in\mathbb{M}_{4}(\mathbb{R})$ on $\Span(e_1, e_2, e_3) = \{x \in \mathbb{R}^4 \mid x_4 = 0\}$ by setting:
		\begin{eqnarray}
		\nonumber De_{1}&=&Ae_{1}\\
		\nonumber De_{2}&=&Ae_{2}-p_{1,2}e_{1}=Ae_{2}\\
		\nonumber De_{3}&=&Ae_{3}-p_{1,3}e_{1}-p_{2,3}e_{2}=Ae_{3}.
		\end{eqnarray}
		Hence, $D$ agrees with $A$ in the first three columns. To complete our construction of $D$, we have to choose a non-zero vector $u$ such that  $u\perp e_{i}$ for $i=1,2,3$, and set $Du=cu$ for some $c>0$. To this end, we set $u=e_{4}$, so that  $D$ becomes:
		\begin{equation}
			\nonumber D=\left(\begin{array}{rrrr}
			1&1&0&0\\
			-1&1&1&0\\
			0&-1&1&0\\
			1&-4.24&-4.24&c
			\end{array}\right).
		\end{equation}
		It follows that
		\begin{equation}
			\nonumber A-D=\left(\begin{array}{cccc}
			0&0&0&0\\
			0&0&0&0\\
			0&0&0&16.94\\
			0&0&0&-17.94-c
			\end{array}\right)
		\end{equation}
	which has a zero eigenvalue with geometric multiplicity $3$ and a negative eigenvalue $(-17.94-c)$ equal to its trace. Moreover, by the Lemma \ref{Schurlemma}, $D$ is positive-definite for large enough $c>0.$ Indeed, in this case, we numerically found that for all $c\geq 9.24$ is enough.
	\end{ex}
			
	\begin{ex}
		The matrix
		\begin{equation}
			A=\left(\begin{array}{rrrr}
			-6&2&1&-3\\
			 2&-8&-1&-2\\
			 1&-1&-3.4&0\\
			-3&-2&0&-6
			\end{array}\right)
		\end{equation}
		is Hurwitz, but symmetric. Thus, there are no vectors $x 
	\in \mathbb{R}^4$ such that $\langle x,Ax\rangle >0.$
	\end{ex}

To control a bifurcation in the system \eqref{systeminchapter3}, we need to rule out additional eigenvalues laying on the imaginary axis. Recall that the eigenvalues of the linearization \eqref{systeminchapter3lin} are given by those of $A -\alpha \lambda_pD$ with $\lambda_p \geq 0$ an eigenvalue of $L_G$. Lemma \eqref{transv} below shows that generically only one of the matrices $A -\alpha \lambda_pD$ is non-hyperbolic. In what follows we denote by $\|\cdot\|$ the operator norm  induced by the Euclidean norm on $\mathbb{R}^n.$
	
\begin{lem}\label{transv}
Let $A,D\in\mathbb{M}_{n}(\mathbb{R})$ be two given matrices  with $D$ positive-definite, and let $\alpha^* \in \mathbb{R}$ be a positive scalar. We furthermore assume $\{\lambda_1, \dots, \lambda_K\}$ is a set of real numbers and consider the matrices $A  - \alpha^*\lambda_iD$ for $i \in \{1, \dots, K\}$. Given any $\varepsilon > 0$, there exist a matrix $\tilde{A}$ and a positive-definite matrix $\tilde{D}$ such that $\|A  - \tilde{A}\|, \|D  - \tilde{D}\| < \varepsilon$ and $\tilde{A}  - \alpha^*\lambda_K\tilde{D} = A  - \alpha^*\lambda_KD$. Moreover, for $i \in \{1, \dots, K-1\}$ the matrix $\tilde{A}  - \alpha^*\lambda_i\tilde{D}$ has a purely hyperbolic spectrum (i.e. no eigenvalues on the imaginary axis).
\end{lem}
\begin{obs}\label{obs4}
		
From $\|A  - \tilde{A}\|, \|D  - \tilde{D}\| < \varepsilon$ we get 

$$\|({A}  - \alpha^*\lambda_i{D}) - (\tilde{A}  - \alpha^*\lambda_i\tilde{D}))\| \leq \|A - \tilde{A}\| + \alpha^*|\lambda_i|\|D  - \tilde{D}\| < \varepsilon(1+\alpha^*|\lambda_i|)\, ,$$
so that we may arrange for $\tilde{A}  - \alpha^*\lambda_i\tilde{D}$ to be arbitrarily close to the original ${A}  - \alpha^*\lambda_i{D}$ for all $i$. Moreover, if $A$ is hyperbolic then for $\varepsilon$ small enough so is $\tilde{A}$, with the same number of stable and unstable eigenvalues. In particular, $\tilde{A}$ may be assumed Hurwitz if $A$ is.  
	\end{obs}

	\begin{demtranslem}
		Let $\delta \not= 0$ be given and set 
		\begin{align}
			\tilde{A}_{\delta} &:= A + \delta \id \\
			\tilde{D}_{\delta} &:= D + \frac{\delta}{\alpha^* \lambda_K} \id \, .
		\end{align}
		Note that the symmetric parts of $\tilde{D}_{\delta}$ and $D$ differ by $\frac{\delta}{\alpha^* \lambda_K} \id$ as well, so that $\tilde{D}_{\delta}$ remains positive-definite for $|\delta|$ small enough. It is also clear that 
		$$\lim\limits_{\delta \to 0}\|A  - \tilde{A}_{\delta}\| = \lim\limits_{\delta \to 0}\|D  - \tilde{D}_{\delta}\| = 0\, .$$
		A direct calculation shows that 
		\begin{align}
			\tilde{A}_{\delta}-\alpha^*\lambda_i\tilde{D}_{\delta}&=A+\delta \id-\alpha^* \lambda_i \left(D+\frac{\delta}{\alpha^*\lambda_K}\id\right)\\ \nonumber
			&=A+\delta \id-\alpha^*\lambda_i D-\frac{\delta\lambda_i}{\lambda_K}\id\\ \nonumber
			&=(A-\alpha^*\lambda_i D)+\left(1-\frac{\lambda_i}{\lambda_K}\right)\delta  \id\, ,
		\end{align}
		for all $i \in \{1, \dots, K\}$. It follows that we have $\tilde{A}_{\delta}  - \alpha^*\lambda_K\tilde{D}_{\delta} = {A}  - \alpha^*\lambda_K{D}.$ For $i \not= K$ we see that $\tilde{A}_{\delta}  - \alpha^*\lambda_i\tilde{D}_{\delta}$ differs from ${A}  - \alpha^*\lambda_i{D}$ by a non-zero scalar multiple of the identity. It follows that for $\delta \not= 0$ small enough all the matrices $\tilde{A}_{\delta}  - \alpha^*\lambda_i\tilde{D}_{\delta}$ for $i \in \{1, \dots, K-1\}$ have their eigenvalues away from the imaginary axis. Setting $\tilde{A} := \tilde{A}_{\delta}$ and  $\tilde{D} := \tilde{D}_{\delta}$ with $\delta = \delta(\varepsilon)$ small enough then finishes the proof.
	\end{demtranslem}

\subsection{Center manifold reduction}\label{Centermanifoldreduction}
	
	Let us now assume $A$, $D$ and $\alpha^*$ are given such that for a particular eigenvalue $\lambda>0$ of $L_G$ the matrix $A-\alpha^*\lambda D$ has an $m$-dimensional center subspace. We moreover assume $\lambda$ is simple and, motivated by Lemma \ref{transv}, that the matrices $A-\alpha^*\kappa D$ are hyperbolic for any other eigenvalue $\kappa$ of $L_G$. It follows that the linearization 
	\begin{equation}\label{linnowhere}
	    \id_{N} \otimes A-\alpha^* L_G\otimes D: \mathbb{R}^N \otimes \mathbb{R}^n \rightarrow \mathbb{R}^N \otimes \mathbb{R}^n
	\end{equation}
	of \eqref{systeminchapter3lin} has an $m$-dimensional center subspace as well. 
	
	In what follows, we write $\hat{I}_{s}$ for the indices of all remaining eigenvalues of $L_G$ except the index $s$. In other words, writing $0 = \lambda_1 < \lambda_2 \leq \dots \leq \lambda_N$ for the eigenvalues of $L_G$,  we have $\lambda = \lambda_s$ for some $s \in \{2, \dots, N\}$ and we set $\hat{I}_{s} = \{1, \dots, N\}\setminus \{s\}$. We will likewise fix an orthonormal set of eigenvectors $v_1, \dots, v_N$ for $L_G$ and simply write $v =v_s$ for the eigenvector corresponding to our fixed eigenvalue $\lambda = \lambda_s$.  Arguably the most natural situation is given by $s = N$, corresponding to the situation where $\alpha$ is increased until the eigenvalues of $A-\alpha \lambda_N D$ first hit the imaginary axis for $\alpha = \alpha^*$. However, we will not need this assumption here.
	
	Next, given a vector $u \in \mathbb{R}^N$ we denote by $\phi_u: \mathbb{R}^N \rightarrow \mathbb{R}^N$ the linear map defined by 
	\begin{equation}
	    \phi_u(w) = \langle u, w \rangle u\, .
	\end{equation}
	Note that  $\phi_u$ is a projection if $\|u\|=1$.
	Finally, we write $E^c, E^h \subseteq \mathbb{R}^n$ for the center- and hyperbolic subspaces of $A-\alpha^*\lambda D$, respectively. It follows that
	\begin{equation}
	    \mathbb{R}^n = E^c \oplus E^h\, ,
	\end{equation}
	and we denote the projections onto the first and second component  by $\pi^c$ and $\pi^h = \id_n - \pi^c$, respectively. Likewise, we denote the center- and hyperbolic subspaces of the map \eqref{linnowhere} by $\mathcal{E}^c, \mathcal{E}^h \subseteq \mathbb{R}^N \otimes \mathbb{R}^n$. Their projections are denoted by $\Pi^c$ and $\Pi^h$. The following lemma establishes some important relations between the aforementioned maps and spaces.
	
	\begin{lem}\label{centersubspaacesproj}
	The spaces $\mathcal{E}^c$ and $E^c$ are related by
	\begin{equation}\label{ecandec}
	    \mathcal{E}^c = v \otimes E^c\, ,
	\end{equation}
	and we have
	\begin{equation}\label{ehandeh}
	    \mathcal{E}^h = (v \otimes E^h) \bigoplus_{i \in \hat{I}_{s}} (v_i \otimes \mathbb{R}^n)\, .
	\end{equation}
	Moreover, it holds that
	\begin{equation}
	    \Pi^c = \phi_v \otimes \pi^c\, .
	\end{equation}
	\end{lem}
	
\begin{dem}
The identities \eqref{ecandec} and \eqref{ehandeh} follow directly from the fact that the linear map \eqref{linnowhere} sends a vector $v_i \otimes x$ to $v_i \otimes (A-\alpha^*\lambda_iD)(x)$ for all $i \in \{1, \dots, N\}$ and $x \in \mathbb{R}^n$. To show that $\Pi^c$ is indeed given by $\phi_v \otimes \pi^c$, we have to show that the latter vanishes on $\mathcal{E}^h$ and restricts to the identity on $\mathcal{E}^c$. To this end, note that for all $i \in \hat{I}_{s}$ and $x \in \mathbb{R}^n$ we have  
	\begin{equation}
	    (\phi_v \otimes \pi^c)(v_i \otimes x) = \phi_v(v_i) \otimes \pi^c(x) = \langle v,v_i \rangle v \otimes \pi^c(x) = 0\, .
	\end{equation}
	Given $x_h \in E^h$ and $x_c \in E^c$, we find
		\begin{align}
		 (\phi_v \otimes \pi^c)(v \otimes x_h) &= \phi_v(v) \otimes \pi^c(x_h) = 0 \quad \text{ and } \\ \nonumber
	    (\phi_v \otimes \pi^c)(v \otimes x_c) &= \phi_v(v) \otimes \pi^c(x_c) = \langle v,v \rangle v \otimes x_c = v \otimes x_c \, ,
	\end{align}
	so that indeed $(\phi_v \otimes \pi^c)|_{\mathcal{E}^h} = 0$ and $(\phi_v \otimes \pi^c)|_{\mathcal{E}^c} = \id_{\mathcal{E}^c}$. This completes the proof.
	\end{dem}

Next, we investigate the dynamics on a center manifold of the system 
	\begin{equation}\label{systeminchapter4}
 	 		\dot{X}=F(X;\varepsilon)-\alpha^*(L_G\otimes D)X \, ,
 	 	\end{equation}
which will lead to a proof of Theorem \ref{mainformal}.  


Recall that center manifold theory predicts a locally defined map $\Psi: \mathcal{E}^c \times \Omega \rightarrow \mathcal{E}^h$ whose graph $M^c$ is invariant for the system \eqref{systeminchapter4} and locally contains all bounded solutions. The map $\Psi$  satisfies $\Psi(0;0) = 0$ and $D\Psi(0;0) = 0$. In fact, as we assume $F(0;\varepsilon) = 0$ for all $\varepsilon \in \Omega$, it follows that $(0;\varepsilon) \in M^c$, as these are bounded solutions. This shows that $\Psi(0;\varepsilon) = 0$ for all $\varepsilon \in \Omega$. 
	
In light of Lemma \ref{centersubspaacesproj}, we may write
\begin{align}\label{formofpsii}
\Psi(v \otimes x_c; \varepsilon) = v \otimes \psi(x_c;\varepsilon) + \sum_{i \in \hat{I}_{s}}v_i \otimes \psi_i(x_c;\varepsilon)\, ,
\end{align}
for certain maps $\psi: E^c \times \Omega \rightarrow E^h$ and $\psi_i: E^c \times \Omega \rightarrow \mathbb{R}^n$. We then likewise have $\psi(0;\varepsilon) = 0, D\psi(0;0) = 0$ and $\psi_i(0;\varepsilon) = 0, D\psi_i(0;0) = 0$ for all $\varepsilon \in \Omega$ and $i \in \hat{I}_{s}$.

The dynamics on the center manifold $M^c$ is conjugate to that of a vector field on $\mathcal{E}^c \times \Omega$ given by
\begin{align}
    \tilde{R}(X_c ; \varepsilon) = \Pi^c S(X_c, \Psi(X_c;\varepsilon);\varepsilon)\, ,
\end{align}
where we write 
\begin{equation}
    S(X_c, X_h; \varepsilon) = F(X_c+X_h;\varepsilon)-\alpha^*(L_G\otimes D)(X_c+X_h) 
\end{equation}
for the vector field on the right hand side of \eqref{systeminchapter4}, with $X_c \in \mathcal{E}^c$, $X_h \in \mathcal{E}^h$ and $\varepsilon \in \Omega$. We further conjugate to a vector field $R$ on $E^c \times \Omega$ by setting 
\begin{align}\label{RinRtilde}
 v \otimes R(x_c;\varepsilon) = \tilde{R}(v \otimes x_c ; \varepsilon)\, .
\end{align}
In order to describe $R$, we first introduce some useful notation. Given $X \in \mathbb{R}^N \otimes \mathbb{R}^n$, we may write 
\begin{equation}\label{decompi23}
    X = \sum_{p=1}^N e_p \otimes x_p
\end{equation}
with $e_1, \dots, e_N$ the canonical basis of $\mathbb{R}^N$ and for some unique vectors $x_p \in \mathbb{R}^n$. In general, given $p \in \{1, \dots, N\}$ we will denote by $x_p \in \mathbb{R}^n$ the $p$th component of $X$ as in the decomposition \eqref{decompi23}. Recall that $v = (\nu_1,\dots,\nu_N)$ is the eigenvector associated with $\lambda$.
Using this notation, we have the following result. 

\begin{prop}\label{non-linearterms}
Denote by $h: \mathbb{R}^n \times \Omega \rightarrow \mathbb{R}^n$ the non-linear part of $f$. That is, we have $f(x;\varepsilon) = Ax + h(x;\varepsilon)$. The reduced vector field $R$ is given explicitly by \begin{align}
R(x_c; \varepsilon) = (A-\alpha^*\lambda D)x_c + \sum_{p=1}^N \nu_p \pi^ch(\nu_px_c +\Psi(x_c \otimes v;\varepsilon)_p;\varepsilon)\, ,
\end{align}
for $x_c \in E^c$ and $\varepsilon \in \Omega$.
\end{prop}

\begin{dem}
We write $S(X; \varepsilon) = S(X_c, X_h; \varepsilon)= JX + H(X; \varepsilon)$ with $J = D_X S(0;0)$ and where $H$ denotes higher order terms. It follows that 
\begin{equation}\label{piGintwoterms}
\Pi^c S(X_c, \Psi(X_c;\varepsilon);\varepsilon) = \Pi^c J(X_c+\Psi(X_c;\varepsilon)) + \Pi^cH(X_c, \Psi(X_c;\varepsilon);\varepsilon)\, .
\end{equation}
We start by focusing on the first term. As $J$ sends $\mathcal{E}^c$ to $\mathcal{E}^c$ and $\mathcal{E}^h$ to $\mathcal{E}^h$, we conclude that $\Pi^c J = J \Pi^c$. We therefore find 
\begin{equation}
 \Pi^c J (X_c+\Psi(X_c;\varepsilon)) =   J \Pi^c(X_c+\Psi(X_c;\varepsilon)) = J X_c\, .
\end{equation}
Writing $X_c = v \otimes x_c$ and using Expression \eqref{systeminchapter3lin} for $J$, we conclude that the linear part of $\tilde{R}$ is given by
\begin{equation}\label{linearreduced44}
    J(v \otimes x_c) = v \otimes (A-\alpha^*\lambda D)x_c\, .
\end{equation}

We next focus on the second term in Equation \eqref{piGintwoterms}. Note that we have 
\begin{equation}
    H(X;\varepsilon)_p = h(x_p;\varepsilon)\, \text{ for all } \, p \in \{1, \dots, N\}\, .
\end{equation}
Now, by Lemma \ref{centersubspaacesproj} it follows that we may write 
\begin{align}
    \Pi^c(X) &= \sum_{p=1}^N\Pi^c( e_p \otimes x_p) = \sum_{p=1}^N \phi_v(e_p) \otimes \pi^c(x_p) = \sum_{p=1}^N \langle v, e_p \rangle v \otimes \pi^c(x_p) = \sum_{p=1}^N \nu_p (v \otimes \pi^c(x_p)) \\ \nonumber
    &= v \otimes \sum_{p=1}^N \nu_p  \pi^c(x_p)\, \text{ for all } X \in \mathbb{R}^N \otimes \mathbb{R}^n.
\end{align}
We therefore find
\begin{align}
\Pi^c(H(X_c, \Psi(X_c;\varepsilon);\varepsilon)) = v \otimes \sum_{p=1}^N \nu_p \pi^ch([X_c+ \Psi(X_c;\varepsilon)]_p;\varepsilon)\, .
\end{align}
Next, we have $(X_c)_p = (v \otimes x_c)_p = \nu_p x_c $, so that we find
\begin{align}\label{nonlinearreduced44}
\Pi^c(H(X_c, \Psi(X_c;\varepsilon);\varepsilon)) = v \otimes \sum_{p=1}^N \nu_p \pi^ch(v_px_c +\Psi(X_c;\varepsilon)_p;\varepsilon)\, .
\end{align}
Combining equations \eqref{linearreduced44} and \eqref{nonlinearreduced44}, we arrive at
\begin{align}\label{totalreduced44}
\tilde{R}(X_c; \varepsilon) = \Pi^c( S (X_c, \Psi(X_c;\varepsilon);\varepsilon)) = v \otimes (A-\alpha^*\lambda D)x_c + v \otimes \sum_{p=1}^N \nu_p \pi^ch(\nu_px_c +\Psi(X_c;\varepsilon)_p;\varepsilon)\, .
\end{align}
Finally, from Equation \ref{RinRtilde} we get
\begin{align}
R(x_c; \varepsilon) = (A-\alpha^*\lambda D)x_c + \sum_{p=1}^N \nu_p \pi^ch(\nu_px_c +\Psi(v \otimes x_c;\varepsilon)_p;\varepsilon)\, ,
\end{align}
which completes the proof.
\end{dem}

To further investigate the Taylor expansion of $R(x_c; \varepsilon)$, we need to know more about how the coefficients of $\Psi: \mathcal{E}^c \times \Omega \rightarrow \mathcal{E}^h$ depend on those of $F$.

To this end, let us consider for a moment the general situation  where $S$ is some vector field on $\mathbb{R}^k$ (in our case $k = nN$) satisfying $S(X) = JX + H(X)$ for some $H: \mathbb{R}^n \rightarrow \mathbb{R}^n$ satisfying $H(0)= 0$, $DH(0) = 0$. We furthermore let $\hat{\mathcal{E}}^c$ and $\hat{\mathcal{E}}^h$ denote the center- and hyperbolic subspaces of $J$, respectively, and write $\Pi^c, \Pi^h$ for the corresponding projections. 
Suppose $\Psi: \hat{\mathcal{E}}^c \to \hat{\mathcal{E}}^h$ is a locally defined map whose graph is a center manifold $M^c$ for the system $\dot{X} = S(X)$. Recall that we have $\Psi(0) = 0$ and $D\Psi(0) = 0$. Moreover, as $M^c$ is a flow-invariant manifold, we see that $S|_{M^c}$ takes values in the tangent bundle of $M^c$. This can be used to iteratively solve for the higher order coefficients of an expansion of $\Psi$ around $0$. 

More precisely, the tangent space at $X_c + \Psi(X_c) \in M^c$ is given by all vectors of the form $(V_c, D\Psi(X_c)V_c) \in \hat{\mathcal{E}}^c \oplus \hat{\mathcal{E}}^h$, with $V_c \in \hat{\mathcal{E}}^c$. Invariance under the flow of $S$ then translates to the identity
\begin{align}\label{identitytaylorcoeffs}
    \Pi^h J \Psi(X_c) + \Pi^h H(X_c, \Psi(X_c)) = D\Psi(X_c)(\Pi^c J X_c + \Pi^cH(X_c, \Psi(X_c)))\, ,
\end{align}
for $X_c$ in some open neighborhood of the origin in $\hat{\mathcal{E}}^c$. Equation \eqref{identitytaylorcoeffs} can be used to show that $D\Psi(0) = 0$. Equation~\eqref{identitytaylorcoeffs} can also be arranged to
\begin{align}\label{identitytaylorcoeffs3}
	\Pi^{h} J \Psi(X_{c})-D\Psi(X_{c})\Pi^{c}J X_{c}=D\Psi(X_{c})\Pi^{c}H(X_{c},\Psi(X_{c}))-\Pi^{h}H(X_{c},\Psi(X_{c}))\, .
\end{align}
As $\Psi$ only has terms of degree $2$ and higher, the same holds for both sides of  Equation~\eqref{identitytaylorcoeffs3}, which depend on $\Psi$ and $H$. More generally, using $\Psi^{\rho}$ to denote the terms of order $\rho \geq 2$ in the Taylor expansion of $\Psi$ around the origin, Equation \eqref{identitytaylorcoeffs3} is readily seen to imply for each $\rho$ 
\begin{align}\label{identitytaylorcoeffs2}
    \Pi^h J \Psi^{\rho}(X_c) - D\Psi^{\rho}(X_c)\Pi^c J X_c = P_{\rho}(X_c)\, ,
\end{align}
    for some homogeneous polynomial $P_{\rho}$ of order $\rho$. Moreover, $P_{\rho}$ depends only on $\Psi^{2}(X_c) \dots, \Psi^{\rho-1}(X_c)$ and on the Taylor expansion of $H$ up to order $\rho$. It can be shown that for fixed $J$ and $P_{\rho}$, Equation \eqref{identitytaylorcoeffs2} has a unique solution $\Psi^{\rho}$ in the form of a homogeneous polynomial of order $\rho$, see \cite{wimmer1979equation}. As a result, we get the following important observation:

    \begin{lem}\label{remarkgeneraldepntaylor}
    We may iteratively solve for the terms $\Psi^{\rho}$ using expression \eqref{identitytaylorcoeffs2}. Moreover, for fixed linearity $J$, the terms of order $\rho$ and less of $\Psi$ are fully determined by the terms of order $\rho$ and less of $H$.
    \end{lem} 
   
We return to our main setting where $R(x_c;\varepsilon)$ is the reduced vector field of the system \eqref{systeminchapter4} as described in Proposition \ref{non-linearterms}. Note that the presence of a parameter $\varepsilon$ means that the center subspace $\hat{\mathcal{E}}^c$ in the observations for general vector fields above is now given by $\mathcal{E}^c \times \Omega$.
    
    \begin{lem}\label{transvch4}
    Let $\rho > 1$ be given, and suppose the vector $v = (\nu_1, \dots, \nu_N) \in \mathbb{R}^N$ satisfies 	
    	\begin{equation}\label{specialpropertych4}
		\sum_{p=1}^{N}\nu^{\ell}_{p}\neq0\, ,\quad\forall~\ell=2,\dots,\rho+1.
	\end{equation}
	Then the reduced vector field $R(x_c;\varepsilon)$ as described in Proposition \ref{non-linearterms} can have any Taylor expansion around $0$ of order $2$ to $\rho$, subject to $R(0;\varepsilon) = 0$, if no conditions are put on the nonlinear part of $f$ other than $f(0;\varepsilon) = 0$ and sufficient smoothness.
    \end{lem}
    
    \begin{dem}
    From Proposition \ref{non-linearterms} we know that 
    \begin{align}\label{reduced1inprooff}
    R(x_c; \varepsilon) = (A-\alpha^*\lambda D)x_c + \sum_{p=1}^N \nu_p \pi^ch(\nu_px_c +\Psi(v \otimes x_c;\varepsilon)_p;\varepsilon)\, .
    \end{align}
    As we have $\Psi(0;\varepsilon) = 0$ and $h(0;\varepsilon) = 0$, we conclude that likewise $R(0; \varepsilon) = 0$ for all $\varepsilon \in \Omega$. In particular, we see that $D_{\varepsilon}R(0;0) = 0$, whereas Equation \eqref{reduced1inprooff} tells us that $D_{x_c}R(0;0) = (A-\alpha^*\lambda D)|_{E^c}$.
    
   \noindent As a warm-up, we start by investigating the second order terms of $R.$ To this end, we write 
    \begin{equation}\label{hotofR}
    h(x;\varepsilon) = Q_{1,1}(x; \varepsilon) + Q_{2,0}(x) + \mathcal{O}(|(x,\varepsilon)|^3)\, ,
    \end{equation}
    where $Q_{1,1}$ is linear in both components and $Q_{2,0}$ is a quadratic form. It follows that 
\begin{eqnarray}\label{hofRinsecondorder}
h(\nu_px_c +\Psi(v \otimes x_c;\varepsilon)_p;\varepsilon) &=&Q_{1,1}(\nu_px_c +\Psi(v \otimes x_c;\varepsilon)_p; \varepsilon) + Q_{2,0}(\nu_px_c +\Psi(v \otimes x_c;\varepsilon)_p) + \mathcal{O}(|(x_c,\varepsilon)|^3) \nonumber \\
     &=&Q_{1,1}(\nu_px_c; \varepsilon) + Q_{2,0}(\nu_px_c) + \mathcal{O}(|(x_c,\varepsilon)|^3)\,, \nonumber
    \end{eqnarray}
    where we use that $\Psi(v \otimes x_c;\varepsilon)$ has no constant or linear terms in $(x_c; \varepsilon)$. From Equation \eqref{hotofR} we obtain
\begin{eqnarray}\label{reduced1inprooff3}
 \sum_{p=1}^N \nu_p \pi^ch(v_px_c +\Psi(v \otimes x_c;\varepsilon)_p;\varepsilon)  
 &=& \pi^c \sum_{p=1}^N \nu_p Q_{1,1}(\nu_px_c; \varepsilon) + \pi^c \sum_{p=1}^N \nu_pQ_{2,0}(\nu_px_c) + \mathcal{O}(|(x_c,\varepsilon)|^3)  \nonumber \\
& = & \pi^c \sum_{p=1}^N \nu_p^2 Q_{1,1}(x_c; \varepsilon) + \pi^c \sum_{p=1}^N \nu_p^3Q_{2,0}(x_c) + \mathcal{O}(|(x_c,\varepsilon)|^3) \, . \nonumber
   \end{eqnarray}
  As we assume  $\sum_{p=1}^N \nu_p^3 \not= 0$, we see that the second order Taylor coefficients of $R(x_c; \varepsilon)$ can be chosen freely (except for the $\mathcal{O}(|\varepsilon|^2)$ term).
  
  Now suppose we are given a polynomial map $P: E^c \times \Omega \to E^c$ of degree $\rho$ satisfying $DP(0) = ((A-\alpha^*\lambda D)|_{E^c};0)$ and $P(0;\varepsilon) = 0$ for all $\varepsilon$. We will prove by induction that we may choose the terms in the Taylor expansion of $h$ up to order $\rho$ in the variables $x_c$ and $\varepsilon$ such that the Taylor expansion up to order $\rho$ of $R$ agrees with $P$. To this end, suppose some choice of $h$ gives agreement between $P$ and the Taylor expansion of $R$ up to order $2\leq k < \rho$. By the foregoing, this can be arranged for $k=2$.
  
  We start by remarking that a change to $h$ that does not influence its Taylor expansion up to order $k$ does not change the Taylor expansion of $\Psi$ up to order $k$. This follows directly from Lemma \ref{remarkgeneraldepntaylor}. As a result, such a change does not influence the Taylor expansion of $R$ up to order $k$ as well. We write
  \begin{align}
      \tilde{h}(x_c; \varepsilon) = h(x_c; \varepsilon) + \sum_{i=1}^{k+1} Q_{i,k+1-i}(x_c; \varepsilon)\,  
  \end{align}
  for an order $k+1$ change to $h$, where each component of $Q_{i,j}: \mathbb{R}^n \times \Omega \to \mathbb{R}^n$ is a homogeneous polynomial of degree $i$ in $x_c$ and degree $j$ in $\varepsilon$. The $(k+1)$-order terms of $R$ in $(x_c; \varepsilon)$ are given by the $(k+1)$-order terms of 
\begin{align}\label{reduced1inprooff1}
  &\sum_{p=1}^N \nu_p \pi^c\tilde{h}(\nu_px_c +\Psi(v \otimes x_c;\varepsilon)_p;\varepsilon) \\ \nonumber
  = \,  &\sum_{p=1}^N \nu_p \pi^c\bigg(h(\nu_px_c +\Psi(v \otimes x_c;\varepsilon)_p;\varepsilon) + \sum_{i=1}^{k+1} Q_{i,k+1-i}(\nu_px_c +\Psi(v \otimes x_c;\varepsilon)_p; \varepsilon)\bigg)\, .
 \end{align}
  As both $h$ and $\Psi$ have no constant and linear terms, we see that  the $(k+1)$-order terms of $R$ are also given by those of
  \begin{align}\label{reduced1inprooff2}
  \sum_{p=1}^N \nu_p \pi^c\bigg(h(\nu_px_c +\Psi^k(v \otimes x_c;\varepsilon)_p;\varepsilon) + \sum_{i=1}^{k+1} Q_{i,k+1-i}(\nu_p x_c ; \varepsilon)\bigg)\, ,
 \end{align}
  where $\Psi^{k}$ denotes the terms of $\Psi$ up to order $k$. As we have previously argued, $\Psi^{k}$ is independent of the additional terms $Q_{i,k+1-i}$. Hence, we may write the order $k+1$ terms in Expression \eqref{reduced1inprooff2} as 
  \begin{align}
 W(x_c ; \varepsilon) +  \sum_{p=1}^N \nu_p \pi^c \sum_{i=1}^{k+1} Q_{i,k+1-i}(\nu_p x_c ; \varepsilon) &= W(x_c ; \varepsilon) +  \sum_{p=1}^N \nu_p \pi^c \sum_{i=1}^{k+1} \nu_p^iQ_{i,k+1-i}(x_c ; \varepsilon) \\ \nonumber
 &=  W(x_c ; \varepsilon) +  \sum_{i=1}^{k+1} \bigg(\sum_{p=1}^{N} \nu_p^{i+1}\bigg) \pi^c  Q_{i,k+1-i}(x_c ; \varepsilon)\, ,
 \end{align}
  where $W(x_c ; \varepsilon)$ denotes the order $k+1$ terms of 
  \[    \sum_{p=1}^N \nu_p \pi^c h(\nu_p x_c +\Psi^k(v \otimes x_c;\varepsilon)_p;\varepsilon) \, .   \] 
  As $\sum_{p=1}^{N} v_p^{j} \not= 0$ for all $j \in \{2, \dots, \rho+1\}$, we see that the order $k+1$ terms of $R$ may be freely chosen. In other words, we may arrange for the Taylor expansion up to order $k+1$ of $R$ to agree with that of $P$ up to order $k+1$. This completes the proof by induction.
 \end{dem}
    
\begin{demmain}
If there exist $m$ mutually orthogonal vectors $x_1,\dots,x_m$ such that $\langle x_i, Ax_i \rangle > 0$, then Proposition \ref{jordan} guarantees the existence of a positive-definite matrix $D$ such that  $A-D$ has a center subspace of dimension $m$ or higher. Given any non-zero eigenvalue $\lambda$ of $L_G$, we may set $\alpha^* = 1/\lambda$ and conclude that $A-\alpha^*\lambda D$ has a center subspace of dimension at least $m$. 
As the eigenvalues of the linearization of \eqref{systeminainthr} around the origin are given by those of the maps $A-\alpha \lambda D$ for $\lambda$ an eigenvalue of $L_G$, we see that the system \eqref{systeminainthr} has a local parameterized center manifold of dimension at least $m$ for some choices of $D$ and $\alpha = \alpha^*$.

If the graph $G$ is $\rho$-versatile for the pair $(\lambda, v)$, then a choice of $D$ as above together with $\alpha^* = 1/\lambda$ guarantees $A- D$ has a center subspace of dimension at least $m$. By Corollary \ref{obsaboutremainingeig} we may assume this center subspace to be of dimension precisely $m$. 
Moreover, by Lemma \ref{transv} we may assume $A-\alpha^*\lambda_i D$ to have a hyperbolic spectrum for all other eigenvalues $\lambda_i \not= \lambda$ of $L_G$, after an arbitrarily small perturbation to $A$ and $D$ if necessary. 
It follows that the system \eqref{systeminainthr} has a local parameterized center manifold of dimension exactly $m$. 
We argue in the proof of Lemma \ref{transvch4} that $R(0;\varepsilon) = 0$ for all $\varepsilon$, and that $DR(0;0) = ((A-\alpha^*\lambda D)|_{E^c};0)$. This latter map is nilpotent by the statement of Proposition \ref{jordan}. Finally, Lemma \ref{transvch4}  shows that any Taylor expansion can be realized for $R$ up to order $\rho$, subject to the aforementioned restrictions.
\end{demmain}

	\section{Stability of the center manifold}
	
    In this section we investigate the stability of the center manifold of the full network system. 
    We know that the spectrum of the linearization of this system is fully understood if we know the spectrum of the matrices $A-\alpha^*\lambda_i D$ for $\lambda_i$ an eigenvalue of $L_G$.
    Proposition \ref{jordan} gives conditions on $A$ that guarantee the existence of a positive-definite matrix $D$ such that $A-\alpha^*\lambda D$ has an $m$-dimensional generalized kernel for some fixed eigenvalue $\lambda>0$ of $L_G$.
    Moreover, by Corollary \ref{obsaboutremainingeig} we may assume that the non-zero eigenvalues of $A-\alpha^*\lambda D$ have negative real parts.
    Lemma \ref{transv} in turn shows that --after a small perturbation of $A$ and $D$ if necessary-- we may assume $A-\alpha^* \lambda_{i} D$ to have a hyperbolic spectrum for all remaining eigenvalues $\lambda_i \not= \lambda$ of $L_G$.
    Thus, if the matrices $A-\alpha^*\lambda_i D$ for these remaining eigenvalues are all Hurwitz, then the $m$-dimensional center manifold of Theorem \ref{mainformal} may be assumed stable.
    
\begin{figure}[H]
		    \begin{minipage}[t]{0.45\textwidth}
			\centering

\tikzset{every picture/.style={line width=0.75pt}} 

\begin{tikzpicture}[x=0.75pt,y=0.75pt,yscale=-1,xscale=1]

\draw [color={rgb, 255:red, 155; green, 155; blue, 155 }  ,draw opacity=1 ] (62,140.23) -- (238,140.23)(150,52.5) -- (150,212) (231,135.23) -- (238,140.23) -- (231,145.23) (145,59.5) -- (150,52.5) -- (155,59.5)  ;
\draw  [color={rgb, 255:red, 208; green, 2; blue, 27 }  ,draw opacity=1 ][fill={rgb, 255:red, 208; green, 2; blue, 27 }  ,fill opacity=1 ] (102.38,94.25) .. controls (102.38,92.8) and (103.55,91.63) .. (105,91.63) .. controls (106.45,91.63) and (107.63,92.8) .. (107.63,94.25) .. controls (107.63,95.7) and (106.45,96.88) .. (105,96.88) .. controls (103.55,96.88) and (102.38,95.7) .. (102.38,94.25) -- cycle ;
\draw [color={rgb, 255:red, 208; green, 2; blue, 27 }  ,draw opacity=1 ][line width=1.5]  [dash pattern={on 5.63pt off 4.5pt}]  (105,94.25) .. controls (139.78,97.57) and (144.84,118.07) .. (149.11,136.43) ;
\draw [shift={(150,140.23)}, rotate = 256.54] [fill={rgb, 255:red, 208; green, 2; blue, 27 }  ,fill opacity=1 ][line width=0.08]  [draw opacity=0] (15.6,-3.9) -- (0,0) -- (15.6,3.9) -- cycle    ;

\draw  [color={rgb, 255:red, 245; green, 166; blue, 35 }  ,draw opacity=1 ][fill={rgb, 255:red, 245; green, 166; blue, 35 }  ,fill opacity=1 ] (87.48,140.1) .. controls (87.48,138.65) and (88.65,137.48) .. (90.1,137.48) .. controls (91.55,137.48) and (92.73,138.65) .. (92.73,140.1) .. controls (92.73,141.55) and (91.55,142.73) .. (90.1,142.73) .. controls (88.65,142.73) and (87.48,141.55) .. (87.48,140.1) -- cycle ;
\draw [color={rgb, 255:red, 245; green, 166; blue, 35 }  ,draw opacity=1 ][line width=1.5]  [dash pattern={on 5.63pt off 4.5pt}]  (90,140.1) -- (51.33,140.16) ;
\draw [shift={(47.33,140.17)}, rotate = 359.91] [fill={rgb, 255:red, 245; green, 166; blue, 35 }  ,fill opacity=1 ][line width=0.08]  [draw opacity=0] (13.4,-6.43) -- (0,0) -- (13.4,6.44) -- (8.9,0) -- cycle    ;
\draw  [color={rgb, 255:red, 126; green, 211; blue, 33 }  ,draw opacity=1 ][fill={rgb, 255:red, 126; green, 211; blue, 33 }  ,fill opacity=1 ] (121.3,140.28) .. controls (121.3,138.83) and (122.48,137.65) .. (123.93,137.65) .. controls (125.37,137.65) and (126.55,138.83) .. (126.55,140.28) .. controls (126.55,141.72) and (125.37,142.9) .. (123.93,142.9) .. controls (122.48,142.9) and (121.3,141.72) .. (121.3,140.28) -- cycle ;
\draw [color={rgb, 255:red, 126; green, 211; blue, 33 }  ,draw opacity=1 ][fill={rgb, 255:red, 126; green, 211; blue, 33 }  ,fill opacity=1 ][line width=1.5]  [dash pattern={on 5.63pt off 4.5pt}]  (123.93,140.28) -- (146,140.23) ;
\draw [shift={(150,140.23)}, rotate = 179.89] [fill={rgb, 255:red, 126; green, 211; blue, 33 }  ,fill opacity=1 ][line width=0.08]  [draw opacity=0] (15.6,-3.9) -- (0,0) -- (15.6,3.9) -- cycle    ;
\draw [color={rgb, 255:red, 74; green, 144; blue, 226 }  ,draw opacity=1 ][line width=1.5]  [dash pattern={on 5.63pt off 4.5pt}]  (104.91,185.36) .. controls (139.69,182.1) and (144.75,161.99) .. (149.02,143.97) ;
\draw [shift={(149.91,140.25)}, rotate = 103.7] [fill={rgb, 255:red, 74; green, 144; blue, 226 }  ,fill opacity=1 ][line width=0.08]  [draw opacity=0] (15.6,-3.9) -- (0,0) -- (15.6,3.9) -- cycle    ;
\draw  [color={rgb, 255:red, 74; green, 144; blue, 226 }  ,draw opacity=1 ][fill={rgb, 255:red, 74; green, 144; blue, 226 }  ,fill opacity=1 ] (102.29,185.36) .. controls (102.29,186.78) and (103.46,187.93) .. (104.91,187.93) .. controls (106.36,187.93) and (107.54,186.78) .. (107.54,185.36) .. controls (107.54,183.94) and (106.36,182.78) .. (104.91,182.78) .. controls (103.46,182.78) and (102.29,183.94) .. (102.29,185.36) -- cycle ;
\draw  [color={rgb, 255:red, 0; green, 0; blue, 0 }  ,draw opacity=1 ][fill={rgb, 255:red, 0; green, 0; blue, 0 }  ,fill opacity=1 ] (148.77,140.23) .. controls (148.77,140.89) and (149.32,141.43) .. (150,141.43) .. controls (150.68,141.43) and (151.23,140.89) .. (151.23,140.23) .. controls (151.23,139.56) and (150.68,139.02) .. (150,139.02) .. controls (149.32,139.02) and (148.77,139.56) .. (148.77,140.23) -- cycle ;

\draw (225,142.4) node [anchor=north west][inner sep=0.75pt]    {$x$};
\draw (133.5,56.9) node [anchor=north west][inner sep=0.75pt]    {$iy$};

\end{tikzpicture}

            \caption{Sketch of a situation where the $m$-dimensional center manifold of Theorem \ref{mainformal} may be assumed stable. Depicted are the eigenvalues of $A-\beta D$ as $\beta$ is varied. Colored dots denote starting points where $\beta=0,$ the blue and red dashed paths form a conjugate pair of complex eigenvalues and the green and orange dashed paths are real eigenvalues. The arrows indicate how the eigenvalues evolve as $\beta$ increases to $1.$ Three of them go to the  origin, whereas one moves away from it. None of them touches the imaginary axis before $\beta=1.$}
		    \label{fig:11}
		    \end{minipage}
		    \hfill
		    \begin{minipage}[t]{0.45\textwidth}
			\centering

		     \centering
		     \includegraphics[width=6cm]{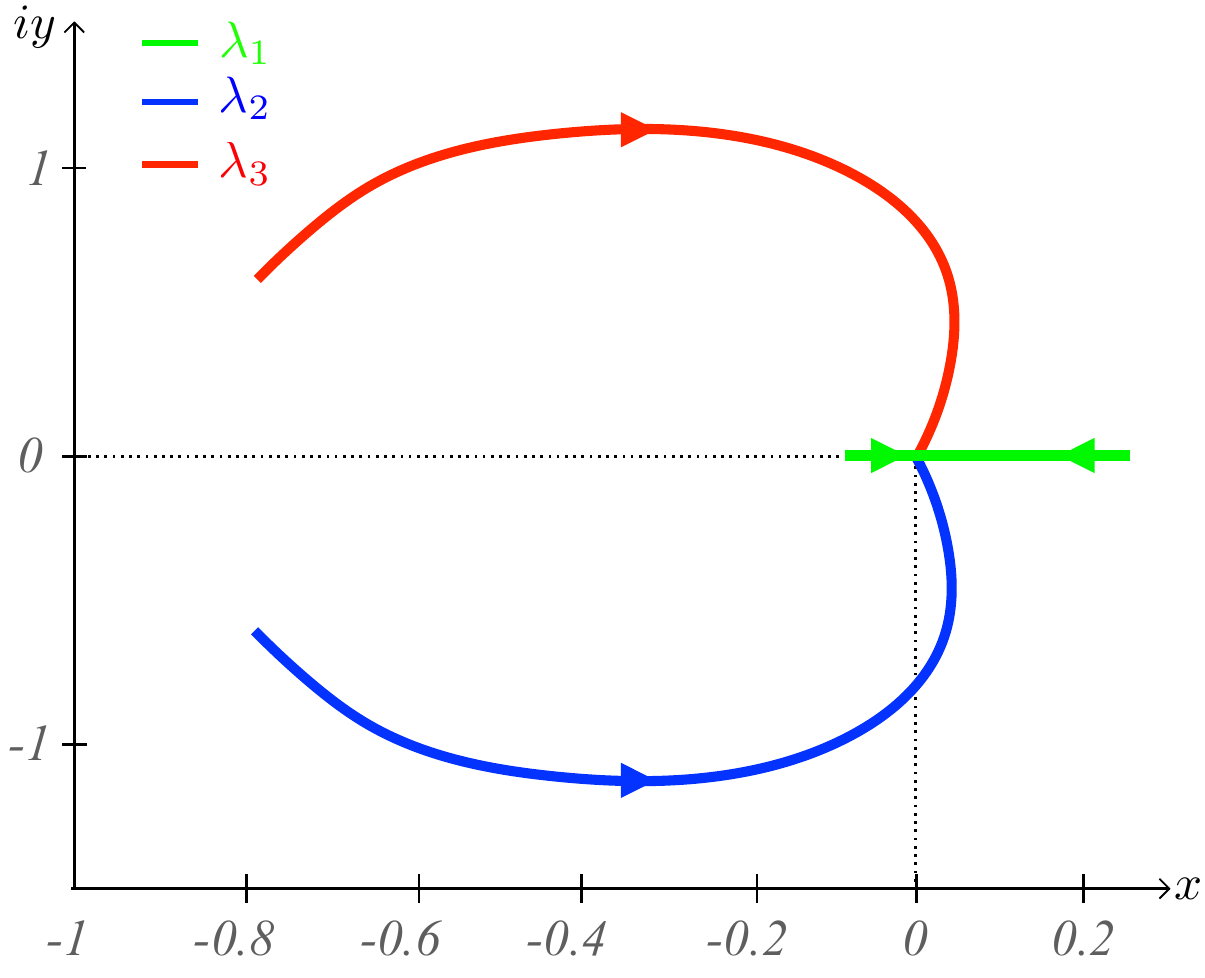}
\caption{Numerically computed behavior of three of the four eigenvalues of the family of matrices from Example \ref{ex:10}. As $\beta$ increases from $0$ to $1$, three eigenvalues move to the origin, whereas a fourth stays to the left of the imaginary axis. For some value of $\beta \in (0,1)$, two complex conjugate eigenvalues already cross the imaginary axis away from the origin. Likewise, a real eigenvalue crosses the origin for some $\beta \in (0,1)$. Data was simulated using Octave. }
\label{fig:12}
\end{minipage}
\end{figure}

    This seems most reasonable to expect when $\lambda$ is the (simple) largest eigenvalue of $L_G$, as the matrices $A-\alpha^*\lambda_i D$ for the remaining eigenvalues of $L_G$ then "lie between" the Hurwitz matrix $A$ and the non-invertible matrix $A-\alpha^*\lambda D$. 
    More precisely, suppose $D$ is scaled such that $A-D = A-\alpha^*\lambda D$. 
    If we let $\alpha$ vary from $0$ to $\alpha^* = 1/\lambda$, then for each eigenvalue $\lambda_i$ of $L_G$, the matrix $A-\alpha \lambda_i D$ is of the form $A-\beta D$ for some $\beta$ in $[0,1]$. 
    Let us therefore denote by $\beta \mapsto \gamma_i(\beta)$ for $i \in \{1, \dots, n\}$ a number of curves through the complex plane capturing the eigenvalues of $A-\beta D$.
    As $\alpha$ varies from $0$ to $\alpha^* = 1/\lambda$, the eigenvalues of $A-\alpha \lambda D$ traverse $\gamma_i$, with the "front runners" given by those of $A-\alpha \lambda D$. 
    In contrast,  for $\lambda_1 = 0$ the eigenvalues of $A-\alpha \lambda_1 D$ of course remain at  $\gamma_i(0)$.
    When $\alpha = \alpha^*$ is reached, the eigenvalues of $A-\alpha^*\lambda_i D$ end up in different places on the curves $\gamma_i$. 
    Hence, if the situation is as in Figure \ref{fig:11}, where each $\gamma_i$ hits the imaginary axis only for $\beta = 1$, or not at all, then we are guaranteed that each of the matrices $A-\alpha^*\lambda_i D$ is Hurwitz for $\lambda_i \not= \lambda$.
    Hence, the center manifold is then stable.
    
    Of course $\beta = 1$ may not be the first value for which a curve $\gamma_i$ hits the imaginary axis, see Figure \ref{fig:12}. 
    Note that, if the matrix $A- \beta D$ indeed has a non-trivial center subspace for some value $\beta \in (0,1)$, then a bifurcation is expected to occur as $\alpha$ is increased, before it hits $\alpha^*$.
	
The next example shows that some of the eigenvalues of $A-\beta D$ might cross the imaginary axis before a high-dimensional kernel emerges at $\beta=1$, see Figure \ref{fig:12}. 
\begin{ex}\label{ex:10}
We consider the matrices $A$ and $D$ constructed in Example \ref{ex8}. Here we choose $c=21$ in order to guarantee that $D$ is a positive-definite matrix. We therefore have the family of matrices:
		\begin{equation}\label{exapositparameterc}
			A-\beta D=\left(\begin{array}{cccc}
				1-\beta & 1-\beta         & 0              & 0\\
				\beta-1 & 1-\beta         & 1 - \beta      & 0\\
				0    & \beta-1         & 1 - \beta      & 16.94\\
				1-\beta & 4.24(\beta-1)  & 4.24(\beta-1) & -21\beta-17.94
			\end{array}\right)\, ,
		\end{equation}
parameterized by the real number $\beta\in[0,1].$ We are interested in the eigenvalue behavior as $\beta$ is varied. We know that for $\beta=0$ we have the Hurwitz matrix $A$, so that all eigenvalues have negative real parts. We would like to know if the family $A-\beta D$ has all eigenvalues with a negative real part for all $\beta\in(0,1).$ However, if $\beta=\frac{1}{2}$ we have $3$ eigenvalues with positive real part. By continuity of the eigenvalues, it means that each of these crossed the imaginary axis for some $\beta<\frac{1}{2}.$ Only after this, for $\beta=1$, do we have the bifurcation studied in the previous chapter, due to the appearance of a triple zero eigenvalue. Figure \ref{fig:12} shows the numerically computed behavior of these three eigenvalue-branches. 
	\end{ex}
	
If we are in the situation of Figure \ref{fig:12}, then the center manifold can still be stable. This occurs when the largest eigenvalue $\mu$ of $L_{G}$ is significantly larger than all other eigenvalues. In that case, we have 
$$A-\alpha^*\lambda D = A-\lambda/ \mu D \approx A$$
for all eigenvalues $\lambda$ of $L_{G}$ unequal to $\mu$. For small enough values of $\lambda/ \mu $ the matrix $A-\lambda/ \mu D$ is therefore still Hurwitz. As it turns out, this can be achieved in the situation explored in Subsection \ref{stargraphs}. More precisely, we have the following result.
	
\begin{prop}\label{eigenvaluegaptozero}
Let $r < C$ be positive integers and suppose $G$ is a connected graph with at least two nodes, consisting of one node of degree $C$ and with all other nodes of degree at most $r$. Let $\mu$ and $\kappa$ denote the largest and second-largest eigenvalue of the Laplacian $L_G$, respectively. Then the value $\kappa/ \mu$ goes to zero as $C/r$ goes to infinity, uniformly in all graphs $G$ satisfying the above conditions. 
\end{prop}

\begin{dem}
	See Appendix.
\end{dem}

\section{Bifurcations in diffusely coupled stable systems}
        
Using our results so far, we show what bifurcations to expect in diffusely coupled stable systems in $1$, $2$ or $3$ bifurcation parameters. Note that Theorem \ref{mainformal} tells us that the dynamics on the center manifold is conjugate to that of a reduced vector field $R: \mathbb{R}^m \times \Omega \rightarrow \mathbb{R}^m$, satisfying $R(0; \varepsilon) = 0$ for all $\varepsilon \in \Omega$. By Corollary \ref{obsaboutremainingeig2sdf} we may furthermore assume the linearization $D_xR(0;0)$ to be nilpotent with a one-dimensional kernel. Other than that, no restrictions apply to the Taylor expansion of $R$.
       
       Assuming that an $m$-parameter bifurcation can generically generate an $m$-dimensional generalized kernel, we each time consider $m$ parameter bifurcations for a system on $\mathbb{R}^m$. In Subsection \ref{onetwopar} we briefly investigate the cases $m=1$ and $m=2$.  Our main result is presented in Subsection \ref{threepar}, where we show the emergence of chaos for $m=3$. In most cases, the main difficulty lies in adapting existing results on generic unfoldings to the setting where $R(0, \varepsilon) = 0$ for all $\varepsilon \in \Omega$.

    \subsection{One and two parameters}\label{onetwopar}
        Motivated by our results so far, we describe the generic $m$ parameter bifurcations for systems $R$ on $\mathbb{R}^m$, where $m = 1,2$, subject to the condition $R(0;\varepsilon) = 0$ for all $\varepsilon \in \Omega$. We each time assume a nilpotent Jacobian with a  one-dimensional kernel. We start with the case $m=1$.
        
        \begin{obs}[The case $m = 1$]
        A map $R: \mathbb{R} \times \mathbb{R} \rightarrow \mathbb{R}$ satisfying $R(0;\varepsilon) = 0$ for all $\varepsilon \in \mathbb{R}$ and $D_xR(0;0) = 0$ has the general form 
        $$R(x;\varepsilon) = x(ax+b\varepsilon + \mathcal{O}((|x, \varepsilon)|^2)) \, , \text{ for } a, b \in \mathbb{R}\, .$$
        Under the generic assumption that $a,b \not= 0$, we find a transcritical bifurcation. Returning to the setting of our network system, this corresponds to a loss of stability of the fully synchronous solution.
        \end{obs}
        
         \begin{obs}[The case $m = 2$]
         Consider first a two-parameter vector field $R: \mathbb{R}^2 \times \mathbb{R}^2 \rightarrow \mathbb{R}^2$ satisfying $R(0;0) = 0$ and with non-zero nilpotent Jacobian $D_xR(0;0)$. Such a system generically displays a Bogdanov-Takens bifurcation. However, in this bifurcation scenario there are parameter values for which there is no fixed point. Hence, if we impose the additional condition $R(0;\varepsilon) = 0$ for all $\varepsilon \in \mathbb{R}^2$, then another (generic) bifurcation scenario has to occur. This latter situation is worked out in \cite{hirschberg1991unfolding}. The corresponding generic bifurcation involves multiple fixed points, heteroclinic as well as homoclinic connections, and periodic orbits. One striking feature is the presence of a homoclinic orbit from the origin, which is approached as a limit of stable periodic solutions. In our network setting, such a periodic solution means a cyclic time-evolution of the system from full synchrony to less synchrony and back. The time at which the system is indistinguishably close to full synchrony can moreover be made arbitrarily long.
         \end{obs}
        
	\subsection{Three parameters: chaotic behavior}\label{threepar}
	
	In this section we will prove Corollary \ref{chaos}, which allows us to conclude that chaotic behavior occurs in diffusely coupled stable systems. To this end, we will apply the theory developed so far. Before this, we will present a detailed background on how we will achieve chaos.
	
	We expect to find chaos in the network through the existence of a Shilnikov homoclinic orbit on a three-dimensional center manifold.
	
\begin{figure}[H]
	    \centering
    \tikzset{every picture/.style={line width=0.75pt}} 
    \begin{tikzpicture}[x=0.75pt,y=0.75pt,yscale=-1,xscale=1]

\draw [color={rgb, 255:red, 0; green, 0; blue, 0 }  ,draw opacity=1 ][line width=0.75]    (211.5,228.8) .. controls (211.67,161.67) and (209.67,109) .. (230.33,91) .. controls (251,73) and (286.33,89.67) .. (304.33,109) .. controls (322.33,128.33) and (331.67,161.67) .. (327,185.67) .. controls (322.38,209.43) and (317.11,226) .. (275.97,241.95) ;
\draw [shift={(274.71,242.43)}, rotate = 339.16] [color={rgb, 255:red, 0; green, 0; blue, 0 }  ,draw opacity=1 ][line width=0.75]    (10.93,-3.29) .. controls (6.95,-1.4) and (3.31,-0.3) .. (0,0) .. controls (3.31,0.3) and (6.95,1.4) .. (10.93,3.29)   ;
\draw  [color={rgb, 255:red, 0; green, 0; blue, 0 }  ,draw opacity=1 ][fill={rgb, 255:red, 0; green, 0; blue, 0 }  ,fill opacity=1 ] (209.24,228.8) .. controls (209.24,227.55) and (210.25,226.54) .. (211.5,226.54) .. controls (212.75,226.54) and (213.76,227.55) .. (213.76,228.8) .. controls (213.76,230.05) and (212.75,231.06) .. (211.5,231.06) .. controls (210.25,231.06) and (209.24,230.05) .. (209.24,228.8) -- cycle ;
\draw  [line width=0.75]  (213.72,228.8) .. controls (214.55,228.69) and (215.12,228.49) .. (215.4,228.2) .. controls (215.68,227.89) and (215.53,227.57) .. (214.89,227.24) .. controls (214.18,226.88) and (213.05,226.6) .. (211.5,226.4) .. controls (209.77,226.18) and (207.88,226.11) .. (205.83,226.2) .. controls (203.55,226.29) and (201.51,226.56) .. (199.71,226.99) .. controls (197.74,227.47) and (196.35,228.08) .. (195.61,228.8) .. controls (194.8,229.59) and (194.84,230.39) .. (195.76,231.21) .. controls (196.76,232.09) and (198.61,232.86) .. (201.28,233.5) .. controls (204.15,234.19) and (207.54,234.63) .. (211.5,234.83) .. controls (215.68,235.04) and (219.87,234.95) .. (224,234.54) .. controls (228.39,234.12) and (232.08,233.41) .. (235.13,232.42) .. controls (238.34,231.37) and (240.3,230.17) .. (241.06,228.8) .. controls (241.84,227.36) and (241.2,225.96) .. (239.08,224.58) .. controls (236.87,223.13) and (233.36,221.93) .. (228.56,220.96) .. controls (223.56,219.95) and (217.87,219.35) .. (211.5,219.14) .. controls (204.87,218.93) and (198.42,219.19) .. (192.16,219.91) .. controls (185.68,220.67) and (180.29,221.82) .. (176.03,223.37) .. controls (171.62,224.98) and (169.02,226.79) .. (168.27,228.8) .. controls (167.47,230.89) and (168.75,232.9) .. (172.08,234.84) .. controls (175.54,236.84) and (180.7,238.49) .. (187.6,239.78) .. controls (194.73,241.11) and (202.68,241.88) .. (211.5,242.08) .. controls (220.58,242.29) and (229.29,241.87) .. (237.67,240.83) .. controls (246.29,239.75) and (253.34,238.16) .. (258.81,236.05) .. controls (264.42,233.88) and (267.63,231.46) .. (268.41,228.8) .. controls (269.21,226.07) and (267.33,223.45) .. (262.76,220.95) .. controls (258.09,218.39) and (251.23,216.29) .. (242.23,214.68) .. controls (232.99,213.02) and (222.77,212.09) .. (211.5,211.89) .. controls (199.97,211.68) and (188.96,212.26) .. (178.49,213.63) .. controls (167.77,215.03) and (159.07,217.07) .. (152.35,219.74) .. controls (145.5,222.47) and (141.69,225.49) .. (140.92,228.8) .. controls (140.14,232.18) and (142.63,235.4) .. (148.4,238.46) .. controls (154.32,241.59) and (162.8,244.12) .. (173.93,246.06) .. controls (185.26,248.04) and (197.82,249.13) .. (211.5,249.34) .. controls (225.48,249.55) and (238.75,248.8) .. (251.35,247.11) .. controls (260.5,245.88) and (268.39,244.26) .. (275.04,242.27) ;
\end{tikzpicture}
	    \caption{Shilnikov\textquotesingle s homoclinic orbit.}
	    \label{fig:shil}
	\end{figure}
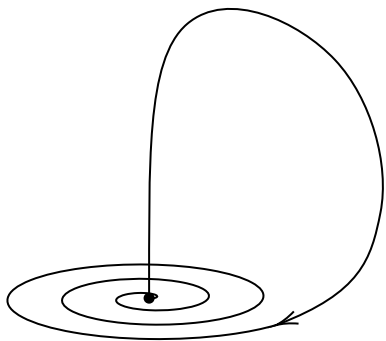

	The Shilnikov configuration can be seen as a combination of linear and nonlinear behavior involving a saddle fixed point. A two-dimensional stable manifold attracts trajectories exponentially fast to the fixed point, where the eigenvalues of the linearization are  $\lambda_{1,2}=-\alpha\pm i\beta$ with $\alpha>0$ and $\beta \not= 0$. Transversal to this there is a one-dimensional unstable manifold repelling away trajectories with real eigenvalue $\gamma>0.$ The Shilnikov homoclinic orbit emerges from the re-injection of the one-dimensional unstable manifold into the two-dimensional stable manifold, see Figure \ref{fig:shil}. Of course this re-injection is a consequence of nonlinear terms. L. P. Shilnikov proved that if $\gamma>\alpha$, there are countably many saddle periodic orbits in a neighborhood of the homoclinic orbit. The proof consists of showing topological equivalence between a Poincar\'e map and the shift map of two symbols. The existence of chaotic behavior is in the sense that Robert L. Devaney defined for deterministic systems, with strongly sensitive dependence on initial conditions, topological transitivity and dense periodic points.

	We next give a brief summary of results contained in the paper \cite{ibanezconfig}. 
	The authors studied the three parameter unfolding of nonlinear vector fields on $\mathbb{R}^3$ with linear part conjugate to a nilpotent singularity of codimension three. After making several coordinate changes, the following normal form is presented:
	\begin{equation}\label{eq:83}
	    y\frac{\partial}{\partial x}+z\frac{\partial}{\partial y}+(\lambda -y+\nu z-\frac{x^{2}}{2}+\mathcal{O}(\kappa))\frac{\partial}{\partial z}
	\end{equation}
	where the parameters are given by $\tau=(\lambda,\nu,\kappa).$ The  parameter $\kappa$ is introduced by means of a blow-up technique, and the term
	\begin{equation*}
	    y\frac{\partial}{\partial x}+z\frac{\partial}{\partial y}
	\end{equation*}
	denotes the nilpotent singularity of codimension three on $\mathbb{R}^{3}$.  
	Equation (\ref{eq:83}) has two hyperbolic fixed points for $\lambda>0$ and $\nu =0$, namely $p_{1}=(-\sqrt{2\lambda},0,0)$ with local behavior given by a two-dimensional stable and one-dimensional unstable manifold and $p_{2}=(+\sqrt{2\lambda},0,0)$ with local behavior given by a two-dimensional unstable and one-dimensional stable manifold. 
	Knowing there is a solution $x(t)$ for a specific positive parameter $\lambda=\lambda^{*}$ such that $x(t)\rightarrow p_{1}$ as $t\rightarrow-\infty$ and $x(t)\rightarrow p_{2}$ as $t\rightarrow+\infty$, the authors proved analytically the existence of another solution, also for the parameter $\lambda^{*},$ connecting both two-dimensional stable and unstable manifolds and thus forming another  heteroclinic orbit. 
	Theorem 4.1 of \cite{ibanezconfig} states that in any neighborhood of the parameter $\tau_{0}=(\lambda_{0},\nu_{0},\kappa_{0})=(\lambda^{*},0,0)$, where the heteroclinic orbits exist, there are parameters $\tau=(\lambda,\nu,\kappa)$ such that the heteroclinic orbit breaks and a Shilnikov homoclinic orbit appears.
	
	For completeness, we state Theorem 4.1 below in a slightly altered form.
	
	\begin{teo}[Theorem 4.1 \cite{ibanezconfig}.]\label{ibanez}
	In every neighborhood of the parameter $\tau_{0}=(\lambda_{0},\nu_{0},\kappa_{0})=(\lambda^{*},0,0)$ there exist parameter values $\tau=(\lambda,\nu,\kappa)$ such that the equation
	\begin{equation*}
	    y\frac{\partial}{\partial x}+z\frac{\partial}{\partial y}+\left(\lambda-y+\nu z-\frac{x^{2}}{2}+\mathcal{O}(\kappa)\right)\frac{\partial}{\partial z}
	\end{equation*}
	has a homoclinic orbit given by the intersection of the two-dimensional stable and one-dimensional unstable invariant manifolds at the hyperbolic fixed point $p_{1}.$
	\end{teo}
	
	As was the case for $m=2$, we cannot immediately use this existing result, as the parameter dependent systems on the center manifold of our network ODE satisfy $R(0,\varepsilon) = 0$ for all $\varepsilon \in \mathbb{R}^3$. 
	It remains to show that with this existing restriction, we may still reduce our system to the family given by Equation \eqref{eq:83}. 
	This then proves Corollary \ref{chaos} as a consequence of Theorem \ref{mainformal}.

We therefore start with a parameterized vector field $R: \mathbb{R}^3 \times \mathbb{R}^3 \rightarrow \mathbb{R}^3$ satisfying $R(0,\varepsilon) = 0$ for all $\varepsilon \in \mathbb{R}^3$. After a linear coordinate change, we may assume the Jacobian $D_xR(0;0)$ to be given by
	\begin{equation*}
	    J=\left(\begin{array}{ccc}
	         0&1&0\\
	         0&0&1\\
	         0&0&0
	    \end{array}\right).
	\end{equation*}

\noindent	We thus get the system
		 \begin{equation}\label{redunf}
		     \begin{array}{ccccc}
		     \dot{x}_{1}&=&x_{2}&+&h_{1}(x;\varepsilon)\\
		     \dot{x}_{2}&=&x_{3}&+&h_{2}(x;\varepsilon)\\
		     \dot{x}_{3}&=& & &h_{3}(x;\varepsilon) 
		 \end{array}, 
		 \end{equation}
		 with $x=(x_1, x_2, x_3)$, and where $h_{1},h_{2},h_{3}$ are the higher order (nonlinear) terms of $R(x;\varepsilon)$. Note that we have  $h_{1}(0;\varepsilon)=h_{2}(0;\varepsilon)=h_{3}(0;\varepsilon)=0$ for all $\varepsilon\in\Omega$, and $Dh_{1}(0;0)=Dh_{2}(0;0)=Dh_{3}(0;0)=0$.
		 
		 
		 To bring our system in the form (\ref{eq:83}), we will proceed as in the paper \cite{ibanezconfig}. 
		 The first step is to  get rid of the nonlinear terms $h_{1}$ and $h_{2}$ by means of a coordinate transformation.
		 To eliminate $h_{1}$ we consider the following change of coordinates
		 \begin{eqnarray}\label{transform_1}
		 y_{1}&=&x_{1}\nonumber \\ 
		 y_{2}&=&x_{2}+h_{1}(x;\varepsilon)\\
		 y_{3}&=&x_{3} \, . \nonumber
		 \end{eqnarray}
		 Applying it to (\ref{redunf}), we get
		 \begin{eqnarray*}
		 \dot{y}_{1}&=&\dot{x}_{1}=x_{2}+h_{1}(x;\varepsilon)=y_{2}\\
		 \dot{y}_{2}&=&\dot{x}_{2}+D_{x}h_{1}(x;\varepsilon)\dot{x} \\
		 &=& x_{3}+h_{2}(x;\varepsilon)+D_{x_{1}}h_{1}(x;\varepsilon)(x_2 + h_1(x;\varepsilon)) + D_{x_{2}}h_{1}(x;\varepsilon)(x_3 + h_2(x;\varepsilon)) + D_{x_{3}}h_{1}(x;\varepsilon) h_3(x;\varepsilon)\\
		 &=& y_{3}+h_{2}(x;\varepsilon)+D_{x_{1}}h_{1}(x;\varepsilon)(x_2 + h_1(x;\varepsilon)) + D_{x_{2}}h_{1}(x;\varepsilon)(x_3 + h_2(x;\varepsilon)) + D_{x_{3}}h_{1}(x;\varepsilon) h_3(x;\varepsilon)\\
		 \dot{y}_{3}&=&\dot{x}_{3}=h_{3}(x;\varepsilon)
		 \end{eqnarray*}
		 We therefore get the new system
		 \begin{eqnarray}\label{eq1ch8}
		 \nonumber\dot{y}_{1}&=&y_{2}\\
		 \dot{y}_{2}&=&y_{3}+\tilde{h}_{2}(y;\varepsilon)\\
		 \nonumber\dot{y}_{3}&=&\tilde{h}_{3}(y;\varepsilon),
		 \end{eqnarray}
		 with $y = (y_1, y_2, y_3)$, and where $\tilde{h}_{2}$ and $\tilde{h}_{3}$ are uniquely defined by the relations 
		 \begin{align}
	    &\tilde{h}_{2}(y; \varepsilon) = 	     {h}_{2}(x;\varepsilon)+D_{x_{1}}h_{1}(x;\varepsilon)(x_2 + h_1(x;\varepsilon)) + D_{x_{2}}h_{1}(x;\varepsilon)(x_3 + h_2(x;\varepsilon)) + D_{x_{3}}h_{1}(x;\varepsilon) h_3(x;\varepsilon) \\ \nonumber
	    &\tilde{h}_{3}(y;\varepsilon) = h_{3}(x;\varepsilon)\, .
		 \end{align}
		 Note that $\tilde{h}_{2}$ and $\tilde{h}_{3}$ again have vanishing linear terms, and moreover satisfy $\tilde{h}_{2}(0;\varepsilon) = \tilde{h}_{3}(0;\varepsilon) = 0$ for all $\varepsilon$.
		 
		 To eliminate $\tilde{h}_2$ we consider the change of coordinates
		 \begin{eqnarray}\label{transform_2}
		 z_{1}&=&y_{1}\nonumber \\
		 z_{2}&=&y_{2}\\ \nonumber
		 z_{3}&=&y_{3}+\tilde{h}_{2}(y;\varepsilon).
		 \end{eqnarray}
		 Applying it to (\ref{eq1ch8}), we obtain
		 \begin{eqnarray}
		 \dot{z}_{1}&=&\dot{y}_{1}=y_{2}=z_{2}\nonumber \\
		 \dot{z}_{2}&=&\dot{y}_{2}=y_{3}+\tilde{h}_{2}(y;\varepsilon)=z_{3}\\ \nonumber
		 \dot{z}_{3}&=&\dot{y}_{3}+D_{y}\tilde{h}_{2}(y;\varepsilon)\dot{y}\\ \nonumber &=& \tilde{h}_{3}(y;\varepsilon)+D_{y_{1}}\tilde{h}_{2}(y;\varepsilon)y_2 + D_{y_{2}}\tilde{h}_{2}(y;\varepsilon)(y_3 + \tilde{h}_2(y;\varepsilon)) + D_{y_{3}}\tilde{h}_{2}(y;\varepsilon) \tilde{h}_3(y;\varepsilon)\, .
		 \end{eqnarray}
		 We thus get the new system
		 \begin{eqnarray}
		 \dot{z}_{1}&=&z_{2} \nonumber \\
		 \dot{z}_{2}&=&z_{3}\\ \nonumber
		 \dot{z}_{3}&=&\hat{h}_{3}(z;\varepsilon)
		 \end{eqnarray}
		 with $z = (z_1, z_2, z_3)$, and where $\hat{h}_{3}(z;\varepsilon)$ is locally defined by
		 \begin{equation}
		     \hat{h}_{3}(z;\varepsilon) = \tilde{h}_{3}(y;\varepsilon)+D_{y_{1}}\tilde{h}_{2}(y;\varepsilon)y_2 + D_{y_{2}}\tilde{h}_{2}(y;\varepsilon)(y_3 + \tilde{h}_2(y;\varepsilon)) + D_{y_{3}}\tilde{h}_{2}(y;\varepsilon) \tilde{h}_3(y;\varepsilon) \, .
		 \end{equation}
		 Note that  $\hat{h}_{3}$ again has no linear terms and satisfies  $\hat{h}_{3}(0;\varepsilon) = 0$ for all $\varepsilon$. 
		 Moreover, in case of $h_1 = h_2 = 0$ we would find $x=y=z$ and $h_3 = \hat{h}_3$, which shows that no other restrictions apply to $\hat{h}_3$.
		 Writing 
		 \begin{equation}
		     \hat{h}_{3}(z;\varepsilon) = \mathcal{E}_1(\varepsilon)z_1 + \mathcal{E}_2(\varepsilon)z_2 + \mathcal{E}_3(\varepsilon)z_3 + \mathcal{O}(\|z\|^2) 
		 \end{equation}
		 for some locally defined $\mathcal{E}_i: \mathbb{R}^3 \to \mathbb{R}$, we therefore see that generically we may redefine $\varepsilon = (\varepsilon_1, \varepsilon_2, \varepsilon_3)$ so that
			 \begin{equation}
		     \hat{h}_{3}(z;\varepsilon) = \varepsilon_{1}z_{1}+\varepsilon_{2}z_{2}+\varepsilon_{3}z_{3} + \mathcal{O}(\|z\|^2) \, .
		 \end{equation}	 
		 It follows that we may write
		 \begin{eqnarray}\label{eq2}
		 \nonumber\dot{z}_{1}&=&z_{2}\\
		 \dot{z}_{2}&=&z_{3}\\
		 \nonumber\dot{z}_{3}&=&\varepsilon_{1}z_{1}+\varepsilon_{2}z_{2}+\varepsilon_{3}z_{3}+a_{1}z^{2}_{1}+a_{2}z_{2}^{2}+a_{3}z_{3}^{2}\\
		 \nonumber&+&a_{4}z_{1}z_{2}+a_{5}z_{1}z_{3}+a_{6}z_{2}z_{3}+\mathcal{O}(\|z\|^3+\|\varepsilon\|\|z\|^2),
		 \end{eqnarray}
		 for some coefficients $a_1, \dots, a_6 \in \mathbb{R}$. We will assume that $a_{1}\neq0,$ so that the system (\ref{eq2}) locally has two branches of steady states: $z(\varepsilon) = (0,0,0)$ and
		 $z(\varepsilon) = (\tilde{z}_1(\varepsilon),0,0)$. 
		 A straightforward calculation shows that
		 \begin{equation}
		     \tilde{z}_1(\varepsilon) = -\frac{\varepsilon_{1}}{a_{1}} +\mathcal{O}{(\|\varepsilon\|^2)}\, .
		 \end{equation}
		 Motivated by this, we perform the change of coordinates:
		 \begin{eqnarray}
		 w_{1}&=&z_{1}+\frac{\varepsilon_{1}}{2a_{1}}\nonumber \\ 
		 w_{2}&=&z_{2}\\ \nonumber
		 w_{3}&=&z_{3}\,
		 \end{eqnarray}
		 Applying this to (\ref{eq2}), we get
		 \begin{eqnarray}
		 \dot{w}_{1}&=&\dot{z}_{1}=z_{2}=w_{2}\\ \nonumber
		 \dot{w}_{2}&=&\dot{z}_{2}=z_{3}=w_{3}\\ \nonumber
		 \dot{w}_{3}&=&\dot{z}_{3}=\varepsilon_{1}(w_{1}-\frac{\varepsilon_{1}}{2a_{1}})+\varepsilon_{2}w_{2}+\varepsilon_{3}w_{3}
		 +a_{1}(w_{1}-\frac{\varepsilon_{1}}{2a_{1}})^{2}+a_{2}w_{2}^{2}+a_{3}w_{3}^{2}\\ \nonumber
		 &+&a_{4}(w_{1}-\frac{\varepsilon_{1}}{2a_{1}})w_{2}+a_{5}(w_{1}-\frac{\varepsilon_{1}}{2a_{1}})w_{3}+a_{6}w_{2}w_{3}+\mathcal{O}(|z|^3+|\varepsilon||z|^2) \\ \nonumber
		 &=&-\frac{\varepsilon^{2}_{1}}{4a_{1}}+\varepsilon_{2}w_{2}+\varepsilon_{3}w_{3}+a_{1}w^{2}_{1}+a_{2}w_{2}^{2}+a_{3}w_{3}^{2}\\
		 &+&a_{4}w_{1}w_{2}-\frac{a_{4}\varepsilon_{1}w_{2}}{2a_{1}}+a_{5}w_{1}w_{3}-\frac{a_{5}\varepsilon_{1}w_{3}}{2a_{1}}+a_{6}w_{2}w_{3}+\mathcal{O}(|z|^3+|\varepsilon||z|^2)\, . \nonumber
		 \end{eqnarray}
		 We have left the remainder term $\mathcal{O}(|z|^3+|\varepsilon||z|^2)$ as is, which will benefit us later. Rearranging terms, we get the new system
		 \begin{eqnarray}\label{equationsinwaa}
		 \dot{w}_{1}&=&w_{2}\\ \nonumber
		 \dot{w}_{2}&=&w_{3}\\ \nonumber
		 \dot{w}_{3}&=&-\frac{\varepsilon^{2}_{1}}{4a_{1}}+(\varepsilon_{2}-\frac{a_{4}\varepsilon_{1}}{2a_{1}})w_{2}+(\varepsilon_{3}-\frac{a_{5}\varepsilon_{1}}{2a_{1}})w_{3}+a_{1}w^{2}_{1}+a_{2}w_{2}^{2}+a_{3}w_{3}^{2}\\ \nonumber
		 &+&a_{4}w_{1}w_{2}+a_{5}w_{1}w_{3}+a_{6}w_{2}w_{3}+\mathcal{O}(|z|^3+|\varepsilon||z|^2)\, . \nonumber
		 \end{eqnarray}
		 Similar to the paper \cite{ibanezconfig}, we now introduce a blow-up parameter $\kappa\in\mathbb{R}$ and write
		 \begin{eqnarray}
		     w_{1}&=&\kappa^{3}u_{1}\qquad \varepsilon_{1}=\kappa^{3}\gamma_{1}\\ \nonumber
		     w_{2}&=&\kappa^{4}u_{2}\qquad \varepsilon_{2}=\kappa^{2}\gamma_{2}\qquad \bar{t} =\kappa t\\ \nonumber
		     w_{3}&=&\kappa^{5}u_{3}\qquad \varepsilon_{3}=\kappa\gamma_{3}.
		 \end{eqnarray}
		 Note that we get 
		 \begin{align}
		     z_1 &= w_1 - \frac{\varepsilon_1}{2a_1} = \kappa^3\left(u_1 - \frac{\gamma_1}{2a_1}\right) \\ \nonumber 
		     z_2 &= w_2 = \kappa^4u_2 \\ \nonumber
		     z_3 &= w_3 = \kappa^5u_3 \, ,
		 \end{align}
		 so that we may write $\|z\| = \mathcal{O}(\kappa^3)$. Applying it to Equation (\ref{equationsinwaa}), we get
		 \begin{eqnarray*}
		 \frac{du_{1}}{d\bar{t}}&=&\frac{1}{\kappa^{3}}\frac{dw_{1}}{d\bar{t}}=\frac{1}{\kappa^{3}}\frac{dw_{1}}{\kappa dt}=\frac{1}{\kappa ^{4}}\frac{dw_{1}}{dt}=\frac{1}{\kappa ^{4}}w_{2}=u_{2}\\
		 \frac{du_{2}}{d\bar{t}}&=&\frac{1}{\kappa ^{4}}\frac{dw_{2}}{d\bar{t}}=\frac{1}{\kappa ^{4}}\frac{dw_{2}}{\kappa dt}=\frac{1}{\kappa ^{5}}\frac{dw_{2}}{dt}=\frac{1}{\kappa ^{5}}w_{3}=u_{3}\\
		 \frac{du_{3}}{d\bar{t}}&=&\frac{1}{\kappa ^{5}}\frac{dw_{3}}{d\bar{t}}=\frac{1}{\kappa ^{5}}\frac{dw_{3}}{\kappa dt}=\frac{1}{\kappa ^{6}}\frac{dw_{3}}{dt}\, ,
		 \end{eqnarray*}
		 where furthermore
		 \begin{eqnarray*}
		     \frac{dw_{3}}{dt} = \dot{w}_3 = \kappa ^{6}\left(-\frac{\gamma_{1}^2}{4a_1}+\gamma_{2}u_{2}+\gamma_{3}u_{3}+a_{1} u_{1}^{2}\right)+ \mathcal{O}(\kappa^7)\, .
		 \end{eqnarray*}
		 Summarizing, we find
		 \begin{eqnarray}\label{eq4}
		 \nonumber\frac{du_{1}}{d\bar{t}}&=&u_{2}\\
		\frac{du_{2}}{d\bar{t}}&=&u_{3}\\
		 \nonumber\frac{du_{3}}{d\bar{t}}&=&-\frac{\gamma_{1}^2}{4a_1}+\gamma_{2}u_{2}+\gamma_{3}u_{3}+a_{1} u_{1}^{2}+\mathcal{O}(\kappa).
		 \end{eqnarray}
		 We next focus on the parameter $\gamma_{2}$. We assume from here on out that $\gamma_{2}<0$ and perform the following change of coordinates
		 \begin{eqnarray}\label{trabsform_3}
		 v_1&=&-2r^{3}u_1 \qquad v_2=-2r^{4}u_2 \\ \nonumber
		 v_3&=&-2r^{5}u_3  \qquad \tau=r^{-1}\bar{t}  \, ,
		 \end{eqnarray}
		 where 
		 \begin{equation}
		     r=\left(-\frac{1}{\gamma_{2}}\right)^{\frac{1}{2}} > 0\, .
		 \end{equation}
Applying it to (\ref{eq4}), we get
		 \begin{eqnarray}
		     \frac{dv_{1}}{d\tau}&=&(-2r^{3})\frac{du_{1}}{d\tau}=(-2r^{4})\frac{du_{1}}{d\bar{t}}=-2r^{4}u_{2}=v_{2}\\ \nonumber
		     \frac{dv_{2}}{d\tau}&=&(-2r^{4})\frac{du_{2}}{d\tau}=(-2r^{5})\frac{du_{2}}{d\bar{t}}=-2r^{5}u_{3}=v_{3}\\ \nonumber
		      \frac{dv_{3}}{d\tau}&=&(-2r^{5})\frac{du_{3}}{d\tau}=(-2r^{6})\frac{du_{3}}{d\bar{t}}\\ \nonumber
		      &=&(-2r^{6})\left(-\frac{\gamma_{1}^2}{4a_1}+\gamma_{2}u_{2}+\gamma_{3}u_{3}+a_{1} u_{1}^{2}+\mathcal{O}(\kappa)\right)\\ \nonumber
		      &=&(-2r^{6})\left(-\frac{\gamma_{1}^2}{4a_1}+\gamma_{2}\left(\frac{v_{2}}{-2r^4}\right)+\gamma_{3}\left(\frac{v_{3}}{-2r^5}\right)+a_{1} \left(\frac{v_{1}}{-2r^3}\right)^{2}+\mathcal{O}(\kappa)\right)\\ \nonumber
		      &=&\frac{\gamma_{1}^2r^6}{2a_1}+\gamma_{2}v_{2}r^2+\gamma_{3}v_{3}r-a_{1} \frac{v_{1}^2}{2}+\mathcal{O}(\kappa)\\ \nonumber
		       &=&\frac{\gamma_{1}^2r^6}{2a_1}-v_{2}+\gamma_{3}v_{3}r-a_{1} \frac{v_{1}^2}{2}+\mathcal{O}(\kappa)\, .
		 \end{eqnarray}
		 We thus get the new system
		 \begin{eqnarray}\label{eq5}
		    \nonumber v_1' := \frac{dv_{1}}{d\tau} &=& v_2\\
		     v_2' :=\frac{dv_{2}}{d\tau} &=& v_3 \\ \nonumber
		    v_3' :=\frac{dv_{3}}{d\tau} &=&\frac{\gamma_{1}^2r^6}{2a_1}-v_{2}+\gamma_{3}v_{3}r-a_{1} \frac{v_{1}^2}{2}+\mathcal{O}(\kappa)\, .
		 \end{eqnarray}
		 Finally, we make the following change of coordinates:
		 \begin{eqnarray}
		     {\rm{x}}&=&a_{1}v_{1}\nonumber \\ 
		     {\rm{y}}&=&a_{1}v_{2}\\ \nonumber
		     {\rm{z}}&=&a_{1}v_{3}\, .
		 \end{eqnarray}
		 This gives
		 \begin{eqnarray}
		     \rm{x}'&=&a_{1}v_{1}'=a_{1}v_{2} = {\rm{y}}\\ \nonumber
		     \rm{y}'&=&a_{1}v_{2}'=a_{1}v_{3} = {\rm{z}}\\ \nonumber
		     \rm{z}'&=&a_{1}v_{3}'= \frac{\gamma_{1}^2r^6}{2}-a_1v_{2}+\gamma_{3}a_1v_{3}r-a_{1}^2 \frac{v_{1}^2}{2}+\mathcal{O}(\kappa) \\ \nonumber
		     &=&\frac{\gamma_{1}^2r^6}{2}-{\rm{y}}+\gamma_{3}r{\rm{z}}- \frac{{\rm{x}}^2}{2}+\mathcal{O}(\kappa)\, .
		 \end{eqnarray}
		Setting $\lambda := \frac{\gamma_{1}^2r^6}{2}$ and $\nu:= \gamma_{3}r$, we arrive at the vector field
        \begin{eqnarray}
		     ={\rm{y}}\frac{\partial}{\partial {\rm{x}}}+{\rm{z}}\frac{\partial}{\partial {\rm{y}}}+\left(
		     \lambda-{\rm{y}}+\nu{\rm{z}}- \frac{{\rm{x}}^2}{2}+\mathcal{O}(\kappa)\right)\frac{\partial}{\partial {\rm{z}}}
		\end{eqnarray}
		from Theorem \ref{ibanez}. Note that $\lambda = \frac{\gamma_{1}^2r^6}{2}$ is necessarily non-negative. However, this may be assumed in the setting of Theorem \ref{ibanez}, as $\lambda^{*}>0$. This theorem thus predicts chaos in the setting of our coupled cell system, provided $m=3$ and the network in question is at least $2$-versatile.
		
		
\section*{Acknowledgments}
We thank  Edmilson Roque and Jeroen Lamb for enlightening discussions. TP was  supported in part by FAPESP Cemeai Grant No. 2013/07375-0 and is a Newton Advanced Fellow of the Royal Society NAF$\backslash$R1$\backslash$180236. EN was supported by FAPESP grant 2020/01100-2. TP and EN 
were partially supported by Serrapilheira Institute (Grant No. Serra-1709-16124). FCQ was supported by CAPES. DT was supported by Leverhulme Trust grant RPG-2021-072.

 \section{Appendix: Examples}	
 	
		\subsection{Versatile graphs by means of the complement graph}\label{means_of_the_complement_graph}
    We next introduce a method for generating $\rho$-versatile graphs, for any $\rho \in \mathbb{N}$. Our construction involves the definition of the complement graph, given below.
    
	
	\begin{defn}
		Given an undirected graph $G$, we define the complement graph $G^{\circ}$ as the graph obtained from $G$ by leaving out all existing edges and adding all edges between distinct vertices that were not there in $G.$ 
	\end{defn}
	
    Since our graphs don't have self-loop, we have that $G^{\circ\circ} = G.$
	
	\begin{teo}\label{versatilefromcomp}
		Let $G$ be a graph consisting of precisely two disconnected components of different order. Then $G^{\circ}$ is a connected graph whose Laplacian has a simple, largest eigenvalue whose eigenvector $v$ satisfies $\sum_{i=1}^{|G^{\circ}|} \nu_i^\ell \not= 0$ for all $\ell>1.$
Moreover, suppose the two disconnected components of $G$ have number of vertices $s$ and $t$. Then the largest eigenvalue of the Laplacian of $G^{\circ}$ is equal to $s+t$ and a corresponding eigenvector is given by 
		$$(\underbrace{t,t, \dots, t}_{s \text{ times}}, \underbrace{-s, -s, \dots, -s}_{t \text{ times}})\, .$$ 
		Here the entries are ordered so that the vertices of the first component of $G$ (which has $s$ vertices) are enumerated first, after which those of the second component of $G$ (which has $t$ vertices) are listed.
	\end{teo}
	
    \begin{dem}
    The proof uses a result that relates the eigenvalues and eigenvectors of the Laplacian of a graph to those of the Laplacian of its complement graph. This result is known,  
   but incorporated here for completeness.   Suppose the two components of $G$ have $s\neq 0$ and $t\neq 0$ vertices, where $s\neq t$ and $|G|=t+s=N$.

Recall that 
%
%
    the dimension of the kernel of $L_{G}$ equals the number of connected components of $G$, we see that $\lambda_{2}=\lambda_{1} = 0$.  Thus, 
    \begin{equation}
        \Span(v_1, v_2) = \{(\nu_1,\dots,\nu_{N}) \in \mathbb{R}^N\mid \nu_1 = \dots = \nu_s, \nu_{s+1} = \dots = \nu_t\} \, ,
    \end{equation}
 where we grouped the entries according to the connected components, we see that may choose
    \begin{align}
        v_1 = (1, \dots, 1) \, \text{ and } \, v_2 = (\underbrace{t,t, \dots, t}_{s \text{ times}}, \underbrace{-s, -s, \dots, -s}_{t \text{ times}})\, 
    \end{align}
    which we assume from here on out. Note that indeed $v_1 \perp v_2$. Let $L_{G^{\circ}}$ be the Laplacian matrix associated to the complement graph $G^{\circ}$. We note that we have the identity
    \begin{equation}\label{63}
        L_{G}+L_{G^{\circ}}=
        N\cdot \id-E \, ,
    \end{equation}
    where $E$ is a matrix where every element equals $1$. 
    As we have $v_1 \perp v_{i}$ for all $i=2,\dots,N$, it follows that $Ev_i = 0$ for all $i=2,\dots,N$.
    From Equation (\ref{63}) we get
    \begin{equation}
        L_{G^{\circ}}=-L_{G}+N\cdot \id-E
    \end{equation}
    and evaluating at the eigenvectors $v_{i}$ for $i=2,\dots,N$ gives
    \begin{eqnarray}
        \nonumber L_{G^{\circ}}v_{i}&=&-L_{G}v_{i}+N\cdot\id v_{i}-E v_{i}\\
        &=&(N-\lambda_{i})v_{i}.
    \end{eqnarray}
    Thus, for each $i=2,\dots,N$ we find that $(N-\lambda_{i})$ is an eigenvalue of $L_{G^{\circ}}$, with a corresponding eigenvector given by $v_i$. As we also have $L_{G^{\circ}}v_1 = 0$, we see that the spectrum of $L_{G^{\circ}}$ is given by
    \begin{equation}
        N-\lambda_{2}\geq N-\lambda_{3}\geq\cdots\geq N-\lambda_{N} \geq \lambda_{1}=0
    \end{equation}
    The largest eigenvalue of $L_{G^{\circ}}$ is therefore equal to $N-0 = N = s+t$, with an eigenvector given by
    \begin{equation}
        v = v_2 =(s,s,\dots,s,-t,-t,\dots,-t)\, .
    \end{equation}
Since $\lambda_3 > 0$, so that the eigenvalue $N$ is  simple. 
    Using that $st\neq 0$ and $s\neq t$, we find for all $\ell>1$
    \begin{equation}
        \sum_{i=1}^{N}\nu_{i}^{\ell}=\sum_{i=1}^{t}s^{\ell}\pm\sum_{i=1}^{s}t^{\ell}=t(s^{\ell})\pm s(t^{\ell})=st(s^{\ell-1}\pm t^{\ell-1})\neq 0\, .
    \end{equation}

    Finally, we argue that $G^{\circ}$ is a connected graph. Indeed, if $x,y\in G$ are in different connected components, then they share an edge in $G^{\circ}$ by definition of this latter graph. If on the other hand $x$ and $y$ are in the same component of $G$, then in $G^{\circ}$ they both share an edge with some node $z$ from the other component of $G$. This completes the proof.
    \end{dem}

		\begin{ex}
		Let $G=(V,E)$ be the undirected graph with $V=\{1,2,3\}$ and two disconnected components of a different order $s=1$ and $t=2$, shown in Figure \ref{f1}. Then $G^{\circ}$ is a connected, non-regular graph with
		\begin{equation}
		\nonumber L_{G^{\circ}}=\left(\begin{array}{ccc}
		1&0&-1\\
		0&1&-1\\
		-1&-1&2
		\end{array}\right)_{3\times 3}.
		\end{equation}
		
		We have $\spec(L_{G^{\circ}})=\{3,1,0\}$ with \textbf{simple} and \textbf{largest eigenvalue} $\lambda=s+t=3$ whose corresponding eigenvector $v=(1,1,-2)$ satisfies $\sum_{i=1}^{3}\nu_{i}^{\ell}\neq0$ for all $\ell>1.$ 
		\end{ex}
	\begin{ex}
		Let $G=(V,E)$ be the undirected graph with $V=\{1,2,3,4,5\}$ and two disconnected components of a different order $s=2$ and $t=3,$ shown in Figure \ref{f2}. Then $G^{\circ}$ is a connected graph with
		\begin{equation}
			\nonumber L_{G^{\circ}}=\left(\begin{array}{ccccc}
				3&0&-1&-1&-1\\
				0&2&0&-1&-1\\
				-1&0&3&-1&-1\\
				-1&-1&-1&3&0\\
				-1&-1&-1&0&3
			\end{array}\right)_{5\times 5}.
		\end{equation}
		
		Here $\spec(L_{G^{\circ}})=\{5,4,3,2,0\}$ with \textbf{simple} and \textbf{largest eigenvalue} $\lambda=s+t=5$. Its corresponding eigenvector $v=(2,2,2,-3,-3)$ satisfies $\sum_{i=1}^{5}\nu_{i}^{\ell}\neq0$ for all $\ell>1.$ 
	\end{ex}
	
	Example \ref{exStar} below shows that the standard star graphs are $\rho$-versatile for any $\rho > 0$. These graphs consist of a single hub-node connected to all other nodes, which in turn have  degree $1$, shown in Figure \ref{f3}. We will explore the $\rho$-versatility of more general star graphs in Subsection \ref{stargraphs}.
	
	\begin{ex}[Star graphs]\label{exStar}
		Let $G=(V,E)$ be the undirected graph with $V=\{1,\dots,N+1\}$ and two disconnected components of order $s=1$ and $t=N,$ shown in Figure \ref{f3}. If the largest component of $G$ is complete, then $G^{\circ}$ is a connected graph with Laplacian matrix given by 
		\begin{equation}
			\nonumber L_{G^{\circ}}=\left(\begin{array}{cccc}
				N+1&-1&\cdots&-1\\
				-1&1&\cdots&0\\
				\vdots&\vdots&\ddots&\vdots\\
				-1&0&\cdots&1
			\end{array}\right)_{N+1\times N+1}.
		\end{equation}
		
		The spectrum $\spec(L_{G^{\circ}})=\{N+1,1,\dots,1,0\}$ has one \textbf{simple} and \textbf{largest eigenvalue} $\lambda=s+t=N+1.$ The  corresponding eigenvector is given by $v=(N+1,-1,\dots,-1)$ which satisfies the property $\sum_{i=1}^{N}\nu_{i}^{\ell}\neq0$ for all powers $\ell>1$. We can generate more examples of graphs $G^{\circ}$ with the same simple largest eigenvalue $\lambda$ and with corresponding eigenvector $v=(N+1,-1,\dots,-1)$, by allowing the largest component of $G$ to be merely connected, instead of complete.
	\end{ex}
	
	\begin{figure}[H]
	
		\begin{minipage}[t]{0.45\textwidth}
		\centering
		
		
		\tikzset {_gzg2y2pdm/.code = {\pgfsetadditionalshadetransform{ \pgftransformshift{\pgfpoint{89.1 bp } { -128.7 bp }  }  \pgftransformscale{1.32 }  }}}
		\pgfdeclareradialshading{_v671db3dj}{\pgfpoint{-72bp}{104bp}}{rgb(0bp)=(1,1,1);
		rgb(0bp)=(1,1,1);
		rgb(25bp)=(0.29,0.56,0.89);
		rgb(400bp)=(0.29,0.56,0.89)}
		
		
		\tikzset {_rkqmm5jab/.code = {\pgfsetadditionalshadetransform{ \pgftransformshift{\pgfpoint{89.1 bp } { -128.7 bp }  }  \pgftransformscale{1.32 }  }}}
		\pgfdeclareradialshading{_rgic83cud}{\pgfpoint{-72bp}{104bp}}{rgb(0bp)=(1,1,1);
		rgb(0bp)=(1,1,1);
		rgb(25bp)=(0.29,0.56,0.89);
		rgb(400bp)=(0.29,0.56,0.89)}
		
		
		\tikzset {_hsr5qgbx8/.code = {\pgfsetadditionalshadetransform{ \pgftransformshift{\pgfpoint{89.1 bp } { -128.7 bp }  }  \pgftransformscale{1.32 }  }}}
		\pgfdeclareradialshading{_adpvps5zo}{\pgfpoint{-72bp}{104bp}}{rgb(0bp)=(1,1,1);
		rgb(0bp)=(1,1,1);
		rgb(25bp)=(0.29,0.56,0.89);
		rgb(400bp)=(0.29,0.56,0.89)}
		
		
		\tikzset {_5u3rn63fj/.code = {\pgfsetadditionalshadetransform{ \pgftransformshift{\pgfpoint{89.1 bp } { -128.7 bp }  }  \pgftransformscale{1.32 }  }}}
		\pgfdeclareradialshading{_1hwh6sttn}{\pgfpoint{-72bp}{104bp}}{rgb(0bp)=(1,1,1);
		rgb(0bp)=(1,1,1);
		rgb(25bp)=(0.29,0.56,0.89);
		rgb(400bp)=(0.29,0.56,0.89)}
		
		
		\tikzset {_9p2k5vr2n/.code = {\pgfsetadditionalshadetransform{ \pgftransformshift{\pgfpoint{89.1 bp } { -128.7 bp }  }  \pgftransformscale{1.32 }  }}}
		\pgfdeclareradialshading{_7zlcojqh9}{\pgfpoint{-72bp}{104bp}}{rgb(0bp)=(1,1,1);
		rgb(0bp)=(1,1,1);
		rgb(25bp)=(0.29,0.56,0.89);
		rgb(400bp)=(0.29,0.56,0.89)}
		
		
		\tikzset {_ipc1gofdg/.code = {\pgfsetadditionalshadetransform{ \pgftransformshift{\pgfpoint{89.1 bp } { -128.7 bp }  }  \pgftransformscale{1.32 }  }}}
		\pgfdeclareradialshading{_t69gd3kot}{\pgfpoint{-72bp}{104bp}}{rgb(0bp)=(1,1,1);
		rgb(0bp)=(1,1,1);
		rgb(25bp)=(0.29,0.56,0.89);
		rgb(400bp)=(0.29,0.56,0.89)}
		
		
		\tikzset {_ic3yhb760/.code = {\pgfsetadditionalshadetransform{ \pgftransformshift{\pgfpoint{89.1 bp } { -128.7 bp }  }  \pgftransformscale{1.32 }  }}}
		\pgfdeclareradialshading{_vfzfj1env}{\pgfpoint{-72bp}{104bp}}{rgb(0bp)=(1,1,1);
		rgb(0bp)=(1,1,1);
		rgb(25bp)=(0.29,0.56,0.89);
		rgb(400bp)=(0.29,0.56,0.89)}
		
		
		\tikzset {_re8uvmsfa/.code = {\pgfsetadditionalshadetransform{ \pgftransformshift{\pgfpoint{89.1 bp } { -128.7 bp }  }  \pgftransformscale{1.32 }  }}}
		\pgfdeclareradialshading{_0y43pqj98}{\pgfpoint{-72bp}{104bp}}{rgb(0bp)=(1,1,1);
		rgb(0bp)=(1,1,1);
		rgb(25bp)=(0.29,0.56,0.89);
		rgb(400bp)=(0.29,0.56,0.89)}
		
		
		\tikzset {_zgu2pk5b7/.code = {\pgfsetadditionalshadetransform{ \pgftransformshift{\pgfpoint{89.1 bp } { -128.7 bp }  }  \pgftransformscale{1.32 }  }}}
		\pgfdeclareradialshading{_pwxlj618f}{\pgfpoint{-72bp}{104bp}}{rgb(0bp)=(1,1,1);
		rgb(0bp)=(1,1,1);
		rgb(25bp)=(0.29,0.56,0.89);
		rgb(400bp)=(0.29,0.56,0.89)}
		\tikzset{every picture/.style={line width=0.75pt}} 
		
		\begin{tikzpicture}[x=0.75pt,y=0.75pt,yscale=-1,xscale=1]
		
		\draw [color={rgb, 255:red, 0; green, 0; blue, 0 }  ,draw opacity=1 ][shading=_v671db3dj,_gzg2y2pdm]   (179.66,87.28) -- (208.51,30) ;
		\draw [color={rgb, 255:red, 0; green, 0; blue, 0 }  ,draw opacity=1 ][shading=_rgic83cud,_rkqmm5jab][line width=0.75]    (88.01,30) -- (30.3,30) ;
		\draw [color={rgb, 255:red, 0; green, 0; blue, 0 }  ,draw opacity=1 ][shading=_adpvps5zo,_hsr5qgbx8]   (150.8,30) -- (179.66,87.28) ;
		\draw  [draw opacity=0][shading=_1hwh6sttn,_5u3rn63fj] (19.6,30) .. controls (19.6,24.09) and (24.39,19.3) .. (30.3,19.3) .. controls (36.21,19.3) and (41,24.09) .. (41,30) .. controls (41,35.91) and (36.21,40.7) .. (30.3,40.7) .. controls (24.39,40.7) and (19.6,35.91) .. (19.6,30) -- cycle ;
		\draw  [draw opacity=0][shading=_7zlcojqh9,_9p2k5vr2n] (77.31,30) .. controls (77.31,24.09) and (82.1,19.3) .. (88.01,19.3) .. controls (93.92,19.3) and (98.71,24.09) .. (98.71,30) .. controls (98.71,35.91) and (93.92,40.7) .. (88.01,40.7) .. controls (82.1,40.7) and (77.31,35.91) .. (77.31,30) -- cycle ;
		\draw  [draw opacity=0][shading=_t69gd3kot,_ipc1gofdg] (46.8,87.6) .. controls (46.8,81.69) and (51.59,76.9) .. (57.5,76.9) .. controls (63.41,76.9) and (68.2,81.69) .. (68.2,87.6) .. controls (68.2,93.51) and (63.41,98.3) .. (57.5,98.3) .. controls (51.59,98.3) and (46.8,93.51) .. (46.8,87.6) -- cycle ;
		\draw  [draw opacity=0][shading=_vfzfj1env,_ic3yhb760] (140.1,30) .. controls (140.1,24.09) and (144.89,19.3) .. (150.8,19.3) .. controls (156.71,19.3) and (161.5,24.09) .. (161.5,30) .. controls (161.5,35.91) and (156.71,40.7) .. (150.8,40.7) .. controls (144.89,40.7) and (140.1,35.91) .. (140.1,30) -- cycle ;
		\draw  [draw opacity=0][shading=_0y43pqj98,_re8uvmsfa] (197.81,30) .. controls (197.81,24.09) and (202.6,19.3) .. (208.51,19.3) .. controls (214.42,19.3) and (219.21,24.09) .. (219.21,30) .. controls (219.21,35.91) and (214.42,40.7) .. (208.51,40.7) .. controls (202.6,40.7) and (197.81,35.91) .. (197.81,30) -- cycle ;
		\draw  [draw opacity=0][shading=_pwxlj618f,_zgu2pk5b7] (168.96,87.28) .. controls (168.96,81.37) and (173.75,76.58) .. (179.66,76.58) .. controls (185.57,76.58) and (190.36,81.37) .. (190.36,87.28) .. controls (190.36,93.19) and (185.57,97.98) .. (179.66,97.98) .. controls (173.75,97.98) and (168.96,93.19) .. (168.96,87.28) -- cycle ;
		
		\draw (57.66,112.41) node  [font=\footnotesize]  {$G$};
		\draw (179.66,112.41) node  [font=\footnotesize]  {$G^{\circ }$};
		\draw (88.01,30) node  [font=\footnotesize] [align=left] {2};
		\draw (30.3,30) node  [font=\footnotesize] [align=left] {1};
		\draw (57.5,87.6) node  [font=\footnotesize] [align=left] {3};
		\draw (150.8,30) node  [font=\footnotesize] [align=left] {1};
		\draw (208.51,30) node  [font=\footnotesize] [align=left] {2};
		\draw (179.66,87.28) node  [font=\footnotesize] [align=left] {3};

		\end{tikzpicture}
		\caption{$G$ consists of precisely two disconnected components of order $1$ and $2$. The complement $G^{\circ}$ is a connected and non-regular graph.\label{f1}}
		\end{minipage}
		\hfill
		\begin{minipage}[t]{0.45\textwidth}
		\centering

		
		\tikzset {_ubji9p54x/.code = {\pgfsetadditionalshadetransform{ \pgftransformshift{\pgfpoint{89.1 bp } { -128.7 bp }  }  \pgftransformscale{1.32 }  }}}
		\pgfdeclareradialshading{_sj90jmkqz}{\pgfpoint{-72bp}{104bp}}{rgb(0bp)=(1,1,1);
		rgb(0bp)=(1,1,1);
		rgb(25bp)=(0.29,0.56,0.89);
		rgb(400bp)=(0.29,0.56,0.89)}
		
		
		\tikzset {_ll8386mw6/.code = {\pgfsetadditionalshadetransform{ \pgftransformshift{\pgfpoint{89.1 bp } { -128.7 bp }  }  \pgftransformscale{1.32 }  }}}
		\pgfdeclareradialshading{_ui1pbqvc3}{\pgfpoint{-72bp}{104bp}}{rgb(0bp)=(1,1,1);
		rgb(0bp)=(1,1,1);
		rgb(25bp)=(0.29,0.56,0.89);
		rgb(400bp)=(0.29,0.56,0.89)}
		
		
		\tikzset {_jug04fvhf/.code = {\pgfsetadditionalshadetransform{ \pgftransformshift{\pgfpoint{89.1 bp } { -128.7 bp }  }  \pgftransformscale{1.32 }  }}}
		\pgfdeclareradialshading{_agmsxonzg}{\pgfpoint{-72bp}{104bp}}{rgb(0bp)=(1,1,1);
		rgb(0bp)=(1,1,1);
		rgb(25bp)=(0.29,0.56,0.89);
		rgb(400bp)=(0.29,0.56,0.89)}
		
		
		\tikzset {_q5508zun0/.code = {\pgfsetadditionalshadetransform{ \pgftransformshift{\pgfpoint{89.1 bp } { -128.7 bp }  }  \pgftransformscale{1.32 }  }}}
		\pgfdeclareradialshading{_fbv7rp9k3}{\pgfpoint{-72bp}{104bp}}{rgb(0bp)=(1,1,1);
		rgb(0bp)=(1,1,1);
		rgb(25bp)=(0.29,0.56,0.89);
		rgb(400bp)=(0.29,0.56,0.89)}
		
		
		\tikzset {_w4z5ihz0n/.code = {\pgfsetadditionalshadetransform{ \pgftransformshift{\pgfpoint{89.1 bp } { -128.7 bp }  }  \pgftransformscale{1.32 }  }}}
		\pgfdeclareradialshading{_wabzuu8mz}{\pgfpoint{-72bp}{104bp}}{rgb(0bp)=(1,1,1);
		rgb(0bp)=(1,1,1);
		rgb(25bp)=(0.29,0.56,0.89);
		rgb(400bp)=(0.29,0.56,0.89)}
		
		
		\tikzset {_65wt3c76k/.code = {\pgfsetadditionalshadetransform{ \pgftransformshift{\pgfpoint{89.1 bp } { -128.7 bp }  }  \pgftransformscale{1.32 }  }}}
		\pgfdeclareradialshading{_h46mmbqgu}{\pgfpoint{-72bp}{104bp}}{rgb(0bp)=(1,1,1);
		rgb(0bp)=(1,1,1);
		rgb(25bp)=(0.29,0.56,0.89);
		rgb(400bp)=(0.29,0.56,0.89)}
		
		
		\tikzset {_n5wybcimk/.code = {\pgfsetadditionalshadetransform{ \pgftransformshift{\pgfpoint{89.1 bp } { -128.7 bp }  }  \pgftransformscale{1.32 }  }}}
		\pgfdeclareradialshading{_t01mkk2u0}{\pgfpoint{-72bp}{104bp}}{rgb(0bp)=(1,1,1);
		rgb(0bp)=(1,1,1);
		rgb(25bp)=(0.29,0.56,0.89);
		rgb(400bp)=(0.29,0.56,0.89)}
		
		
		\tikzset {_vfdtnkhow/.code = {\pgfsetadditionalshadetransform{ \pgftransformshift{\pgfpoint{89.1 bp } { -128.7 bp }  }  \pgftransformscale{1.32 }  }}}
		\pgfdeclareradialshading{_ic2xpbu5h}{\pgfpoint{-72bp}{104bp}}{rgb(0bp)=(1,1,1);
		rgb(0bp)=(1,1,1);
		rgb(25bp)=(0.29,0.56,0.89);
		rgb(400bp)=(0.29,0.56,0.89)}
		
		
		\tikzset {_8jgcxhl75/.code = {\pgfsetadditionalshadetransform{ \pgftransformshift{\pgfpoint{89.1 bp } { -128.7 bp }  }  \pgftransformscale{1.32 }  }}}
		\pgfdeclareradialshading{_ycyjm4405}{\pgfpoint{-72bp}{104bp}}{rgb(0bp)=(1,1,1);
		rgb(0bp)=(1,1,1);
		rgb(25bp)=(0.29,0.56,0.89);
		rgb(400bp)=(0.29,0.56,0.89)}
		
		
		\tikzset {_b22dwhzrj/.code = {\pgfsetadditionalshadetransform{ \pgftransformshift{\pgfpoint{89.1 bp } { -128.7 bp }  }  \pgftransformscale{1.32 }  }}}
		\pgfdeclareradialshading{_x9ibavta4}{\pgfpoint{-72bp}{104bp}}{rgb(0bp)=(1,1,1);
		rgb(0bp)=(1,1,1);
		rgb(25bp)=(0.29,0.56,0.89);
		rgb(400bp)=(0.29,0.56,0.89)}
		
		
		\tikzset {_twxnfn77m/.code = {\pgfsetadditionalshadetransform{ \pgftransformshift{\pgfpoint{89.1 bp } { -128.7 bp }  }  \pgftransformscale{1.32 }  }}}
		\pgfdeclareradialshading{_78mgvslo4}{\pgfpoint{-72bp}{104bp}}{rgb(0bp)=(1,1,1);
		rgb(0bp)=(1,1,1);
		rgb(25bp)=(0.29,0.56,0.89);
		rgb(400bp)=(0.29,0.56,0.89)}
		
		
		\tikzset {_mgxnnidv7/.code = {\pgfsetadditionalshadetransform{ \pgftransformshift{\pgfpoint{89.1 bp } { -128.7 bp }  }  \pgftransformscale{1.32 }  }}}
		\pgfdeclareradialshading{_2c71ei65g}{\pgfpoint{-72bp}{104bp}}{rgb(0bp)=(1,1,1);
		rgb(0bp)=(1,1,1);
		rgb(25bp)=(0.29,0.56,0.89);
		rgb(400bp)=(0.29,0.56,0.89)}
		
		
		\tikzset {_invu2fmpz/.code = {\pgfsetadditionalshadetransform{ \pgftransformshift{\pgfpoint{89.1 bp } { -128.7 bp }  }  \pgftransformscale{1.32 }  }}}
		\pgfdeclareradialshading{_t1owgf7bd}{\pgfpoint{-72bp}{104bp}}{rgb(0bp)=(1,1,1);
		rgb(0bp)=(1,1,1);
		rgb(25bp)=(0.29,0.56,0.89);
		rgb(400bp)=(0.29,0.56,0.89)}
		
		
		\tikzset {_xtb5ewpht/.code = {\pgfsetadditionalshadetransform{ \pgftransformshift{\pgfpoint{89.1 bp } { -128.7 bp }  }  \pgftransformscale{1.32 }  }}}
		\pgfdeclareradialshading{_foro80n4s}{\pgfpoint{-72bp}{104bp}}{rgb(0bp)=(1,1,1);
		rgb(0bp)=(1,1,1);
		rgb(25bp)=(0.29,0.56,0.89);
		rgb(400bp)=(0.29,0.56,0.89)}
		
		
		\tikzset {_6y46m3st1/.code = {\pgfsetadditionalshadetransform{ \pgftransformshift{\pgfpoint{89.1 bp } { -128.7 bp }  }  \pgftransformscale{1.32 }  }}}
		\pgfdeclareradialshading{_cajcmbu71}{\pgfpoint{-72bp}{104bp}}{rgb(0bp)=(1,1,1);
		rgb(0bp)=(1,1,1);
		rgb(25bp)=(0.29,0.56,0.89);
		rgb(400bp)=(0.29,0.56,0.89)}
		
		
		\tikzset {_3b7b9wyqi/.code = {\pgfsetadditionalshadetransform{ \pgftransformshift{\pgfpoint{89.1 bp } { -128.7 bp }  }  \pgftransformscale{1.32 }  }}}
		\pgfdeclareradialshading{_09qt8iv7s}{\pgfpoint{-72bp}{104bp}}{rgb(0bp)=(1,1,1);
		rgb(0bp)=(1,1,1);
		rgb(25bp)=(0.29,0.56,0.89);
		rgb(400bp)=(0.29,0.56,0.89)}
		
		
		\tikzset {_rtsvbcwiw/.code = {\pgfsetadditionalshadetransform{ \pgftransformshift{\pgfpoint{89.1 bp } { -128.7 bp }  }  \pgftransformscale{1.32 }  }}}
		\pgfdeclareradialshading{_sd65j3o8h}{\pgfpoint{-72bp}{104bp}}{rgb(0bp)=(1,1,1);
		rgb(0bp)=(1,1,1);
		rgb(25bp)=(0.29,0.56,0.89);
		rgb(400bp)=(0.29,0.56,0.89)}
		
		
		\tikzset {_t9b2kqg99/.code = {\pgfsetadditionalshadetransform{ \pgftransformshift{\pgfpoint{89.1 bp } { -128.7 bp }  }  \pgftransformscale{1.32 }  }}}
		\pgfdeclareradialshading{_8ulr93yth}{\pgfpoint{-72bp}{104bp}}{rgb(0bp)=(1,1,1);
		rgb(0bp)=(1,1,1);
		rgb(25bp)=(0.29,0.56,0.89);
		rgb(400bp)=(0.29,0.56,0.89)}
		
		
		\tikzset {_b5obh7jbm/.code = {\pgfsetadditionalshadetransform{ \pgftransformshift{\pgfpoint{89.1 bp } { -128.7 bp }  }  \pgftransformscale{1.32 }  }}}
		\pgfdeclareradialshading{_ra2qydpsc}{\pgfpoint{-72bp}{104bp}}{rgb(0bp)=(1,1,1);
		rgb(0bp)=(1,1,1);
		rgb(25bp)=(0.29,0.56,0.89);
		rgb(400bp)=(0.29,0.56,0.89)}
		
		
		\tikzset {_9307skan2/.code = {\pgfsetadditionalshadetransform{ \pgftransformshift{\pgfpoint{89.1 bp } { -128.7 bp }  }  \pgftransformscale{1.32 }  }}}
		\pgfdeclareradialshading{_uoyb8479p}{\pgfpoint{-72bp}{104bp}}{rgb(0bp)=(1,1,1);
		rgb(0bp)=(1,1,1);
		rgb(25bp)=(0.29,0.56,0.89);
		rgb(400bp)=(0.29,0.56,0.89)}
		\tikzset{every picture/.style={line width=0.75pt}} 
		
		\begin{tikzpicture}[x=0.75pt,y=0.75pt,yscale=-1,xscale=1]
		
		\draw [color={rgb, 255:red, 0; green, 0; blue, 0 }  ,draw opacity=1 ][shading=_sj90jmkqz,_ubji9p54x]   (280.71,98.4) -- (191.75,86.03) ;
		\draw [color={rgb, 255:red, 0; green, 0; blue, 0 }  ,draw opacity=1 ][shading=_ui1pbqvc3,_ll8386mw6]   (280.71,98.4) -- (191.75,34.37) ;
		\draw [color={rgb, 255:red, 0; green, 0; blue, 0 }  ,draw opacity=1 ][shading=_agmsxonzg,_jug04fvhf]   (237.68,114.3) -- (191.75,33.92) ;
		\draw [color={rgb, 255:red, 0; green, 0; blue, 0 }  ,draw opacity=1 ][shading=_fbv7rp9k3,_q5508zun0]   (80.22,113.42) -- (123.25,97.52) ;
		\draw [color={rgb, 255:red, 0; green, 0; blue, 0 }  ,draw opacity=1 ][shading=_wabzuu8mz,_w4z5ihz0n]   (34.29,85.15) -- (94.24,33.04) ;
		\draw [color={rgb, 255:red, 0; green, 0; blue, 0 }  ,draw opacity=1 ][shading=_h46mmbqgu,_65wt3c76k]   (34.29,33.04) -- (94.24,33.04) ;
		\draw [color={rgb, 255:red, 0; green, 0; blue, 0 }  ,draw opacity=1 ][shading=_t01mkk2u0,_n5wybcimk]   (191.75,86.03) -- (237.68,114.3) ;
		\draw [color={rgb, 255:red, 0; green, 0; blue, 0 }  ,draw opacity=1 ][shading=_ic2xpbu5h,_vfdtnkhow]   (237.68,114.3) -- (251.7,33.92) ;
		\draw [color={rgb, 255:red, 0; green, 0; blue, 0 }  ,draw opacity=1 ][shading=_ycyjm4405,_8jgcxhl75]   (280.71,98.4) -- (251.7,33.92) ;
		\draw [color={rgb, 255:red, 0; green, 0; blue, 0 }  ,draw opacity=1 ][shading=_x9ibavta4,_b22dwhzrj]   (191.75,86.03) -- (191.75,33.92) ;
		\draw  [draw opacity=0][shading=_78mgvslo4,_twxnfn77m] (23.59,33.04) .. controls (23.59,27.13) and (28.38,22.34) .. (34.29,22.34) .. controls (40.2,22.34) and (44.99,27.13) .. (44.99,33.04) .. controls (44.99,38.95) and (40.2,43.74) .. (34.29,43.74) .. controls (28.38,43.74) and (23.59,38.95) .. (23.59,33.04) -- cycle ;
		\draw  [draw opacity=0][shading=_2c71ei65g,_mgxnnidv7] (23.59,85.15) .. controls (23.59,79.24) and (28.38,74.45) .. (34.29,74.45) .. controls (40.2,74.45) and (44.99,79.24) .. (44.99,85.15) .. controls (44.99,91.06) and (40.2,95.85) .. (34.29,95.85) .. controls (28.38,95.85) and (23.59,91.06) .. (23.59,85.15) -- cycle ;
		\draw  [draw opacity=0][shading=_t1owgf7bd,_invu2fmpz] (83.54,33.04) .. controls (83.54,27.13) and (88.33,22.34) .. (94.24,22.34) .. controls (100.15,22.34) and (104.94,27.13) .. (104.94,33.04) .. controls (104.94,38.95) and (100.15,43.74) .. (94.24,43.74) .. controls (88.33,43.74) and (83.54,38.95) .. (83.54,33.04) -- cycle ;
		\draw  [draw opacity=0][shading=_foro80n4s,_xtb5ewpht] (69.52,113.42) .. controls (69.52,107.51) and (74.31,102.72) .. (80.22,102.72) .. controls (86.13,102.72) and (90.92,107.51) .. (90.92,113.42) .. controls (90.92,119.32) and (86.13,124.11) .. (80.22,124.11) .. controls (74.31,124.11) and (69.52,119.32) .. (69.52,113.42) -- cycle ;
		\draw  [draw opacity=0][shading=_cajcmbu71,_6y46m3st1] (112.55,97.52) .. controls (112.55,91.61) and (117.34,86.82) .. (123.25,86.82) .. controls (129.15,86.82) and (133.95,91.61) .. (133.95,97.52) .. controls (133.95,103.43) and (129.15,108.22) .. (123.25,108.22) .. controls (117.34,108.22) and (112.55,103.43) .. (112.55,97.52) -- cycle ;
		\draw  [draw opacity=0][shading=_09qt8iv7s,_3b7b9wyqi] (181.05,33.92) .. controls (181.05,28.01) and (185.85,23.22) .. (191.75,23.22) .. controls (197.66,23.22) and (202.45,28.01) .. (202.45,33.92) .. controls (202.45,39.83) and (197.66,44.62) .. (191.75,44.62) .. controls (185.85,44.62) and (181.05,39.83) .. (181.05,33.92) -- cycle ;
		\draw  [draw opacity=0][shading=_sd65j3o8h,_rtsvbcwiw] (241.01,33.92) .. controls (241.01,28.01) and (245.8,23.22) .. (251.7,23.22) .. controls (257.61,23.22) and (262.4,28.01) .. (262.4,33.92) .. controls (262.4,39.83) and (257.61,44.62) .. (251.7,44.62) .. controls (245.8,44.62) and (241.01,39.83) .. (241.01,33.92) -- cycle ;
		\draw  [draw opacity=0][shading=_8ulr93yth,_t9b2kqg99] (181.05,86.03) .. controls (181.05,80.13) and (185.85,75.34) .. (191.75,75.34) .. controls (197.66,75.34) and (202.45,80.13) .. (202.45,86.03) .. controls (202.45,91.94) and (197.66,96.73) .. (191.75,96.73) .. controls (185.85,96.73) and (181.05,91.94) .. (181.05,86.03) -- cycle ;
		\draw  [draw opacity=0][shading=_ra2qydpsc,_b5obh7jbm] (226.98,114.3) .. controls (226.98,108.39) and (231.78,103.6) .. (237.68,103.6) .. controls (243.59,103.6) and (248.38,108.39) .. (248.38,114.3) .. controls (248.38,120.21) and (243.59,125) .. (237.68,125) .. controls (231.78,125) and (226.98,120.21) .. (226.98,114.3) -- cycle ;
		\draw  [draw opacity=0][shading=_uoyb8479p,_9307skan2] (270.01,98.4) .. controls (270.01,92.49) and (274.8,87.7) .. (280.71,87.7) .. controls (286.62,87.7) and (291.41,92.49) .. (291.41,98.4) .. controls (291.41,104.31) and (286.62,109.1) .. (280.71,109.1) .. controls (274.8,109.1) and (270.01,104.31) .. (270.01,98.4) -- cycle ;

		\draw (225.6,144.21) node  [font=\footnotesize]  {$G^{\circ }$};
		\draw (64.26,144.21) node  [font=\footnotesize]  {$G$};
		\draw (123.25,97.52) node   [align=left] {5};
		\draw (80.22,113.42) node   [align=left] {4};
		\draw (34.29,85.15) node   [align=left] {3};
		\draw (94.24,33.04) node   [align=left] {2};
		\draw (34.29,33.04) node   [align=left] {1};
		\draw (191.75,86.03) node   [align=left] {3};
		\draw (251.7,33.92) node   [align=left] {2};
		\draw (237.68,114.3) node   [align=left] {4};
		\draw (280.71,98.4) node   [align=left] {5};
		\draw (191.75,33.92) node   [align=left] {1};

		\end{tikzpicture}

		\caption{$G$ consists of precisely two disconnected components of order $2$ and $3$ and its complement $G^{\circ}$ is connected.\label{f2}}
		\end{minipage}
	
    \centering		
  
    \tikzset {_7ywv8fedf/.code = {\pgfsetadditionalshadetransform{ \pgftransformshift{\pgfpoint{89.1 bp } { -128.7 bp }  }  \pgftransformscale{1.32 }  }}}
    \pgfdeclareradialshading{_ug958oisi}{\pgfpoint{-72bp}{104bp}}{rgb(0bp)=(1,1,1);
    rgb(0bp)=(1,1,1);
    rgb(25bp)=(0.29,0.56,0.89);
    rgb(400bp)=(0.29,0.56,0.89)}

  
    \tikzset {_ldokxi9rq/.code = {\pgfsetadditionalshadetransform{ \pgftransformshift{\pgfpoint{89.1 bp } { -128.7 bp }  }  \pgftransformscale{1.32 }  }}}
    \pgfdeclareradialshading{_stir83qji}{\pgfpoint{-72bp}{104bp}}{rgb(0bp)=(1,1,1);
    rgb(0bp)=(1,1,1);
    rgb(25bp)=(0.29,0.56,0.89);
    rgb(400bp)=(0.29,0.56,0.89)}

  
    \tikzset {_jayfo11z4/.code = {\pgfsetadditionalshadetransform{ \pgftransformshift{\pgfpoint{89.1 bp } { -128.7 bp }  }  \pgftransformscale{1.32 }  }}}
    \pgfdeclareradialshading{_hdvbwfrgt}{\pgfpoint{-72bp}{104bp}}{rgb(0bp)=(1,1,1);
    rgb(0bp)=(1,1,1);
    rgb(25bp)=(0.29,0.56,0.89);
    rgb(400bp)=(0.29,0.56,0.89)}

  
    \tikzset {_pz14a4jga/.code = {\pgfsetadditionalshadetransform{ \pgftransformshift{\pgfpoint{89.1 bp } { -128.7 bp }  }  \pgftransformscale{1.32 }  }}}
    \pgfdeclareradialshading{_hy0ac92es}{\pgfpoint{-72bp}{104bp}}{rgb(0bp)=(1,1,1);
    rgb(0bp)=(1,1,1);
    rgb(25bp)=(0.29,0.56,0.89);
    rgb(400bp)=(0.29,0.56,0.89)}

  
    \tikzset {_6dn4x0pvw/.code = {\pgfsetadditionalshadetransform{ \pgftransformshift{\pgfpoint{89.1 bp } { -128.7 bp }  }  \pgftransformscale{1.32 }  }}}
    \pgfdeclareradialshading{_xofz3wx6h}{\pgfpoint{-72bp}{104bp}}{rgb(0bp)=(1,1,1);
    rgb(0bp)=(1,1,1);
    rgb(25bp)=(0.29,0.56,0.89);
    rgb(400bp)=(0.29,0.56,0.89)}

  
    \tikzset {_fimdqawn2/.code = {\pgfsetadditionalshadetransform{ \pgftransformshift{\pgfpoint{89.1 bp } { -128.7 bp }  }  \pgftransformscale{1.32 }  }}}
    \pgfdeclareradialshading{_5q5tht23e}{\pgfpoint{-72bp}{104bp}}{rgb(0bp)=(1,1,1);
    rgb(0bp)=(1,1,1);
    rgb(25bp)=(0.29,0.56,0.89);
    rgb(400bp)=(0.29,0.56,0.89)}

  
    \tikzset {_wikkj9663/.code = {\pgfsetadditionalshadetransform{ \pgftransformshift{\pgfpoint{89.1 bp } { -128.7 bp }  }  \pgftransformscale{1.32 }  }}}
\pgfdeclareradialshading{_9mzif6sm7}{\pgfpoint{-72bp}{104bp}}{rgb(0bp)=(1,1,1);
rgb(0bp)=(1,1,1);
rgb(25bp)=(0.29,0.56,0.89);
rgb(400bp)=(0.29,0.56,0.89)}

  
\tikzset {_0frnzby3m/.code = {\pgfsetadditionalshadetransform{ \pgftransformshift{\pgfpoint{89.1 bp } { -128.7 bp }  }  \pgftransformscale{1.32 }  }}}
\pgfdeclareradialshading{_o1dqiupkg}{\pgfpoint{-72bp}{104bp}}{rgb(0bp)=(1,1,1);
rgb(0bp)=(1,1,1);
rgb(25bp)=(0.48,0.15,0.15);
rgb(400bp)=(0.48,0.15,0.15)}

  
\tikzset {_oebqfz9lv/.code = {\pgfsetadditionalshadetransform{ \pgftransformshift{\pgfpoint{89.1 bp } { -128.7 bp }  }  \pgftransformscale{1.32 }  }}}
\pgfdeclareradialshading{_5ks0gwk9h}{\pgfpoint{-72bp}{104bp}}{rgb(0bp)=(1,1,1);
rgb(0bp)=(1,1,1);
rgb(25bp)=(0.29,0.56,0.89);
rgb(400bp)=(0.29,0.56,0.89)}

  
\tikzset {_m4jev9wte/.code = {\pgfsetadditionalshadetransform{ \pgftransformshift{\pgfpoint{89.1 bp } { -128.7 bp }  }  \pgftransformscale{1.32 }  }}}
\pgfdeclareradialshading{_9p6jii84o}{\pgfpoint{-72bp}{104bp}}{rgb(0bp)=(1,1,1);
rgb(0bp)=(1,1,1);
rgb(25bp)=(0.29,0.56,0.89);
rgb(400bp)=(0.29,0.56,0.89)}

  
\tikzset {_fjqf1fz4l/.code = {\pgfsetadditionalshadetransform{ \pgftransformshift{\pgfpoint{89.1 bp } { -128.7 bp }  }  \pgftransformscale{1.32 }  }}}
\pgfdeclareradialshading{_na9scb270}{\pgfpoint{-72bp}{104bp}}{rgb(0bp)=(1,1,1);
rgb(0bp)=(1,1,1);
rgb(25bp)=(0.29,0.56,0.89);
rgb(400bp)=(0.29,0.56,0.89)}

  
\tikzset {_at2nz4284/.code = {\pgfsetadditionalshadetransform{ \pgftransformshift{\pgfpoint{89.1 bp } { -128.7 bp }  }  \pgftransformscale{1.32 }  }}}
\pgfdeclareradialshading{_q9j1vvuzs}{\pgfpoint{-72bp}{104bp}}{rgb(0bp)=(1,1,1);
rgb(0bp)=(1,1,1);
rgb(25bp)=(0.29,0.56,0.89);
rgb(400bp)=(0.29,0.56,0.89)}

  
\tikzset {_l56pski6r/.code = {\pgfsetadditionalshadetransform{ \pgftransformshift{\pgfpoint{89.1 bp } { -128.7 bp }  }  \pgftransformscale{1.32 }  }}}
\pgfdeclareradialshading{_nccmja4vd}{\pgfpoint{-72bp}{104bp}}{rgb(0bp)=(1,1,1);
rgb(0bp)=(1,1,1);
rgb(25bp)=(0.29,0.56,0.89);
rgb(400bp)=(0.29,0.56,0.89)}

  
\tikzset {_sxhfbk87s/.code = {\pgfsetadditionalshadetransform{ \pgftransformshift{\pgfpoint{89.1 bp } { -128.7 bp }  }  \pgftransformscale{1.32 }  }}}
\pgfdeclareradialshading{_iw8tse0rr}{\pgfpoint{-72bp}{104bp}}{rgb(0bp)=(1,1,1);
rgb(0bp)=(1,1,1);
rgb(25bp)=(0.29,0.56,0.89);
rgb(400bp)=(0.29,0.56,0.89)}

  
\tikzset {_vtdr5hjle/.code = {\pgfsetadditionalshadetransform{ \pgftransformshift{\pgfpoint{89.1 bp } { -128.7 bp }  }  \pgftransformscale{1.32 }  }}}
\pgfdeclareradialshading{_l3b8fm862}{\pgfpoint{-72bp}{104bp}}{rgb(0bp)=(1,1,1);
rgb(0bp)=(1,1,1);
rgb(25bp)=(0.29,0.56,0.89);
rgb(400bp)=(0.29,0.56,0.89)}

  
\tikzset {_ctkgni188/.code = {\pgfsetadditionalshadetransform{ \pgftransformshift{\pgfpoint{89.1 bp } { -128.7 bp }  }  \pgftransformscale{1.32 }  }}}
\pgfdeclareradialshading{_e9i7mwtlb}{\pgfpoint{-72bp}{104bp}}{rgb(0bp)=(1,1,1);
rgb(0bp)=(1,1,1);
rgb(25bp)=(0.48,0.15,0.15);
rgb(400bp)=(0.48,0.15,0.15)}
\tikzset{every picture/.style={line width=0.75pt}} 

\begin{tikzpicture}[x=0.75pt,y=0.75pt,yscale=-1,xscale=1]

\draw    (189.25,189.75) -- (177,196.73) ;
\draw    (145.67,214.47) -- (129.25,223.75) ;
\draw    (140.75,87) -- (176.25,118.75) ;
\draw    (85.25,71.75) -- (140.75,87) ;
\draw    (85.25,71.75) -- (176.25,118.75) ;
\draw    (129.25,223.75) -- (85.25,71.75) ;
\draw    (54.25,117.75) -- (140.75,87) ;
\draw    (85.25,71.75) -- (189.25,189.75) ;
\draw    (85.25,71.75) -- (65.25,186.75) ;
\draw    (54.25,117.75) -- (85.25,71.75) ;
\draw    (189.25,189.75) -- (54.25,117.75) ;
\draw    (54.25,117.75) -- (129.25,223.75) ;
\draw    (65.25,186.75) -- (176.25,118.75) ;
\draw    (65.25,186.75) -- (140.75,87) ;
\draw    (189.25,189.75) -- (65.25,186.75) ;
\draw    (176.25,118.75) -- (54.25,117.75) ;
\draw    (65.25,186.75) -- (129.25,223.75) ;
\draw    (129.25,223.75) -- (176.25,118.75) ;
\draw    (189.25,189.75) -- (140.75,87) ;
\draw    (129.25,223.75) -- (140.75,87) ;
\draw    (65.25,186.75) -- (54.25,117.75) ;
\draw    (176.25,118.75) -- (189.25,189.75) ;
\draw  [draw opacity=0][shading=_ug958oisi,_7ywv8fedf] (76,71.75) .. controls (76,66.64) and (80.14,62.5) .. (85.25,62.5) .. controls (90.36,62.5) and (94.5,66.64) .. (94.5,71.75) .. controls (94.5,76.86) and (90.36,81) .. (85.25,81) .. controls (80.14,81) and (76,76.86) .. (76,71.75) -- cycle ;
\draw  [draw opacity=0][shading=_stir83qji,_ldokxi9rq] (45,117.75) .. controls (45,112.64) and (49.14,108.5) .. (54.25,108.5) .. controls (59.36,108.5) and (63.5,112.64) .. (63.5,117.75) .. controls (63.5,122.86) and (59.36,127) .. (54.25,127) .. controls (49.14,127) and (45,122.86) .. (45,117.75) -- cycle ;
\draw  [draw opacity=0][shading=_hdvbwfrgt,_jayfo11z4] (56,186.75) .. controls (56,181.64) and (60.14,177.5) .. (65.25,177.5) .. controls (70.36,177.5) and (74.5,181.64) .. (74.5,186.75) .. controls (74.5,191.86) and (70.36,196) .. (65.25,196) .. controls (60.14,196) and (56,191.86) .. (56,186.75) -- cycle ;
\draw  [draw opacity=0][shading=_hy0ac92es,_pz14a4jga] (120,223.75) .. controls (120,218.64) and (124.14,214.5) .. (129.25,214.5) .. controls (134.36,214.5) and (138.5,218.64) .. (138.5,223.75) .. controls (138.5,228.86) and (134.36,233) .. (129.25,233) .. controls (124.14,233) and (120,228.86) .. (120,223.75) -- cycle ;
\draw  [draw opacity=0][shading=_xofz3wx6h,_6dn4x0pvw] (180,189.75) .. controls (180,184.64) and (184.14,180.5) .. (189.25,180.5) .. controls (194.36,180.5) and (198.5,184.64) .. (198.5,189.75) .. controls (198.5,194.86) and (194.36,199) .. (189.25,199) .. controls (184.14,199) and (180,194.86) .. (180,189.75) -- cycle ;
\draw  [draw opacity=0][shading=_5q5tht23e,_fimdqawn2] (167,118.75) .. controls (167,113.64) and (171.14,109.5) .. (176.25,109.5) .. controls (181.36,109.5) and (185.5,113.64) .. (185.5,118.75) .. controls (185.5,123.86) and (181.36,128) .. (176.25,128) .. controls (171.14,128) and (167,123.86) .. (167,118.75) -- cycle ;
\draw  [draw opacity=0][shading=_9mzif6sm7,_wikkj9663] (131.5,87) .. controls (131.5,81.89) and (135.64,77.75) .. (140.75,77.75) .. controls (145.86,77.75) and (150,81.89) .. (150,87) .. controls (150,92.11) and (145.86,96.25) .. (140.75,96.25) .. controls (135.64,96.25) and (131.5,92.11) .. (131.5,87) -- cycle ;
\draw  [draw opacity=0][shading=_o1dqiupkg,_0frnzby3m] (110.6,135.2) .. controls (110.6,129.29) and (115.39,124.5) .. (121.3,124.5) .. controls (127.21,124.5) and (132,129.29) .. (132,135.2) .. controls (132,141.11) and (127.21,145.9) .. (121.3,145.9) .. controls (115.39,145.9) and (110.6,141.11) .. (110.6,135.2) -- cycle ;
\draw  [dash pattern={on 0.84pt off 2.51pt}]  (145.67,214.47) -- (177,196.73) ;

\draw    (357.8,137.95) -- (412.75,121.5) ;
\draw    (321.75,74.5) -- (357.8,137.95) ;
\draw    (365.75,226.5) -- (357.8,137.95) ;
\draw    (357.8,137.95) -- (377.25,89.75) ;
\draw    (357.8,137.95) -- (425.75,192.5) ;
\draw    (290.75,120.5) -- (357.8,137.95) ;
\draw    (301.75,189.5) -- (357.8,137.95) ;
\draw  [draw opacity=0][shading=_5ks0gwk9h,_oebqfz9lv] (312.5,74.5) .. controls (312.5,69.39) and (316.64,65.25) .. (321.75,65.25) .. controls (326.86,65.25) and (331,69.39) .. (331,74.5) .. controls (331,79.61) and (326.86,83.75) .. (321.75,83.75) .. controls (316.64,83.75) and (312.5,79.61) .. (312.5,74.5) -- cycle ;
\draw  [draw opacity=0][shading=_9p6jii84o,_m4jev9wte] (281.5,120.5) .. controls (281.5,115.39) and (285.64,111.25) .. (290.75,111.25) .. controls (295.86,111.25) and (300,115.39) .. (300,120.5) .. controls (300,125.61) and (295.86,129.75) .. (290.75,129.75) .. controls (285.64,129.75) and (281.5,125.61) .. (281.5,120.5) -- cycle ;
\draw  [draw opacity=0][shading=_na9scb270,_fjqf1fz4l] (292.5,189.5) .. controls (292.5,184.39) and (296.64,180.25) .. (301.75,180.25) .. controls (306.86,180.25) and (311,184.39) .. (311,189.5) .. controls (311,194.61) and (306.86,198.75) .. (301.75,198.75) .. controls (296.64,198.75) and (292.5,194.61) .. (292.5,189.5) -- cycle ;
\draw  [draw opacity=0][shading=_q9j1vvuzs,_at2nz4284] (356.5,226.5) .. controls (356.5,221.39) and (360.64,217.25) .. (365.75,217.25) .. controls (370.86,217.25) and (375,221.39) .. (375,226.5) .. controls (375,231.61) and (370.86,235.75) .. (365.75,235.75) .. controls (360.64,235.75) and (356.5,231.61) .. (356.5,226.5) -- cycle ;
\draw  [draw opacity=0][shading=_nccmja4vd,_l56pski6r] (416.5,192.5) .. controls (416.5,187.39) and (420.64,183.25) .. (425.75,183.25) .. controls (430.86,183.25) and (435,187.39) .. (435,192.5) .. controls (435,197.61) and (430.86,201.75) .. (425.75,201.75) .. controls (420.64,201.75) and (416.5,197.61) .. (416.5,192.5) -- cycle ;
\draw  [draw opacity=0][shading=_iw8tse0rr,_sxhfbk87s] (403.5,121.5) .. controls (403.5,116.39) and (407.64,112.25) .. (412.75,112.25) .. controls (417.86,112.25) and (422,116.39) .. (422,121.5) .. controls (422,126.61) and (417.86,130.75) .. (412.75,130.75) .. controls (407.64,130.75) and (403.5,126.61) .. (403.5,121.5) -- cycle ;
\draw  [draw opacity=0][shading=_l3b8fm862,_vtdr5hjle] (368,89.75) .. controls (368,84.64) and (372.14,80.5) .. (377.25,80.5) .. controls (382.36,80.5) and (386.5,84.64) .. (386.5,89.75) .. controls (386.5,94.86) and (382.36,99) .. (377.25,99) .. controls (372.14,99) and (368,94.86) .. (368,89.75) -- cycle ;
\draw  [draw opacity=0][shading=_e9i7mwtlb,_ctkgni188] (347.1,137.95) .. controls (347.1,132.04) and (351.89,127.25) .. (357.8,127.25) .. controls (363.71,127.25) and (368.5,132.04) .. (368.5,137.95) .. controls (368.5,143.86) and (363.71,148.65) .. (357.8,148.65) .. controls (351.89,148.65) and (347.1,143.86) .. (347.1,137.95) -- cycle ;
\draw  [dash pattern={on 0.84pt off 2.51pt}]  (382.17,217.22) -- (413.5,199.48) ;

\draw (120,255.75) node  [font=\footnotesize]  {$G$};
\draw (67.25,185.95) node  [font=\footnotesize]  {$2$};
\draw (56.25,116.95) node  [font=\footnotesize]  {$3$};
\draw (87.25,70.55) node  [font=\footnotesize]  {$4$};
\draw (143.15,86.2) node  [font=\footnotesize]  {$5$};
\draw (178.25,117.55) node  [font=\footnotesize]  {$6$};
\draw (191.25,189.35) node  [font=\footnotesize]  {$7$};
\draw (129.25,223.75) node  [font=\footnotesize]  {$N$};
\draw (122.1,133.2) node    {$1$};
\draw (356.5,258.5) node  [font=\footnotesize]  {$G^{\circ }$};
\draw (303.75,188.7) node  [font=\footnotesize]  {$2$};
\draw (292.75,119.7) node  [font=\footnotesize]  {$3$};
\draw (323.75,73.3) node  [font=\footnotesize]  {$4$};
\draw (379.65,88.95) node  [font=\footnotesize]  {$5$};
\draw (414.75,120.3) node  [font=\footnotesize]  {$6$};
\draw (427.75,192.1) node  [font=\footnotesize]  {$7$};
\draw (365.75,226.5) node  [font=\footnotesize]  {$N$};
\draw (358.6,135.95) node    {$1$};

        \end{tikzpicture}

        \caption{$G$ consists of precisely two disconnected components of different sizes $N$ and $1.$ The complement $G^{\circ}$ is a connected  graph.\label{f3}}
		\end{figure}
	
	In what follows we turn to negative examples. The first of them shows us the importance of starting with connected components of different order, whereas the second one shows us what goes wrong if we start with more than $2$ components.
	
	\begin{ex}
		Let $G=(V,E)$ be the undirected graph with $V=\{1,2,3,4\}$ and two disconnected components, this time of the same order $s=t=2$, shown in Figure \ref{f4}. Then $G^{\circ}$ is a connected graph with
		\begin{equation}
			\nonumber L_{G^{\circ}}=\left(\begin{array}{cccc}
				2&0&-1&-1\\
				0&2&-1&-1\\
				-1&-1&2&0\\
				-1&-1&0&2
			\end{array}\right).
		\end{equation}
		
		Here $\spec(L_{G^{\circ}})=\{4,2,2,0\}$ with \textbf{simple} and \textbf{largest eigenvalue} $\lambda=4.$ However, there are no eigenvectors satisfying $\sum_{i=1}^{4}\nu_{i}^{\ell}\neq0$ for all $\ell>1,$ except multiples of $\mathbf{1}=(1,1,1,1).$ The eigenvectors for the other eigenvalues satisfy  $\sum_{i=1}^{4}\nu_i^{\ell}= 0$ whenever $\ell$ is odd. This example indicates that symmetry can be an obstruction for $\rho$-versatility.
	\end{ex}


%
	
	\begin{ex}
		Let $G=(V,E)$ be the undirected graph with $V=\{1,2,3,4,5,6\}$ and three disconnected components of order $s=1,$ $t=2$ and $r=3,$ shown in Figure \ref{f6}. Then $G^{\circ}$ is a connected graph with
		\begin{equation}
			\nonumber L_{G^{\circ}}=\left(\begin{array}{cccccc}
				4&0&-1&-1&-1&-1\\
				0&3&0&-1&-1&-1\\
				-1&0&4&-1&-1&-1\\
				-1&-1&-1&4&0&-1\\
				-1&-1&-1&0&4&-1\\
				-1&-1&-1&-1&-1&5
			\end{array}\right)_{6\times 6}.
		\end{equation}
		
		We have $\spec(L_{G^{\circ}})=\{6,6,5,4,3,0\}$ with \textbf{non-simple} and \textbf{largest eigenvalue} $\lambda_{1,2}=6$. Nevertheless, two corresponding eigenvectors are given by $\{(-1,-1,-1,0,0,3)$ and $(-2,-2,-2,3,3,0)\}$, which both satisfy $\sum_{i=1}^{6}\nu_{i}^{\ell}\neq0$ for all powers $\ell>1.$
	\end{ex}
	
%
		
	\begin{figure}[H]
		\begin{minipage}[t]{0.45\textwidth}
			\centering
			
			\tikzset {_m5t92ypuz/.code = {\pgfsetadditionalshadetransform{ \pgftransformshift{\pgfpoint{89.1 bp } { -128.7 bp }  }  \pgftransformscale{1.32 }  }}}
			\pgfdeclareradialshading{_hr48qomrb}{\pgfpoint{-72bp}{104bp}}{rgb(0bp)=(1,1,1);
				rgb(0bp)=(1,1,1);
				rgb(25bp)=(0.29,0.56,0.89);
				rgb(400bp)=(0.29,0.56,0.89)}
			
			
			\tikzset {_qwkbg84js/.code = {\pgfsetadditionalshadetransform{ \pgftransformshift{\pgfpoint{89.1 bp } { -128.7 bp }  }  \pgftransformscale{1.32 }  }}}
			\pgfdeclareradialshading{_84c4b7as9}{\pgfpoint{-72bp}{104bp}}{rgb(0bp)=(1,1,1);
				rgb(0bp)=(1,1,1);
				rgb(25bp)=(0.29,0.56,0.89);
				rgb(400bp)=(0.29,0.56,0.89)}
			
			
			\tikzset {_5d8z04i89/.code = {\pgfsetadditionalshadetransform{ \pgftransformshift{\pgfpoint{89.1 bp } { -128.7 bp }  }  \pgftransformscale{1.32 }  }}}
			\pgfdeclareradialshading{_5iwkzmjsk}{\pgfpoint{-72bp}{104bp}}{rgb(0bp)=(1,1,1);
				rgb(0bp)=(1,1,1);
				rgb(25bp)=(0.29,0.56,0.89);
				rgb(400bp)=(0.29,0.56,0.89)}
			
			
			\tikzset {_s0h72ppzf/.code = {\pgfsetadditionalshadetransform{ \pgftransformshift{\pgfpoint{89.1 bp } { -128.7 bp }  }  \pgftransformscale{1.32 }  }}}
			\pgfdeclareradialshading{_gca0ma2y6}{\pgfpoint{-72bp}{104bp}}{rgb(0bp)=(1,1,1);
				rgb(0bp)=(1,1,1);
				rgb(25bp)=(0.29,0.56,0.89);
				rgb(400bp)=(0.29,0.56,0.89)}
			
			
			\tikzset {_wm67zkktr/.code = {\pgfsetadditionalshadetransform{ \pgftransformshift{\pgfpoint{89.1 bp } { -128.7 bp }  }  \pgftransformscale{1.32 }  }}}
			\pgfdeclareradialshading{_jouysak7d}{\pgfpoint{-72bp}{104bp}}{rgb(0bp)=(1,1,1);
				rgb(0bp)=(1,1,1);
				rgb(25bp)=(0.29,0.56,0.89);
				rgb(400bp)=(0.29,0.56,0.89)}
			
			
			\tikzset {_2umrx5dr6/.code = {\pgfsetadditionalshadetransform{ \pgftransformshift{\pgfpoint{89.1 bp } { -128.7 bp }  }  \pgftransformscale{1.32 }  }}}
			\pgfdeclareradialshading{_2mexbeuiv}{\pgfpoint{-72bp}{104bp}}{rgb(0bp)=(1,1,1);
				rgb(0bp)=(1,1,1);
				rgb(25bp)=(0.29,0.56,0.89);
				rgb(400bp)=(0.29,0.56,0.89)}
			
			
			\tikzset {_0ibysugxg/.code = {\pgfsetadditionalshadetransform{ \pgftransformshift{\pgfpoint{89.1 bp } { -128.7 bp }  }  \pgftransformscale{1.32 }  }}}
			\pgfdeclareradialshading{_jvsh7ppcb}{\pgfpoint{-72bp}{104bp}}{rgb(0bp)=(1,1,1);
				rgb(0bp)=(1,1,1);
				rgb(25bp)=(0.29,0.56,0.89);
				rgb(400bp)=(0.29,0.56,0.89)}
			
			
			\tikzset {_x2y3pyfpa/.code = {\pgfsetadditionalshadetransform{ \pgftransformshift{\pgfpoint{89.1 bp } { -128.7 bp }  }  \pgftransformscale{1.32 }  }}}
			\pgfdeclareradialshading{_thuovt89s}{\pgfpoint{-72bp}{104bp}}{rgb(0bp)=(1,1,1);
				rgb(0bp)=(1,1,1);
				rgb(25bp)=(0.29,0.56,0.89);
				rgb(400bp)=(0.29,0.56,0.89)}
			
			
			\tikzset {_prd9o6hbr/.code = {\pgfsetadditionalshadetransform{ \pgftransformshift{\pgfpoint{89.1 bp } { -128.7 bp }  }  \pgftransformscale{1.32 }  }}}
			\pgfdeclareradialshading{_qh9idu2gg}{\pgfpoint{-72bp}{104bp}}{rgb(0bp)=(1,1,1);
				rgb(0bp)=(1,1,1);
				rgb(25bp)=(0.29,0.56,0.89);
				rgb(400bp)=(0.29,0.56,0.89)}
			
			
			\tikzset {_832d3j4ow/.code = {\pgfsetadditionalshadetransform{ \pgftransformshift{\pgfpoint{89.1 bp } { -128.7 bp }  }  \pgftransformscale{1.32 }  }}}
			\pgfdeclareradialshading{_a2he04xrr}{\pgfpoint{-72bp}{104bp}}{rgb(0bp)=(1,1,1);
				rgb(0bp)=(1,1,1);
				rgb(25bp)=(0.29,0.56,0.89);
				rgb(400bp)=(0.29,0.56,0.89)}
			
			
			\tikzset {_pr3hpz9l5/.code = {\pgfsetadditionalshadetransform{ \pgftransformshift{\pgfpoint{89.1 bp } { -128.7 bp }  }  \pgftransformscale{1.32 }  }}}
			\pgfdeclareradialshading{_e1kq4nvs0}{\pgfpoint{-72bp}{104bp}}{rgb(0bp)=(1,1,1);
				rgb(0bp)=(1,1,1);
				rgb(25bp)=(0.29,0.56,0.89);
				rgb(400bp)=(0.29,0.56,0.89)}
			
			
			\tikzset {_m0uizghyu/.code = {\pgfsetadditionalshadetransform{ \pgftransformshift{\pgfpoint{89.1 bp } { -128.7 bp }  }  \pgftransformscale{1.32 }  }}}
			\pgfdeclareradialshading{_ukwb1jjne}{\pgfpoint{-72bp}{104bp}}{rgb(0bp)=(1,1,1);
				rgb(0bp)=(1,1,1);
				rgb(25bp)=(0.29,0.56,0.89);
				rgb(400bp)=(0.29,0.56,0.89)}
			
			
			\tikzset {_iwq1safzj/.code = {\pgfsetadditionalshadetransform{ \pgftransformshift{\pgfpoint{89.1 bp } { -128.7 bp }  }  \pgftransformscale{1.32 }  }}}
			\pgfdeclareradialshading{_ylkt88q2j}{\pgfpoint{-72bp}{104bp}}{rgb(0bp)=(1,1,1);
				rgb(0bp)=(1,1,1);
				rgb(25bp)=(0.29,0.56,0.89);
				rgb(400bp)=(0.29,0.56,0.89)}
			
			
			\tikzset {_r0zikvfxd/.code = {\pgfsetadditionalshadetransform{ \pgftransformshift{\pgfpoint{89.1 bp } { -128.7 bp }  }  \pgftransformscale{1.32 }  }}}
			\pgfdeclareradialshading{_4z1rc5g2v}{\pgfpoint{-72bp}{104bp}}{rgb(0bp)=(1,1,1);
				rgb(0bp)=(1,1,1);
				rgb(25bp)=(0.29,0.56,0.89);
				rgb(400bp)=(0.29,0.56,0.89)}
			\tikzset{every picture/.style={line width=0.75pt}} 
			
			\begin{tikzpicture}[x=0.75pt,y=0.75pt,yscale=-1,xscale=1]
				
				\draw [shading=_hr48qomrb,_m5t92ypuz]   (65.41,68.52) -- (65.91,107.66) ;
				\draw [shading=_84c4b7as9,_qwkbg84js]   (32.48,34.38) -- (94.36,34.38) ;
				\draw [shading=_5iwkzmjsk,_5d8z04i89]   (199.16,67.52) -- (228.6,33.38) ;
				\draw [shading=_gca0ma2y6,_s0h72ppzf]   (199.16,67.52) -- (166.72,33.38) ;
				\draw [shading=_jouysak7d,_wm67zkktr]   (228.6,33.38) -- (199.66,106.66) ;
				\draw [shading=_2mexbeuiv,_2umrx5dr6]   (166.72,33.38) -- (199.66,106.66) ;
				\draw  [draw opacity=0][shading=_jvsh7ppcb,_0ibysugxg] (21.78,34.38) .. controls (21.78,28.47) and (26.57,23.68) .. (32.48,23.68) .. controls (38.39,23.68) and (43.18,28.47) .. (43.18,34.38) .. controls (43.18,40.29) and (38.39,45.08) .. (32.48,45.08) .. controls (26.57,45.08) and (21.78,40.29) .. (21.78,34.38) -- cycle ;
				\draw  [draw opacity=0][shading=_thuovt89s,_x2y3pyfpa] (83.66,34.38) .. controls (83.66,28.47) and (88.45,23.68) .. (94.36,23.68) .. controls (100.27,23.68) and (105.06,28.47) .. (105.06,34.38) .. controls (105.06,40.29) and (100.27,45.08) .. (94.36,45.08) .. controls (88.45,45.08) and (83.66,40.29) .. (83.66,34.38) -- cycle ;
				\draw  [draw opacity=0][shading=_qh9idu2gg,_prd9o6hbr] (54.71,68.52) .. controls (54.71,62.61) and (59.5,57.82) .. (65.41,57.82) .. controls (71.32,57.82) and (76.11,62.61) .. (76.11,68.52) .. controls (76.11,74.42) and (71.32,79.21) .. (65.41,79.21) .. controls (59.5,79.21) and (54.71,74.42) .. (54.71,68.52) -- cycle ;
				\draw  [draw opacity=0][shading=_a2he04xrr,_832d3j4ow] (55.21,107.66) .. controls (55.21,101.75) and (60,96.96) .. (65.91,96.96) .. controls (71.82,96.96) and (76.61,101.75) .. (76.61,107.66) .. controls (76.61,113.57) and (71.82,118.36) .. (65.91,118.36) .. controls (60,118.36) and (55.21,113.57) .. (55.21,107.66) -- cycle ;
				\draw  [draw opacity=0][shading=_e1kq4nvs0,_pr3hpz9l5] (156.02,33.38) .. controls (156.02,27.47) and (160.81,22.68) .. (166.72,22.68) .. controls (172.63,22.68) and (177.42,27.47) .. (177.42,33.38) .. controls (177.42,39.29) and (172.63,44.08) .. (166.72,44.08) .. controls (160.81,44.08) and (156.02,39.29) .. (156.02,33.38) -- cycle ;
				\draw  [draw opacity=0][shading=_ukwb1jjne,_m0uizghyu] (217.9,33.38) .. controls (217.9,27.47) and (222.69,22.68) .. (228.6,22.68) .. controls (234.51,22.68) and (239.3,27.47) .. (239.3,33.38) .. controls (239.3,39.29) and (234.51,44.08) .. (228.6,44.08) .. controls (222.69,44.08) and (217.9,39.29) .. (217.9,33.38) -- cycle ;
				\draw  [draw opacity=0][shading=_ylkt88q2j,_iwq1safzj] (188.96,106.66) .. controls (188.96,100.75) and (193.75,95.96) .. (199.66,95.96) .. controls (205.57,95.96) and (210.36,100.75) .. (210.36,106.66) .. controls (210.36,112.57) and (205.57,117.36) .. (199.66,117.36) .. controls (193.75,117.36) and (188.96,112.57) .. (188.96,106.66) -- cycle ;
				\draw  [draw opacity=0][shading=_4z1rc5g2v,_r0zikvfxd] (188.46,67.52) .. controls (188.46,61.61) and (193.25,56.82) .. (199.16,56.82) .. controls (205.07,56.82) and (209.86,61.61) .. (209.86,67.52) .. controls (209.86,73.42) and (205.07,78.21) .. (199.16,78.21) .. controls (193.25,78.21) and (188.46,73.42) .. (188.46,67.52) -- cycle ;
				
				\draw (63.42,137) node  [font=\footnotesize]  {$G$};
				\draw (32.48,34.38) node   [align=left] {1};
				\draw (94.36,34.38) node   [align=left] {2};
				\draw (65.41,68.52) node   [align=left] {3};
				\draw (65.91,107.66) node   [align=left] {4};
				\draw (199.16,136) node  [font=\footnotesize]  {$G^{\circ }$};
				\draw (166.72,33.38) node   [align=left] {1};
				\draw (199.16,67.52) node   [align=left] {3};
				\draw (199.66,106.66) node   [align=left] {4};
				\draw (228.6,33.38) node   [align=left] {2};

			\end{tikzpicture}
			\caption{$G$ consists of two disconnected components both of the same order $2$ and the complement $G^{\circ}$ is connected. The Laplacian matrix $L_{G^{\circ}}$ has a simple and largest eigenvalue. However there are no eigenvectors giving the $\rho$-versatility condition for $\rho > 1$, except for multiples of $\bf{1}$.\label{f4}}
		\end{minipage}
		\hfill
		\begin{minipage}[t]{0.45\textwidth}
			\centering
			
			\tikzset {_p9y0517mo/.code = {\pgfsetadditionalshadetransform{ \pgftransformshift{\pgfpoint{89.1 bp } { -128.7 bp }  }  \pgftransformscale{1.32 }  }}}
			\pgfdeclareradialshading{_z1rji44xf}{\pgfpoint{-72bp}{104bp}}{rgb(0bp)=(1,1,1);
				rgb(0bp)=(1,1,1);
				rgb(25bp)=(0.29,0.56,0.89);
				rgb(400bp)=(0.29,0.56,0.89)}
			
			
			\tikzset {_qz7p1iz4t/.code = {\pgfsetadditionalshadetransform{ \pgftransformshift{\pgfpoint{89.1 bp } { -128.7 bp }  }  \pgftransformscale{1.32 }  }}}
			\pgfdeclareradialshading{_6f3p3cqwj}{\pgfpoint{-72bp}{104bp}}{rgb(0bp)=(1,1,1);
				rgb(0bp)=(1,1,1);
				rgb(25bp)=(0.29,0.56,0.89);
				rgb(400bp)=(0.29,0.56,0.89)}
			
			
			\tikzset {_u4kxf007k/.code = {\pgfsetadditionalshadetransform{ \pgftransformshift{\pgfpoint{89.1 bp } { -128.7 bp }  }  \pgftransformscale{1.32 }  }}}
			\pgfdeclareradialshading{_61cfgojpz}{\pgfpoint{-72bp}{104bp}}{rgb(0bp)=(1,1,1);
				rgb(0bp)=(1,1,1);
				rgb(25bp)=(0.29,0.56,0.89);
				rgb(400bp)=(0.29,0.56,0.89)}
			
			
			\tikzset {_5e3rar0pp/.code = {\pgfsetadditionalshadetransform{ \pgftransformshift{\pgfpoint{89.1 bp } { -128.7 bp }  }  \pgftransformscale{1.32 }  }}}
			\pgfdeclareradialshading{_7xqdwl5d8}{\pgfpoint{-72bp}{104bp}}{rgb(0bp)=(1,1,1);
				rgb(0bp)=(1,1,1);
				rgb(25bp)=(0.29,0.56,0.89);
				rgb(400bp)=(0.29,0.56,0.89)}
			
			
			\tikzset {_lnsnybvf3/.code = {\pgfsetadditionalshadetransform{ \pgftransformshift{\pgfpoint{89.1 bp } { -128.7 bp }  }  \pgftransformscale{1.32 }  }}}
			\pgfdeclareradialshading{_l6w8sclsk}{\pgfpoint{-72bp}{104bp}}{rgb(0bp)=(1,1,1);
				rgb(0bp)=(1,1,1);
				rgb(25bp)=(0.29,0.56,0.89);
				rgb(400bp)=(0.29,0.56,0.89)}
			
			
			\tikzset {_syl0iwqw8/.code = {\pgfsetadditionalshadetransform{ \pgftransformshift{\pgfpoint{89.1 bp } { -128.7 bp }  }  \pgftransformscale{1.32 }  }}}
			\pgfdeclareradialshading{_9e1xyo3wu}{\pgfpoint{-72bp}{104bp}}{rgb(0bp)=(1,1,1);
				rgb(0bp)=(1,1,1);
				rgb(25bp)=(0.29,0.56,0.89);
				rgb(400bp)=(0.29,0.56,0.89)}
			
			
			\tikzset {_9xxwdcoo2/.code = {\pgfsetadditionalshadetransform{ \pgftransformshift{\pgfpoint{89.1 bp } { -128.7 bp }  }  \pgftransformscale{1.32 }  }}}
			\pgfdeclareradialshading{_ksruf6qu9}{\pgfpoint{-72bp}{104bp}}{rgb(0bp)=(1,1,1);
				rgb(0bp)=(1,1,1);
				rgb(25bp)=(0.29,0.56,0.89);
				rgb(400bp)=(0.29,0.56,0.89)}
			
			
			\tikzset {_4hb9jxwco/.code = {\pgfsetadditionalshadetransform{ \pgftransformshift{\pgfpoint{89.1 bp } { -128.7 bp }  }  \pgftransformscale{1.32 }  }}}
			\pgfdeclareradialshading{_bngviui0i}{\pgfpoint{-72bp}{104bp}}{rgb(0bp)=(1,1,1);
				rgb(0bp)=(1,1,1);
				rgb(25bp)=(0.29,0.56,0.89);
				rgb(400bp)=(0.29,0.56,0.89)}
			
			
			\tikzset {_txgho5xgg/.code = {\pgfsetadditionalshadetransform{ \pgftransformshift{\pgfpoint{89.1 bp } { -128.7 bp }  }  \pgftransformscale{1.32 }  }}}
			\pgfdeclareradialshading{_8br9403ul}{\pgfpoint{-72bp}{104bp}}{rgb(0bp)=(1,1,1);
				rgb(0bp)=(1,1,1);
				rgb(25bp)=(0.29,0.56,0.89);
				rgb(400bp)=(0.29,0.56,0.89)}
			
			
			\tikzset {_26mq42a3i/.code = {\pgfsetadditionalshadetransform{ \pgftransformshift{\pgfpoint{89.1 bp } { -128.7 bp }  }  \pgftransformscale{1.32 }  }}}
			\pgfdeclareradialshading{_npyol1ibr}{\pgfpoint{-72bp}{104bp}}{rgb(0bp)=(1,1,1);
				rgb(0bp)=(1,1,1);
				rgb(25bp)=(0.29,0.56,0.89);
				rgb(400bp)=(0.29,0.56,0.89)}
			
			
			\tikzset {_1gta85xns/.code = {\pgfsetadditionalshadetransform{ \pgftransformshift{\pgfpoint{89.1 bp } { -128.7 bp }  }  \pgftransformscale{1.32 }  }}}
			\pgfdeclareradialshading{_pj2iw26ks}{\pgfpoint{-72bp}{104bp}}{rgb(0bp)=(1,1,1);
				rgb(0bp)=(1,1,1);
				rgb(25bp)=(0.29,0.56,0.89);
				rgb(400bp)=(0.29,0.56,0.89)}
			
			
			\tikzset {_96m09gq79/.code = {\pgfsetadditionalshadetransform{ \pgftransformshift{\pgfpoint{89.1 bp } { -128.7 bp }  }  \pgftransformscale{1.32 }  }}}
			\pgfdeclareradialshading{_celg4xd1g}{\pgfpoint{-72bp}{104bp}}{rgb(0bp)=(1,1,1);
				rgb(0bp)=(1,1,1);
				rgb(25bp)=(0.29,0.56,0.89);
				rgb(400bp)=(0.29,0.56,0.89)}
			
			
			\tikzset {_36e20q54r/.code = {\pgfsetadditionalshadetransform{ \pgftransformshift{\pgfpoint{89.1 bp } { -128.7 bp }  }  \pgftransformscale{1.32 }  }}}
			\pgfdeclareradialshading{_5w5rsupxm}{\pgfpoint{-72bp}{104bp}}{rgb(0bp)=(1,1,1);
				rgb(0bp)=(1,1,1);
				rgb(25bp)=(0.29,0.56,0.89);
				rgb(400bp)=(0.29,0.56,0.89)}
			
			
			\tikzset {_k5ie8k293/.code = {\pgfsetadditionalshadetransform{ \pgftransformshift{\pgfpoint{89.1 bp } { -128.7 bp }  }  \pgftransformscale{1.32 }  }}}
			\pgfdeclareradialshading{_dywvvlbhq}{\pgfpoint{-72bp}{104bp}}{rgb(0bp)=(1,1,1);
				rgb(0bp)=(1,1,1);
				rgb(25bp)=(0.29,0.56,0.89);
				rgb(400bp)=(0.29,0.56,0.89)}
			
			
			\tikzset {_mzdlgopv2/.code = {\pgfsetadditionalshadetransform{ \pgftransformshift{\pgfpoint{89.1 bp } { -128.7 bp }  }  \pgftransformscale{1.32 }  }}}
			\pgfdeclareradialshading{_wkbx26jbj}{\pgfpoint{-72bp}{104bp}}{rgb(0bp)=(1,1,1);
				rgb(0bp)=(1,1,1);
				rgb(25bp)=(0.29,0.56,0.89);
				rgb(400bp)=(0.29,0.56,0.89)}
			
			
			\tikzset {_f7im5ijnp/.code = {\pgfsetadditionalshadetransform{ \pgftransformshift{\pgfpoint{89.1 bp } { -128.7 bp }  }  \pgftransformscale{1.32 }  }}}
			\pgfdeclareradialshading{_w7qq44rhu}{\pgfpoint{-72bp}{104bp}}{rgb(0bp)=(1,1,1);
				rgb(0bp)=(1,1,1);
				rgb(25bp)=(0.29,0.56,0.89);
				rgb(400bp)=(0.29,0.56,0.89)}
			
			
			\tikzset {_7gc0s0r3l/.code = {\pgfsetadditionalshadetransform{ \pgftransformshift{\pgfpoint{89.1 bp } { -128.7 bp }  }  \pgftransformscale{1.32 }  }}}
			\pgfdeclareradialshading{_8l9b1c3rw}{\pgfpoint{-72bp}{104bp}}{rgb(0bp)=(1,1,1);
				rgb(0bp)=(1,1,1);
				rgb(25bp)=(0.29,0.56,0.89);
				rgb(400bp)=(0.29,0.56,0.89)}
			
			
			\tikzset {_swcya0164/.code = {\pgfsetadditionalshadetransform{ \pgftransformshift{\pgfpoint{89.1 bp } { -128.7 bp }  }  \pgftransformscale{1.32 }  }}}
			\pgfdeclareradialshading{_m0lhn941h}{\pgfpoint{-72bp}{104bp}}{rgb(0bp)=(1,1,1);
				rgb(0bp)=(1,1,1);
				rgb(25bp)=(0.29,0.56,0.89);
				rgb(400bp)=(0.29,0.56,0.89)}
			
			
			\tikzset {_d6mi6l8mn/.code = {\pgfsetadditionalshadetransform{ \pgftransformshift{\pgfpoint{89.1 bp } { -128.7 bp }  }  \pgftransformscale{1.32 }  }}}
			\pgfdeclareradialshading{_lxdvc5nxz}{\pgfpoint{-72bp}{104bp}}{rgb(0bp)=(1,1,1);
				rgb(0bp)=(1,1,1);
				rgb(25bp)=(0.29,0.56,0.89);
				rgb(400bp)=(0.29,0.56,0.89)}
			
			
			\tikzset {_vsmajtrvh/.code = {\pgfsetadditionalshadetransform{ \pgftransformshift{\pgfpoint{89.1 bp } { -128.7 bp }  }  \pgftransformscale{1.32 }  }}}
			\pgfdeclareradialshading{_x3244a6q2}{\pgfpoint{-72bp}{104bp}}{rgb(0bp)=(1,1,1);
				rgb(0bp)=(1,1,1);
				rgb(25bp)=(0.29,0.56,0.89);
				rgb(400bp)=(0.29,0.56,0.89)}
			
			
			\tikzset {_r6pv6da78/.code = {\pgfsetadditionalshadetransform{ \pgftransformshift{\pgfpoint{89.1 bp } { -128.7 bp }  }  \pgftransformscale{1.32 }  }}}
			\pgfdeclareradialshading{_a1f7l4u7v}{\pgfpoint{-72bp}{104bp}}{rgb(0bp)=(1,1,1);
				rgb(0bp)=(1,1,1);
				rgb(25bp)=(0.29,0.56,0.89);
				rgb(400bp)=(0.29,0.56,0.89)}
			
			
			\tikzset {_5quavshel/.code = {\pgfsetadditionalshadetransform{ \pgftransformshift{\pgfpoint{89.1 bp } { -128.7 bp }  }  \pgftransformscale{1.32 }  }}}
			\pgfdeclareradialshading{_x1vtzop57}{\pgfpoint{-72bp}{104bp}}{rgb(0bp)=(1,1,1);
				rgb(0bp)=(1,1,1);
				rgb(25bp)=(0.29,0.56,0.89);
				rgb(400bp)=(0.29,0.56,0.89)}
			
			
			\tikzset {_ziwpnzx8m/.code = {\pgfsetadditionalshadetransform{ \pgftransformshift{\pgfpoint{89.1 bp } { -128.7 bp }  }  \pgftransformscale{1.32 }  }}}
			\pgfdeclareradialshading{_vaisfccag}{\pgfpoint{-72bp}{104bp}}{rgb(0bp)=(1,1,1);
				rgb(0bp)=(1,1,1);
				rgb(25bp)=(0.29,0.56,0.89);
				rgb(400bp)=(0.29,0.56,0.89)}
			
			
			\tikzset {_dlpuhisvl/.code = {\pgfsetadditionalshadetransform{ \pgftransformshift{\pgfpoint{89.1 bp } { -128.7 bp }  }  \pgftransformscale{1.32 }  }}}
			\pgfdeclareradialshading{_azq210mar}{\pgfpoint{-72bp}{104bp}}{rgb(0bp)=(1,1,1);
				rgb(0bp)=(1,1,1);
				rgb(25bp)=(0.29,0.56,0.89);
				rgb(400bp)=(0.29,0.56,0.89)}
			
			
			\tikzset {_xfusv1unj/.code = {\pgfsetadditionalshadetransform{ \pgftransformshift{\pgfpoint{89.1 bp } { -128.7 bp }  }  \pgftransformscale{1.32 }  }}}
			\pgfdeclareradialshading{_pvpvcag25}{\pgfpoint{-72bp}{104bp}}{rgb(0bp)=(1,1,1);
				rgb(0bp)=(1,1,1);
				rgb(25bp)=(0.29,0.56,0.89);
				rgb(400bp)=(0.29,0.56,0.89)}
			
			
			\tikzset {_i5o292jc3/.code = {\pgfsetadditionalshadetransform{ \pgftransformshift{\pgfpoint{89.1 bp } { -128.7 bp }  }  \pgftransformscale{1.32 }  }}}
			\pgfdeclareradialshading{_mn9mf1zu8}{\pgfpoint{-72bp}{104bp}}{rgb(0bp)=(1,1,1);
				rgb(0bp)=(1,1,1);
				rgb(25bp)=(0.29,0.56,0.89);
				rgb(400bp)=(0.29,0.56,0.89)}
			
			
			\tikzset {_tfjir20m5/.code = {\pgfsetadditionalshadetransform{ \pgftransformshift{\pgfpoint{89.1 bp } { -128.7 bp }  }  \pgftransformscale{1.32 }  }}}
			\pgfdeclareradialshading{_bnog30q8t}{\pgfpoint{-72bp}{104bp}}{rgb(0bp)=(1,1,1);
				rgb(0bp)=(1,1,1);
				rgb(25bp)=(0.29,0.56,0.89);
				rgb(400bp)=(0.29,0.56,0.89)}
			\tikzset{every picture/.style={line width=0.75pt}} 
			
			\begin{tikzpicture}[x=0.75pt,y=0.75pt,yscale=-1,xscale=1]
				
				\draw  [draw opacity=0][shading=_z1rji44xf,_p9y0517mo] (110.16,58.77) .. controls (110.16,52.86) and (114.95,48.07) .. (120.85,48.07) .. controls (126.76,48.07) and (131.55,52.86) .. (131.55,58.77) .. controls (131.55,64.68) and (126.76,69.47) .. (120.85,69.47) .. controls (114.95,69.47) and (110.16,64.68) .. (110.16,58.77) -- cycle ;
				\draw [color={rgb, 255:red, 0; green, 0; blue, 0 }  ,draw opacity=1 ][shading=_6f3p3cqwj,_qz7p1iz4t]   (292.45,97.21) -- (205.95,85.08) ;
				\draw [color={rgb, 255:red, 0; green, 0; blue, 0 }  ,draw opacity=1 ][shading=_61cfgojpz,_u4kxf007k]   (292.45,97.21) -- (205.95,34.35) ;
				\draw [color={rgb, 255:red, 0; green, 0; blue, 0 }  ,draw opacity=1 ][shading=_7xqdwl5d8,_5e3rar0pp]   (292.45,58.63) -- (205.95,85.08) ;
				\draw [color={rgb, 255:red, 0; green, 0; blue, 0 }  ,draw opacity=1 ][shading=_l6w8sclsk,_lnsnybvf3]   (205.95,33.92) -- (292.45,58.63) ;
				\draw [color={rgb, 255:red, 0; green, 0; blue, 0 }  ,draw opacity=1 ][shading=_9e1xyo3wu,_syl0iwqw8]   (250.61,112.82) -- (292.45,58.63) ;
				\draw [color={rgb, 255:red, 0; green, 0; blue, 0 }  ,draw opacity=1 ][shading=_ksruf6qu9,_9xxwdcoo2]   (250.61,112.82) -- (205.95,33.92) ;
				\draw [color={rgb, 255:red, 0; green, 0; blue, 0 }  ,draw opacity=1 ][shading=_bngviui0i,_4hb9jxwco]   (292.45,97.21) -- (292.45,58.63) ;
				\draw [color={rgb, 255:red, 0; green, 0; blue, 0 }  ,draw opacity=1 ][shading=_8br9403ul,_txgho5xgg]   (292.45,58.63) -- (264.24,33.92) ;
				\draw [color={rgb, 255:red, 0; green, 0; blue, 0 }  ,draw opacity=1 ][shading=_npyol1ibr,_26mq42a3i]   (79.01,111.95) -- (120.85,96.35) ;
				\draw [color={rgb, 255:red, 0; green, 0; blue, 0 }  ,draw opacity=1 ][shading=_pj2iw26ks,_1gta85xns]   (34.35,84.21) -- (92.65,33.05) ;
				\draw [color={rgb, 255:red, 0; green, 0; blue, 0 }  ,draw opacity=1 ][shading=_celg4xd1g,_96m09gq79]   (34.35,33.05) -- (92.65,33.05) ;
				\draw [color={rgb, 255:red, 0; green, 0; blue, 0 }  ,draw opacity=1 ][shading=_5w5rsupxm,_36e20q54r]   (205.95,85.08) -- (250.61,112.82) ;
				\draw [color={rgb, 255:red, 0; green, 0; blue, 0 }  ,draw opacity=1 ][shading=_dywvvlbhq,_k5ie8k293]   (250.61,112.82) -- (264.24,33.92) ;
				\draw [color={rgb, 255:red, 0; green, 0; blue, 0 }  ,draw opacity=1 ][shading=_wkbx26jbj,_mzdlgopv2]   (292.45,97.21) -- (264.24,33.92) ;
				\draw [color={rgb, 255:red, 0; green, 0; blue, 0 }  ,draw opacity=1 ][shading=_w7qq44rhu,_f7im5ijnp]   (205.95,85.08) -- (205.95,33.92) ;
				\draw  [draw opacity=0][shading=_8l9b1c3rw,_7gc0s0r3l] (23.65,33.05) .. controls (23.65,27.14) and (28.44,22.35) .. (34.35,22.35) .. controls (40.26,22.35) and (45.05,27.14) .. (45.05,33.05) .. controls (45.05,38.96) and (40.26,43.75) .. (34.35,43.75) .. controls (28.44,43.75) and (23.65,38.96) .. (23.65,33.05) -- cycle ;
				\draw  [draw opacity=0][shading=_m0lhn941h,_swcya0164] (81.95,33.05) .. controls (81.95,27.14) and (86.74,22.35) .. (92.65,22.35) .. controls (98.56,22.35) and (103.35,27.14) .. (103.35,33.05) .. controls (103.35,38.96) and (98.56,43.75) .. (92.65,43.75) .. controls (86.74,43.75) and (81.95,38.96) .. (81.95,33.05) -- cycle ;
				\draw  [draw opacity=0][shading=_lxdvc5nxz,_d6mi6l8mn] (23.65,84.21) .. controls (23.65,78.3) and (28.44,73.51) .. (34.35,73.51) .. controls (40.26,73.51) and (45.05,78.3) .. (45.05,84.21) .. controls (45.05,90.12) and (40.26,94.91) .. (34.35,94.91) .. controls (28.44,94.91) and (23.65,90.12) .. (23.65,84.21) -- cycle ;
				\draw  [draw opacity=0][shading=_x3244a6q2,_vsmajtrvh] (68.31,111.95) .. controls (68.31,106.05) and (73.1,101.26) .. (79.01,101.26) .. controls (84.92,101.26) and (89.71,106.05) .. (89.71,111.95) .. controls (89.71,117.86) and (84.92,122.65) .. (79.01,122.65) .. controls (73.1,122.65) and (68.31,117.86) .. (68.31,111.95) -- cycle ;
				\draw  [draw opacity=0][shading=_a1f7l4u7v,_r6pv6da78] (110.16,96.35) .. controls (110.16,90.44) and (114.95,85.65) .. (120.85,85.65) .. controls (126.76,85.65) and (131.55,90.44) .. (131.55,96.35) .. controls (131.55,102.26) and (126.76,107.05) .. (120.85,107.05) .. controls (114.95,107.05) and (110.16,102.26) .. (110.16,96.35) -- cycle ;
				\draw  [draw opacity=0][shading=_x1vtzop57,_5quavshel] (195.25,33.92) .. controls (195.25,28.01) and (200.04,23.22) .. (205.95,23.22) .. controls (211.85,23.22) and (216.64,28.01) .. (216.64,33.92) .. controls (216.64,39.83) and (211.85,44.62) .. (205.95,44.62) .. controls (200.04,44.62) and (195.25,39.83) .. (195.25,33.92) -- cycle ;
				\draw  [draw opacity=0][shading=_vaisfccag,_ziwpnzx8m] (195.25,85.08) .. controls (195.25,79.17) and (200.04,74.38) .. (205.95,74.38) .. controls (211.85,74.38) and (216.64,79.17) .. (216.64,85.08) .. controls (216.64,90.98) and (211.85,95.77) .. (205.95,95.77) .. controls (200.04,95.77) and (195.25,90.98) .. (195.25,85.08) -- cycle ;
				\draw  [draw opacity=0][shading=_azq210mar,_dlpuhisvl] (253.54,33.92) .. controls (253.54,28.01) and (258.33,23.22) .. (264.24,23.22) .. controls (270.15,23.22) and (274.94,28.01) .. (274.94,33.92) .. controls (274.94,39.83) and (270.15,44.62) .. (264.24,44.62) .. controls (258.33,44.62) and (253.54,39.83) .. (253.54,33.92) -- cycle ;
				\draw  [draw opacity=0][shading=_pvpvcag25,_xfusv1unj] (239.91,112.82) .. controls (239.91,106.91) and (244.7,102.12) .. (250.61,102.12) .. controls (256.52,102.12) and (261.31,106.91) .. (261.31,112.82) .. controls (261.31,118.73) and (256.52,123.52) .. (250.61,123.52) .. controls (244.7,123.52) and (239.91,118.73) .. (239.91,112.82) -- cycle ;
				\draw  [draw opacity=0][shading=_mn9mf1zu8,_i5o292jc3] (281.75,97.21) .. controls (281.75,91.31) and (286.54,86.51) .. (292.45,86.51) .. controls (298.36,86.51) and (303.15,91.31) .. (303.15,97.21) .. controls (303.15,103.12) and (298.36,107.91) .. (292.45,107.91) .. controls (286.54,107.91) and (281.75,103.12) .. (281.75,97.21) -- cycle ;
				\draw  [draw opacity=0][shading=_bnog30q8t,_tfjir20m5] (281.75,58.63) .. controls (281.75,52.72) and (286.54,47.93) .. (292.45,47.93) .. controls (298.36,47.93) and (303.15,52.72) .. (303.15,58.63) .. controls (303.15,64.54) and (298.36,69.33) .. (292.45,69.33) .. controls (286.54,69.33) and (281.75,64.54) .. (281.75,58.63) -- cycle ;
				
				\draw (241.61,143.17) node  [font=\footnotesize]  {$G^{\circ }$};
				\draw (63.5,143.17) node  [font=\footnotesize]  {$G$};
				\draw (34.35,33.05) node   [align=left] {1};
				\draw (79.01,111.95) node   [align=left] {4};
				\draw (92.65,33.05) node   [align=left] {2};
				\draw (120.85,96.35) node   [align=left] {5};
				\draw (120.85,58.77) node   [align=left] {6};
				\draw (34.35,84.21) node   [align=left] {3};
				\draw (205.95,33.92) node   [align=left] {1};
				\draw (264.24,33.92) node   [align=left] {2};
				\draw (205.95,85.08) node   [align=left] {3};
				\draw (292.45,97.21) node   [align=left] {5};
				\draw (292.45,58.63) node   [align=left] {6};
				\draw (250.61,112.82) node   [align=left] {4};

			\end{tikzpicture}
			\caption{$G$ consists of three disconnected components, each of a different order. The complement $G^{\circ}$ is connected, but the largest eigenvalue of its Laplacian is non-simple. 			\label{f6}}
		\end{minipage}
		\hfill
	\end{figure}

	\subsection{Versatile graphs by means of the degree distribution}\label{stargraphs}
	
	We next investigate another pathway to $\rho$-versatility, namely by looking at the degree distribution of the nodes in the network. To this end, we will prove:
	
	\begin{prop}\label{degreedisgivesvers}
	Let $r<C$ be two positive integers and suppose $G = (V,E)$ is a graph consisting of one node of degree $C$ and $N$ nodes of degree at most $r$, where $N \geq 1$. If 
	\begin{equation}\label{equationfordis}
	    \frac{C+1}{r} > \sqrt[3]{N}+1\, ,
	\end{equation}
	then the largest eigenvalue of $L_{G}$ is simple and every corresponding eigenvector $v$ satisfies
		\begin{equation}
	    \sum_{i=0}^N\nu_i^{\ell} \not= 0\, , \text{ for all } \ell > 1\, .
	\end{equation}
	\end{prop}

The proof of Proposition \ref{degreedisgivesvers} uses the well-known result that for any graph $G$ with at least one edge, the largest eigenvalue $\mu$ of $L_{G}$ satisfies $\lambda \geq d+1$, with $d$ the largest degree of any node in $G$. 
In the setting of Proposition \ref{degreedisgivesvers} we therefore have $\lambda \geq C+1$.
	
\begin{dempropdegreedis}\label{propdegispvers}
Let $\lambda$ be the largest eigenvalue of $L_{G}$ and write $v \in \mathbb{R}^{N+1}$ for a corresponding eigenvector. 
By re-scaling $v$, we may assume that $|\nu_i| \leq 1$ for all $i \in \{0, \dots, n\}$. The condition that $(\lambda, v)$ is an eigenvalue-eigenvector pair for $L_{G}$ gives
	\begin{align}
	    \sum_{j=0}^N(L_{G})_{i,j}\nu_j = \lambda \nu_i\, \text{ for all } i \in \{0, \dots, n\}\, .
	\end{align}
	We therefore find
	\begin{align}\label{equalitydegreeesq}
	   \sum_{\substack{j=0 \\ j\neq i}}^N (L_{G})_{i,j}\nu_j = (\lambda -k_i) \nu_i\, ,
	\end{align}	
	where $k_i$ denotes the degree of node $i$. From our observation that $\lambda \geq C+1$ we see that $(\lambda -k_i)$ is always positive. For $i \neq 0$ we therefore get from Equation \eqref{equalitydegreeesq}
		\begin{align}
	   (C+1 - r) |\nu_i| \leq (\lambda -k_i) |\nu_i|  = \bigg|\sum_{\substack{j=0 \\ j\neq i}}^N (L_{G})_{i,j}\nu_j\bigg| \leq \sum_{\substack{j=0 \\ j\neq i}}^N |(L_{G})_{i,j}||\nu_j| \leq k_i \leq r\, .
	    \end{align}
	    We thus find
	    \begin{align}\label{boundsonviww}
	   |\nu_i| \leq \frac{r}{C+1 - r} = \frac{1}{(C+1)/r - 1} < \frac{1}{ \sqrt[3]{N}} \leq 1\, .
	    \end{align}

		Summarizing, we see that the condition $|\nu_i| \leq 1$ for all $i \in \{0, \dots, n\}$ yields $|\nu_i|<1$ for all $i \in \{1, \dots, n\}$. This is only possible if $\nu_0 \not= 0$, which therefore has to hold for any eigenvector $v$ of $\lambda$.

		Now suppose $v$ and $v'$ are two eigenvectors for $\lambda$. By the foregoing, there exists a nonzero scalar $s$ such that the vector $sv - v'$ has vanishing zeroth component. As nevertheless $L_{G}(sv - v') = \lambda(sv - v')$, we conclude that $sv - v' = 0$ and so $v' = sv$. This shows that the eigenvalue $\lambda$ is simple. 

		To prove the $\rho$-versatility claim, we re-scale the eigenvector $v$ such that $|\nu_i| \leq 1$ for all $i \in \{0, \dots, n\}$, with $\nu_j = 1$ for at least one $j$. By the forgoing, this means that necessarily  $\nu_0 = 1$, with the other $\nu_i$ satisfying Equation \eqref{boundsonviww}. We conclude that for all $\rho \geq 3$ we have 
	   \begin{align}\label{boundsonviwww2}
	  	\bigg|\sum_{i=0}^N \nu_i^{\rho}\bigg| = 	  \bigg|1 + \sum_{i=1}^N \nu_i^{\rho}\bigg| \geq 1 -  \sum_{i=1}^N |\nu_i|^{\rho} >  1 -  \frac{N}{ (\sqrt[3]{N})^{\rho}} = 1 - N^{1-\rho/3} \geq 1-N^0 = 0\, .
	    \end{align}
	    We therefore find $\sum_{i=0}^N \nu_i^{\rho} \not= 0$ for all $\rho \geq 3$. 
	    As we clearly have 
		$\sum_{i=0}^N \nu_i^{2} > 0$, the result follows.
	    \end{dempropdegreedis}

\begin{ex}
Examples of connected networks satisfying the conditions of Proposition \ref{degreedisgivesvers} can easily be constructed. Let $r,N>0$ be given numbers such that 
\begin{equation}\label{conditionwelldefinednodes}
    r \leq \frac{N}{\sqrt[3]{N}+1}\, .
\end{equation}
We first construct a graph $G'$ consisting of $N$ nodes, all of which have degree at most $r-1$. The graph $G$ is then obtained from $G'$ by adding a node $n_0$, together with $C \geq (\sqrt[3]{N}+1)r$ edges between $n_0$ and different nodes of $G'$. Note that condition \eqref{conditionwelldefinednodes} guarantees that $(\sqrt[3]{N}+1)r \leq N$, so that we are not demanding that $n_0$ is connected to more nodes than $G'$ contains. It follows that all nodes in $G$ apart from $n_0$ have degree at most $r$. Finally, the degree $C$ of $n_0$ satisfies 
\begin{align}
    C+1>C\geq (\sqrt[3]{N}+1)r \, ,
\end{align}
so that
\begin{align}
    \frac{C+1}{r}> \sqrt[3]{N}+1\, .
\end{align}
The graph $G$ is connected if an edge was added from $n_0$ to at least one node from every connected component of $G'$.
\end{ex}

The proof of Proposition \ref{eigenvaluegaptozero} uses a result about the effects of adding an edge to the graph on the spectrum of the Laplacian. Given a graph $G$, we denote by $G+e$ the graph obtained from $G$ by adding some edge $e$ that was not there before. If $G$ and $G+e$ have $M$ nodes, then we denote by $0 = \lambda^G_1 \leq \dots \leq \lambda^G_M$ the eigenvalues of $L_G$ and by $0 = \lambda^{G+e}_1 \leq \dots \leq \lambda^{G+e}_M$ the eigenvalues of $L_{G+e}$. It can then be shown that
	\begin{align}
		0 = \lambda^G_1 = \lambda^{G+e}_1 \leq \lambda^{G}_2 \leq \lambda^{G+e}_2 \leq \dots \leq \lambda^G_M \leq \lambda^{G+e}_M \, .
	\end{align}
	This result is sometimes referred to as an interlacing theorem for graphs, see \cite{mohar1991laplacian}. In the proof of Proposition \ref{eigenvaluegaptozero}, we are interested only in the inequality $\lambda^G_{M-1} \leq \lambda^{G+e}_{M-1}$ corresponding to the second-largest eigenvalues. Repeated use of this latter result gives $\lambda^G_{M-1} \leq \lambda^{G'}_{M-1}$, where $G'$ is obtained from $G$ by adding any number of edges.\\	
	\begin{demgap}
		Let us say that $G$ has $N+1$ nodes, and write $n_0$ for the unique node of degree $C$. We denote the eigenvalues of $L_G$ by $0 = \lambda^G_1 \leq \dots \leq  \lambda^G_{N+1}$, so that $\mu = \lambda^G_{N+1}$ and $\kappa = \lambda^G_{N}$. If we have $N=1$ then $\kappa = 0$, so that there is nothing left to prove. Hence, we assume from here on out that $N>1$. Just as in the proof of Proposition \ref{degreedisgivesvers}, we have
		\begin{equation}\label{partofmanycomplements1}
			\mu =  \lambda_{N+1} \geq C+1\, .
		\end{equation}
		Next, let  $G'$ denote the graph obtained from $G$ by adding edges between $n_0$ and other nodes until $n_0$ is connected to every other node. 
		By the observation above we have $\kappa =  \lambda^G_N \leq \lambda^{G'}_N$, where the eigenvalues of $L_{G'}$ are given by $0 = \lambda^{G'}_1 \leq \dots \leq \lambda^{G'}_{N+1}$. 
		Let us consider the graph $G'$ instead. Because $n_0$ is connected to every other node, we see that the complement graph $G'^{\circ}$ consists of two components: $\{n_0\}$ and the remaining part $H$, where we do not claim  $H$ itself is connected.		
		We denote by $0 = \lambda^H_2 \leq \dots \lambda^H_{N+1}$ the eigenvalues of $L_H$, so that those of $L_{G'^{\circ}}$ are given by $$0 = \lambda^{G'^{\circ}}_1 =  \lambda^H_2 \leq \dots \lambda^H_{N+1}.$$ By the techniques used in the proof of Theorem \ref{versatilefromcomp}, we conclude that $\lambda^H_3 = N+1-\lambda^{G'}_{N}$. 		
		Next, consider the complement $H^{\circ}$ of $H$. Again by the techniques used in the proof of Theorem \ref{versatilefromcomp}, we see that an eigenvalue of $L_{H^{\circ}}$ is given by $N - \lambda^H_3 =   N-(N+1-\lambda^{G'}_{N}) =\lambda^{G'}_{N} - 1$. Moreover,  by construction of $H$ and $H^{\circ}$, we see that this latter graph can be obtained from $G'$ (or from $G$), by deleting $n_0$ and every edge connected to this node. In particular, we conclude that every node in $H^{\circ}$ has degree at most $r$. A straightforward application of the Gershgorin disk theorem \cite{gerschgorin31} now tells us that all eigenvalues of $L_{H^{\circ}}$ are bounded from above by $2r$. In particular, we find $\lambda^{G'}_{N} - 1 \leq 2r$ and so
		\begin{equation}\label{partofmanycomplements2}
			\kappa \leq \lambda^{G'}_{N} \leq 2r+1\, .
		\end{equation}
		Combining equations \eqref{partofmanycomplements1} and \eqref{partofmanycomplements2}, we see that
		\begin{equation}\label{partofmanycomplements3}
			\frac{\kappa}{\mu} \leq \frac{2r+1}{C+1} \leq \frac{3r}{C} \, ,
		\end{equation}
		from which the result follows.
\end{demgap}

To conclude our discussion on $\rho$-versatile graphs, we fix values $n, \rho \in \mathbb{N}$ with $n>2$ and define $\mathcal{S}_n^{\rho}$ as the set of all symmetric $(n \times n)$ matrices with a simple largest eigenvalue, whose corresponding eigenvector $\nu$ satisfies
\begin{equation}\label{conditionficedrhoo}
\sum_{i=1}^n \nu_i^{\ell} \not= 0 \text{ for all } \ell \in  \{2, \dots \rho+1\} \, .
\end{equation}
It follows that $\mathcal{S}_n^{\rho}$ is an open subset of the space of symmetric matrices. Heuristically speaking, if $G$ is a graph such that $L_G \in \mathcal{S}_n^{\rho}$, then we therefore expect $L_{G'} \in \mathcal{S}_n^{\rho}$ for any graph $G'$ obtained from $G$ by a small perturbation. 

In fact, as matrices generically have simple eigenvalues, and as Equation \eqref{conditionficedrhoo} is likewise valid for generic (eigen)vectors, we expect $L_G \in \mathcal{S}_n^{\rho}$ for "most" graphs $G$. Of course these statements will have to be made precise, which we do not attempt here. 

One common obstruction to $L_G \in \mathcal{S}_n^{\rho}$ seems to be symmetry in the graph $G$.  An explanation for this is that symmetry often forces eigenvalues with high multiplicity.  Moreover, if an eigenvalue $\lambda$ of $L_G$ is simple, then the span of a corresponding eigenvector $v$ forms a $1$-dimensional representation of the symmetry in question.  For any finite group symmetry, a one-dimensional real representation is either trivial, or generated by $v \mapsto -v$.  In the latter case, the graph symmetry contains a transformation $\alpha$ such that for any node $n$ of $G$, the corresponding coefficients $\nu_{n}$ and $\nu_{\alpha(n)}$ of $v$ are related by $\nu_{n} = -\nu_{\alpha(n)}$. This means that for any value $c \in \mathbb{R}$ there is an equal number of nodes $n$ such that $\nu_{n}=c$ as there are nodes $m$ such that $\nu_{m}=-c$. As a consequence, we then necessarily have
\begin{equation}
\sum_{i=1}^n \nu_i^{\ell} = 0 \text{ for all odd } \ell > 0\, .
\end{equation}
This is a common observation; imposing additional structure on a graph (such as for instance symmetry) induces high dimensional center subspaces and restrictions to the Taylor-coefficients of reduced vector fields in the associated dynamical systems. This generally leads to more elaborate bifurcation scenarios.

	\subsection{Hurwitzness}\label{InternalDynamics}
	
	In this subsection we give examples of Hurwitz matrices $A$ such that there exist $m>0$ mutually orthogonal vectors $x_{1},\dots,x_{m}$ satisfying
		\begin{equation}\label{posit2on33}
		    \langle x_{i},Ax_{i} \rangle>0 \text{ for all } i=1,\dots,m \, .
		\end{equation} 
	Note that any Hurwitz matrix $A$ has a negative trace, as this number equals the sum of its eigenvalues. It follows that Equation \eqref{posit2on33} can then only hold when $m<n$, where $n$ is the size of $A$. We start by looking at the case $n=2$.
	
	\begin{ex}
	       A general $2$ by $2$ matrix $A$ is of the form
	       	\begin{equation}
	    A=\left(\begin{array}{cc}
	         a & b\\
	         c & d 
	    \end{array}\right)\, ,
	\end{equation}
	with $a,b,c,d \in \mathbb{R}$. $A$ is Hurwitz if and only if $a+d < 0$ and $ad-bc > 0$. This can easily be arranged if in addition $a>0$, by first choosing $d<0$ such that $a+d < 0$, and then choosing $b$ and $c$ such that $ad-bc > 0$. We can construct examples of Hurwitz matrices $A$ such that Equation \eqref{posit2on33} holds with $m=1$ and $x_1 = (1,0)^{T}$. The set of all such matrices forms a non-empty open subset of the space of all $2$ by $2$ matrices. A similar observation of course holds when $d>0$.
	\end{ex}
	
		\begin{ex}
		Consider the $3$ by $3$ matrix 
			\begin{equation}\label{circulantexam}
	    A=\left(\begin{array}{ccc}
	         a+b+c & e & d\\
	          c & a+e & b+d\\
	           c & b+e & a+d\\
	    \end{array}\right)\, ,
	\end{equation}
	for $a,b,c,d,e \in \mathbb{R}$. Using the theory of network multipliers, it is shown in \cite{deville2021circulant} that the eigenvalues of $A$ are given by $a+b+c+d+e$, $a-b$ and $a+b$. It is therefore clear that for certain choices of $a$ through $e$ we can arrange for $A$ to be Hurwitz. Moreover, these eigenvalues do not change if we apply the transformation 
	\begin{align}
	&c \mapsto c-2\delta \\ \nonumber
	&d \mapsto d+\delta \\ \nonumber
	&e \mapsto e+\delta 
	\end{align}
	for any $\delta \in \mathbb{R}$, while keeping $a$ and $b$ the same. Hence, if $A$ as given by Equation \eqref{circulantexam} is Hurwitz, then so is the matrix 
	\begin{equation}\label{circulantexamrer}
	    A_{\delta}=\left(\begin{array}{ccc}
	         a+b+c-2\delta & e+\delta & d+\delta\\
	          c-2\delta & a+e+\delta & b+d+\delta\\
	           c-2\delta & b+e+\delta & a+d+\delta\\
	    \end{array}\right)\, 
	\end{equation}
for any $\delta \in \mathbb{R}$. Choosing $\delta$ large enough, we see that Equation \eqref{posit2on33} holds for $m=2$ and $x_1 = (0,1,0)^T$, $x_2 = (0,0, 1)^T$.
		\end{ex}

	Finally, we show: 
	\begin{prop}
	Let $X = \{x_1, \dots, x_{n-1}\}$ be a set of $n-1$ mutually orthogonal vectors in $\mathbb{R}^n$, where $n>1$. Denote by $\mathcal{H}_X^n$ the set of all $(n \times n)$ Hurwitz matrices $A$ such that 
		\begin{equation}\label{posit2on3333}
		    \langle x_{i},Ax_{i} \rangle>0 \text{ for all } i=1,\dots,n-1 \, .
		\end{equation} 
		Then, $\mathcal{H}_X^n$ forms a non-empty open subset of the space of all $(n \times n)$ matrices. 
	\end{prop}
	
\begin{dem}
As the set of all Hurwitz matrices is open, and because the same holds for the set of all matrices $A$ for which Equation \eqref{posit2on3333} holds, we see that $\mathcal{H}_X^n$ is likewise open. It remains to show that $\mathcal{H}_X^n$ is non-empty. 

We will first show this when $X$ consists of the first $n-1$ standard vectors $e_1 = (1, 0, \dots, 0)^T$, $e_2 = (0, 1, \dots, 0)^T$ and so forth, up to $e_{n-1} = (0, \dots, 1, 0)^T$. Given numbers $b_1, \dots, b_n \in \mathbb{R}$, we define the $(n \times n)$ matrix
		\begin{equation}
			\nonumber A_{b_1, \dots, b_n}=\left(\begin{array}{cccc}
				b_1&b_2& \dots &b_n\\
				b_1&b_2& \dots &b_n\\
				\vdots&\vdots&\vdots&\vdots\\
				b_1&b_2& \dots &b_n
			\end{array}\right)\, .
		\end{equation}
    As $A_{b_1, \dots, b_n}$ has rank $1$, we see that it has an $(n-1)$-dimensional kernel. The remaining eigenvalue is given by $b_1+ \dots +b_n$ with eigenvector $(1,\dots, 1)^T$. Let us choose $b_1, \dots, b_{n-1} > 0$ and $b_n < -(b_1+ \dots +b_{n-1})$. We also choose $a \in \mathbb{R}$ such that 
    $$0<a<\min(b_1, \dots, b_{n-1})\, .$$ As a result, we see that the matrix
    \begin{equation}
			\nonumber A_{b_1, \dots, b_n}-a\id_n=\left(\begin{array}{cccc}
				b_1-a&b_2& \dots &b_n\\
				b_1&b_2-a& \dots &b_n\\
				\vdots&\vdots&\vdots&\vdots\\
				b_1&b_2& \dots &b_n-a
			\end{array}\right)\, 
		\end{equation}
		has eigenvalues $b_1+ \dots +b_{n} - a < 0$ and $-a<0$. 
		Moreover, because $b_i-a > 0$ for all $i \in \{1, \dots, n-1\}$, we conclude that $A_{b_1, \dots, b_n}-a\id_n \in \mathcal{H}_E^n$ where $E = \{e_1, \dots, e_{n-1}\}$. 
		
To show that $\mathcal{H}_X^n$ is non-empty for general $X$, we pick $X = \{x_1, \dots, x_{n-1}\}$ and extend it to an orthogonal basis $\{x_1, \dots, x_{n-1}, x_n\}$. Let $U$ be the matrix such that $Ue_i = x_i$ for all $i \in \{1, \dots, n\}$. It follows that $U^TU$ equals a diagonal matrix $D$ with positive diagonal entries given by $\langle x_i, x_i \rangle = \|x_i\|^2$. In particular, we have $U^{T} = DU^{-1}$. Now suppose we pick an element $A \in \mathcal{H}_E^n$. It follows that $UAU^{-1}$ is Hurwitz as well. Moreover, for all $i \in \{1, \dots, n-1\}$ we find
		\begin{align}
		    \langle x_i, UAU^{-1}x_i\rangle =  \langle U^Tx_i, AU^{-1}x_i\rangle =  \langle DU^{-1}x_i, AU^{-1}x_i\rangle =  \langle De_i, Ae_i\rangle =\|x_i\|^2\langle e_i, Ae_i\rangle > 0\, .
		\end{align}
This shows that $\mathcal{H}_X^n$ is likewise non-empty, which concludes the proof.
\end{dem}
	
Note that we have $A \in \mathcal{H}_X^n \implies cA \in \mathcal{H}_X^n$ for all $c \in \mathbb{R}_{>0}$. Taking the union over all the sets $\mathcal{H}_X^n$, we arrive at:
	
	\begin{cor}
	    Given $n>1$, the set of all $(n \times n)$ Hurwitz matrices $A$ for which some orthogonal vectors $x_1, \dots, x_{n-1}$ exist satisfying 
	    	\begin{equation}\label{posit2on3333ddd}
		    \langle x_{i},Ax_{i} \rangle>0 \text{ for all } i = 1, \dots, n-1\, ,
		\end{equation} 
		is open and non-empty.
	\end{cor}

 	\bibliographystyle{plain}
 	\bibliography{bibliography.bib}

\end{document}